\documentclass[12pt]{article}
\usepackage{amssymb,amsmath,amsfonts}
\usepackage{graphicx}
\usepackage[labelsep=period]{caption}
\usepackage{relsize}
\usepackage[sort]{cite}
\usepackage{esvect}
\usepackage{txfonts}
\usepackage{color}

\makeatletter
\renewcommand{\fnum@figure}{Fig. \thefigure}
\makeatother
\newtheorem{theorem}{Theorem}[section]
\newtheorem{definition}[theorem]{Definition}
\newtheorem{corollary}[theorem]{Corollary}
\newtheorem{proposition}[theorem]{Proposition}
\newtheorem{remark}[theorem]{Remark}
\newtheorem{lemma}[theorem]{Lemma}
\newtheorem{claim}[theorem]{Claim}

\newtheorem{assumption}[theorem]{Assumption}
\newtheorem{fnassumption}[theorem]{Finiteness Assumption}
\newtheorem{wplemma}[theorem]{Whitney Partition Lemma}
\newtheorem{ptchlemma}[theorem]{Patching Lemma}
\newtheorem{mainlemma}[theorem]{Main Lemma}
\newtheorem{indhypothesis}[theorem]{Inductive Hypothesis}
\newtheorem{hmla}[theorem]{Hypotheses of the Main Lemma for the Label $\Ac$}
\newtheorem{lassumption}[theorem]{Large $A$ Assumption}
\newtheorem{lassum}[theorem]{Large $A$ Assumption for Whitney Partitions}
\newtheorem{smepsassumption}[theorem]{Small $\ve$ Assumption}
\newtheorem{contheorem}[theorem]{(Conjectured) Theorem}
\newtheorem{hellytheorem}[theorem]{Helly's Theorem}
\newcommand {\Ac}      {{\mathcal A}}
\newcommand {\Bc}      {{\mathcal B}}

\newcommand {\Fc}      {{\mathcal F}}

\newcommand {\Kc}      {{\mathcal K}}

\newcommand {\Mc}      {{\mathcal M}}
\newcommand {\Nc}      {{\mathcal N}}
\newcommand {\R}       {{\mathbb R}}
\newcommand {\N}       {{\mathbb N}}

\newcommand {\tB}      {\widetilde{B}}

\newcommand {\tE}      {\widetilde{E}}
\newcommand {\tF}      {\widetilde{F}}

\newcommand {\tL}      {\widetilde{L}}

\newcommand {\tS}      {\widetilde{S}}
\newcommand {\tT}      {\widetilde{T}}

\newcommand {\tX}      {\widetilde{X}}

\newcommand {\tf}      {\tilde{f}}

\newcommand {\ta}      {\tilde{a}}

\newcommand {\tv}      {\tilde{v }}

\newcommand {\tz}      {\tilde{z}}
\newcommand {\tl}      {\tilde{\ell}}
\newcommand {\trh}     {\tilde{\rho}}
\newcommand {\tpsi}    {\widetilde{\psi}}

\newcommand {\hF}      {\hat{F}}
\newcommand {\hf}      {\hat{f}}
\newcommand {\hS}      {\hat{S}}

\newcommand {\ha}      {\hat{a}}
\newcommand {\hb}      {\hat{b}}
\newcommand {\hc}      {\hat{c}}
\newcommand {\hx}      {\hat{x}}
\newcommand {\hxi}     {\hat{\xi}}
\newcommand {\ve}      {\varepsilon}
\newcommand {\RN}      {\R^n}
\newcommand {\RM}      {\R^m}
\newcommand {\MS}      {\Mc}
\newcommand {\BY}      {B_{\BS}}
\newcommand {\BL}      {B_{\BS}}
\newcommand {\dm}      {\rho}
\newcommand {\dt}      {d}
\newcommand {\BS}      {Y}
\newcommand {\KM}      {\Kc_m(\BS)}
\newcommand {\KMY}     {\Kc_m(\BS)}
\newcommand {\KRM}     {\Kc(\RM)}
\newcommand {\SM}      {{\Mc'}}
\newcommand {\CNMY}    {\Conv_m(Y)}
\newcommand {\SX}      {S_{\hspace{-0.5mm}Y}}
\newcommand {\AM}      {\AFF_m(\BS)}
\newcommand {\FN}      {N(m,\BS)}
\newcommand {\CRF}     {G}
\newcommand {\REL}     {{\sc Rel}~}
\newcommand {\TRL}     {\text{\sc Rel}}
\newcommand {\RLV}     {\text{\sc Rlv}}
\newcommand {\ks}      {k^{\sharp}}
\newcommand {\kn}      {k^*}
\newcommand {\gnd}     {\gamma}
\newcommand {\gz}      {\gamma_0}
\newcommand {\gone}    {\gamma_1}
\newcommand {\GL}      {\Gamma_\ell}
\newcommand {\GLO}     {\Gamma_{\ell-1}}
\newcommand {\GH}      {\hat{\Gamma}}
\newcommand {\hAc}     {\hat{\Ac}}
\newcommand {\dl}      {\delta}
\newcommand {\Gm}      {\Gamma}
\newcommand {\hsg}     {\hat{\sigma}}
\newcommand {\hv}      {\hat{v}}
\newcommand {\hz}      {\hat{\zeta}}
\newcommand {\BXR}     {B(x_0,r_0)}
\newcommand {\GLA}     {\Gm_{\ell(\Ac)}}
\newcommand {\GLAO}    {\Gm_{\ell(\Ac)-1}}
\newcommand {\FF}      {F}
\newcommand {\ff}      {f}
\newcommand {\BO}      {\tB}
\newcommand {\tff}     {\tilde{\ff}}
\newcommand {\MPL}     {m+1}
\newcommand {\wtM}     {\widetilde{\Mc}}
\newcommand {\cng}     {c_{\hspace{-0.3mm}\mathsmaller{\Nc}}}
\newcommand {\CNG}     {D_{\hspace{-0.2mm}\mathsmaller{\Nc}}}
\newcommand {\CLS}     {C_{LS}}
\newcommand {\CWH}     {C_{Wh}}
\newcommand {\CETA}    {C_{\eta}}
\newcommand {\CAGR}    {C^\#}
\newcommand {\CLIP}    {C_{Lip}}
\newcommand {\DS}      {D^*}
\newcommand {\pl}      {10}
\newcommand {\emin}    {\ve^*}

\newcommand {\bsk}     {\bigskip}

\newcommand{\ltfrac}[2]{\mbox{\large$\frac{#1}{#2}$}}
\newcommand {\emp}     {\emptyset}
\newcommand {\vf}      {\varphi}
\newcommand {\je}      {\leftrightarrow}
\newcommand {\ip}[1]   {\langle{#1}\rangle}
\newcommand {\smed}    {\mathlarger{\sum}}

\newcommand {\capsm}   {\mathsmaller{\bigcap}}
\newcommand {\capbig}  {\mathlarger{\mathlarger{\cap}}}

\newcommand {\Lip}     {\operatorname{Lip}}
\newcommand {\diam}    {\operatorname{diam}}
\newcommand {\dist}    {\operatorname{dist}}

\newcommand {\Aff}     {\operatorname{Aff}}
\newcommand {\AFF}     {\operatorname{{\it Aff}}}
\newcommand {\Vect}    {\operatorname{Vect}}
\newcommand {\Conv}    {\operatorname{Conv}}
\newcommand {\dhf}     {\operatorname{d_H}}
\newcommand {\affspan} {\operatorname{{\it affhull\hspace*{0.5mm}}}}
\newcommand {\RELX}    {\operatorname{RELX}}
\newcommand {\cl}      {\operatorname{\,cl}}
\newcommand {\Sel}     {\operatorname{Selec}}
\newcommand {\ST}      {\operatorname{St}}
\newcommand {\DA}      {\operatorname{DA}}
\newcommand {\APT}     {\operatorname{APT}}
\newcommand {\VST}     {\vskip 1mm}
\newcommand {\VSU}     {\vskip 1mm}
\newcommand {\bx}      {\hfill$\Box$}

\newcommand {\BXD}     {\hspace{12mm}\Box}
\newcommand {\nn}      {\nonumber}
\newcommand {\rf}[1]    {(\ref{#1})}      
\newcommand {\reff}[1] {\ref{#1}}         
\newcommand{\lbl}[1]      {\label{#1}}       
\newcommand{\be}          {\begin{eqnarray}}
\newcommand{\bel}[1]      {\begin{eqnarray} \label{#1}}
\newcommand{\ee}           {\end{eqnarray}}
\newcommand {\SECT}[2] {\section*{\centerline{\normalsize
{\bf #1}}} \setcounter{section}{#2}
\setcounter{theorem}{0}\setcounter{equation}{0}}
\begin{document}
\parindent 1em
\parskip 0mm
\medskip
\centerline{\large{\bf\large SHARP FINITENESS PRINCIPLES FOR LIPSCHITZ SELECTIONS}}\vskip 9mm
\vspace*{5mm} \centerline{ \large{\sc Charles Fefferman~~~${\bf\cdot}$~~~Pavel Shvartsman}}
\vspace*{22 mm}
\renewcommand{\thefootnote}{ }
\footnotetext[1]{{\it\hspace{-6mm}Math Subject
Classification:} 46E35\smallskip
\par This research was supported by Grant No 2014055 from the United States-Israel Binational Science Foundation (BSF). The first author was also supported in part by NSF grant DMS-1265524 and AFOSR grant FA9550-12-1-0425.}

\begin{abstract} Let $(\MS,\dm)$ be a metric space and let $\BS$ be a Banach space. Given a positive integer $m$, let $F$ be a set-valued mapping from $\MS$ into the family of all compact convex subsets of $\BS$ of dimension at most $m$. In this paper we prove a finiteness principle for the existence of a Lipschitz selection of $F$ with the sharp value of the finiteness constant.

\end{abstract}
\renewcommand{\contentsname}{ }
\tableofcontents
\addtocontents{toc}{{\centerline{\sc{Contents}}}
\vspace*{10mm}\par}
\SECT{1. Introduction}{1}

\addtocontents{toc}{\hspace*{3.2mm} 1. Introduction\hfill \thepage\par\VST}
\indent\par We prove a finiteness theorem for {\it Lipschitz selection problems}, conjectured by Yu. Brudnyi and Shvartsman \cite{BS,S02} and established in special cases by Fefferman, Israel and Luli \cite{FIL1,FIL2} and Shvartsman \cite{S86,S92,S01,S02,S04}.
\par In its simplest setting, our problem is as follows. We are given a metric space $(\Mc,\rho)$ and a positive integer $m$. For each point $x\in\Mc$, we are given a nonempty compact convex set $F(x)\subset\RM$.
\par We want to find a Lipschitz map $f:\Mc\to\RM$ such that $f(x)\in F(x)$ for all $x\in\Mc$. Such an $f$ is called a {\it Lipschitz selection} of the set-valued map $F:\Mc\to\Kc(\RM)$, where $\KRM$ denotes the family of all nonempty compact convex subsets of $\RM$. If a Lipschitz selection $f$ exists, then we ask how small we can take its Lipschitz seminorm.
\par In this setting, our main result implies the following:
\begin{theorem}\lbl{RM-FP} Let $(\Mc,\rho)$ be a metric space, let $F:\Mc\to\KRM$, and let $\lambda$ be a positive real number. Suppose that for every $\SM\subset\Mc$ consisting of at most $2^m$ points, the restriction $F|_{\SM}$ of $F$ to $\SM$ has a Lipschitz selection $f_{\SM}$ with Lipschitz  seminorm at most $\lambda$.
\smallskip
\par Then $F$ has a Lipschitz selection with Lipschitz  seminorm at most $\gamma\lambda$. Here, $\gamma$ depends only on the dimension $m$.
\end{theorem}
\par Equivalently, we may suppose that $\Mc$ contains at least $2^m$ points and take $\Mc'$ to contain {\it exactly} $2^m$ points.
\par Lipschitz selection problems are closely related to
\medskip
\par\noindent {\bf Whitney's Extension Problem} (\cite{W1})~ {\it Fix $m,n\ge 1$, and let $f$ be a real-valued function defined on a given (arbitrary) closed set $E\subset\RN$. Decide whether $f$ extends to a function $F\in C^m(\RN)$ with a finite $C^m$-norm.
\par If such an extension $F$ exists, then how small can we take its $C^m$-norm?}
\smallskip
\par There is a finiteness principle for Whitney's Extension Problem, e.g. when $E\subset\RN$ is a large finite set. See Brudnyi-Shvartsman
\cite{BS,BS1,BS3,S82,S1,S02,S-Tr} and the later papers of Fefferman, Israel, Klartag and Luli
\cite{F2,F-J,F4,F6,F7,F8,F-Bl,FIL1,FIL2}, as well as A. and Yu. Brudnyi \cite{BB} for that finiteness principle and several related results.
\par The idea of Lipschitz selection first arose in connection with Whitney's extension problem, see \cite{BB,BS1,BS,BS3,S82,S84,S1}. In particular, a variant of a special case of Theorem \reff{RM-FP} was the main ingredient in the proof \cite{BS1,BS3,S82,S1,S02} of the finiteness principle for Whitney's Problem in the simplest non-trivial case, $m=2$.
The later papers \cite{F2,F-J,F4,F6,F7,F8,FIL1,FIL2} didn't explicitly mention Lipschitz selection, but they broadened Whitney's Problem by asking for functions $F\in C^m(\RN)$ that agree with $f$ on $E$ to a given accuracy.
\par Of course, a Lipschitz selection problem may also be regarded as a search for a smooth function that agrees approximately with data.
\medskip
\par Our main result is more general than Theorem \reff{RM-FP}. First of all, we allow $(\Mc,\rho)$ to be a {\it pseudometric space}, i.e., $\rho:\Mc\times\Mc\to [0,+\infty]$, $\rho(x,x)=0$, $\rho(x,y)=\rho(y,x)$, $\rho(x,y)\le \rho(x,z)+\rho(z,y)$ for all $x,y,z\in\Mc$.
Note that $\rho(x,y)=0$ may hold with $x\ne y$, and $\rho(x,y)$ may be $+\infty$.
\par Secondly, the convex sets $F(x)$ needn't sit inside $\RM$. Instead, we fix a Banach space $(Y,\|\cdot\|)$ and let $\KM$ denote the family of all nonempty compact convex subsets $K\subset Y$ of dimension at most $m$. (We say that a convex subset of $Y$ has dimension at most $m$ if it is contained in an affine subspace of $Y$ of dimension at most $m$.)
\par We write
\bel{NMY-1}
N(m,Y)=\min\{2^{m+1},2^{\dim Y}\} ~~~\text{if}~~Y~~\text{is finite-dimensional},
\ee
and
\bel{NMY-2}
N(m,Y)=2^{m+1}~~~\text{if}~~Y~~\text{is infinite-dimensional}.
\ee
\par We define the Lipschitz seminorm of a map $f:\Mc\to Y$ for a Banach space $Y$ and a pseudometric space $(\Mc,\rho)$ by setting
$$
\|f\|_{\Lip(\Mc,Y)}=\inf\{\,\lambda>0:\|f(x)-f(y)\|
\le\lambda\,\rho(x,y)~~~\text{for all} ~~~x,y\in\Mc\,\}.
$$
In particular, $\|f\|_{\Lip(\Mc,Y)}=+\infty$ if no such $\lambda$ exists.
\smallskip
\par We can now state our main result in full generality; Theorem \reff{RM-FP} will be a simple consequence.
\begin{theorem}\lbl{MAIN-FP} Fix $m\ge 1$. Let $(\Mc,\rho)$ be a pseudometric space, and let $F:\Mc\to\KM$ for a Banach space $Y$. Let $\lambda$ be a positive real number.
\par Suppose that for every $\SM\subset\Mc$ consisting of at most $N=N(m,Y)$ points, the restriction $F|_{\SM}$ of $F$ to $\SM$ has a Lipschitz selection $f_{\SM}$ with Lipschitz  seminorm $\|f_{\SM}\|_{\Lip(\SM,Y)}\le \lambda$.
\smallskip
\par Then $F$ has a Lipschitz selection $f$ with Lipschitz  seminorm $\|f\|_{\Lip(\Mc,Y)}\le \gamma\lambda$.
\par Here, $\gamma$ depends only on $m$.
\end{theorem}
\smallskip
\par The ``finiteness constants'' $2^m$ in Theorem \reff{RM-FP} and $N(m,Y)$ in Theorem \reff{MAIN-FP} are optimal; see \cite{S92} and\cite [Theorem 1.4]{S02}. We also refer the reader to the paper \cite[Section 8.1]{FS-2017}, which contains detailed proofs of this statement for $m=1,2$.
\par If the set $\Mc$ is {\it finite} in Theorem \reff{RM-FP} or Theorem \reff{MAIN-FP}, then we can omit the assumption that the convex sets $F(x)$ $(x\in\Mc)$ are compact. In this case, it is enough to assume that $F:\Mc\to\CNMY$, where
\bel{CMY}
\CNMY=\{\,\text{all nonempty convex subsets of}~~Y~~\text{of dimension at most}~~m\,\}.
\ee
See Theorem \reff{FINITE}.
\smallskip
\par For the case of the trivial distance function $\dm\equiv 0$, Theorems \reff{MAIN-FP} and
\reff{FINITE} agree with the classical Helly's Theorem \cite{DGK}, except that the optimal finiteness constant for $\dm\equiv 0$ is
\bel{DM-NMY}
n(m,Y)=\min\{m+2,\dim Y+1\}~~~~\text{in place of}~~~~N(m,Y)=\min\{2^{m+1},2^{\dim Y}\}.
\ee
Thus, our results may be regarded as a generalization of Helly's Theorem. However, we make extensive use of Helly's Theorem in our proofs.
\smallskip
\par Theorem \reff{MAIN-FP} and its variants were previously known in several special cases:
\smallskip
\par \textbullet~ $Y=\R^2$ \cite{S02};
\smallskip
\par \textbullet~ Each $F(x)$ ($x\in\Mc$) is an affine subspace of $Y$ of dimension at most $m$ \cite{S86,S92} ($Y=\RM$), \cite{S01} ($Y$ is a Hilbert space), \cite{S04} ($Y$ is a Banach space). Of course, all $F(x)$ are non-compact in this case;
\smallskip
\par \textbullet~ $(\MS,\rho)=(\RN,\|\cdot\|)$ and $Y=\RM$ with the constant $N$ and the constant $\gamma$ depending on $n$ as well as on $m$ \cite{FIL1}.
\medskip
\par Let us recount how we arrived at our proof of Theorem \reff{MAIN-FP}. P. Shvartsman (unpublished) had already reduced Theorem \reff{MAIN-FP} to the special case of a {\it metric tree} with nodes of bounded degree. We recall the relevant standard definitions.
\smallskip
\par Let $T=(X,E)$ be a finite (graph theoretic) tree, where $X$ denotes the set of nodes of $T$, and $E$ denotes the set of edges. The {\it degree} of a node $x\in X$ is the number of nodes $y$ to which $x$ is joined by an edge.
\par Suppose we assign a positive number $\Delta(e)$ to each edge $e\in E$. Then for $x,y\in X$ we can define their {\it distance} $d(x,y)$ to be the sum of $\Delta(e)$ over all the edges $e$ in the ``minimal path'' joining $x$ to $y$ as in Fig. 1.
\medskip
\begin{figure}[h!]
\hspace*{30mm}
\includegraphics[scale=0.35]{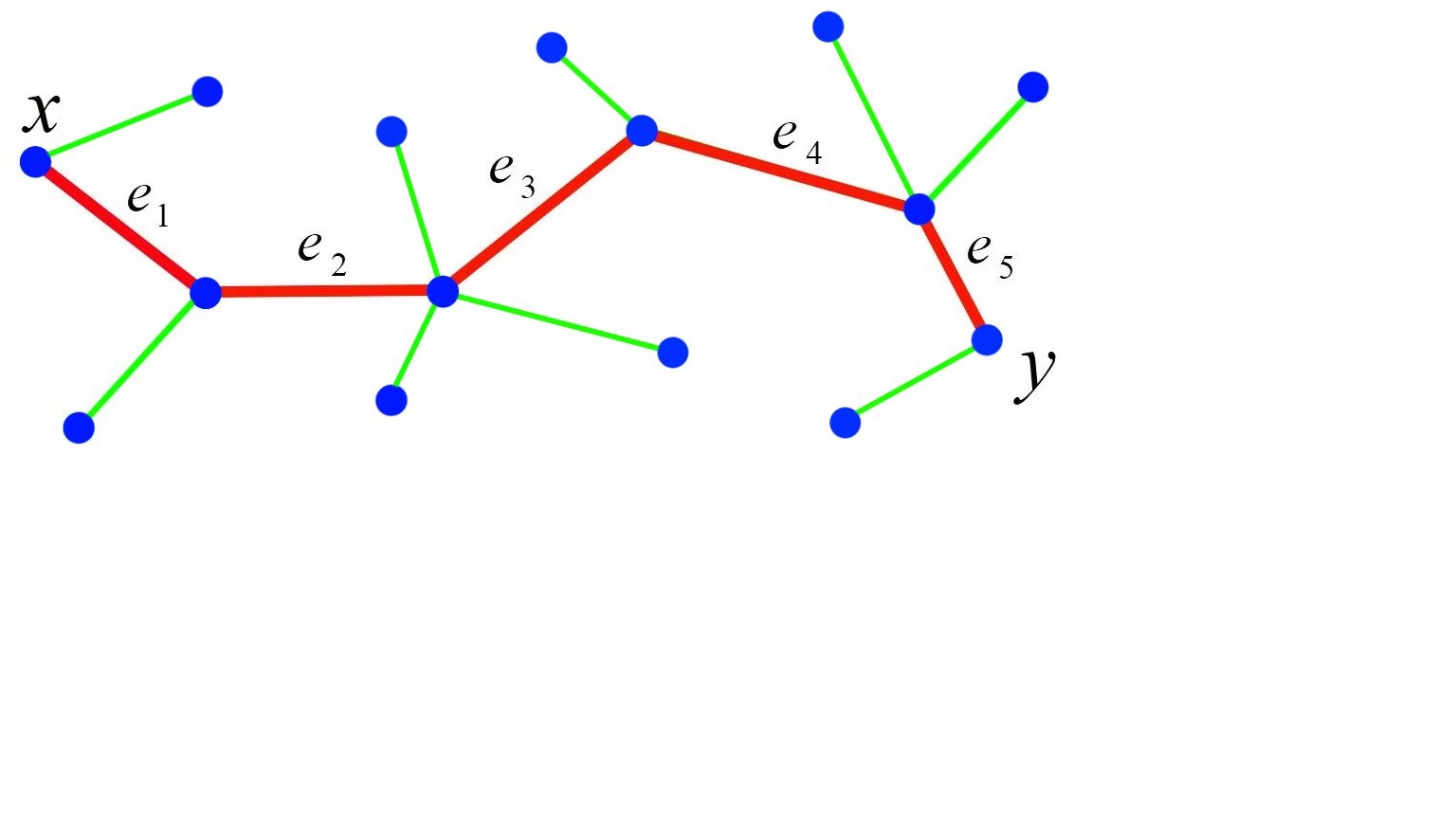}
\vspace*{-25mm}
\caption{A minimal path joining nodes $x$ and $y$ in a tree. In this case, $d(x,y)=\Delta(e_1)+\Delta(e_2)+...+\Delta(e_5).$}
\end{figure}
\par We call $d$ a {\it tree metric}; $(X,d)$ is a {\it metric tree}.
\par Shvartsman's unpublished previous work reduced Theorem \reff{MAIN-FP} to the following weakened form of a special case.
\begin{contheorem}\lbl{CONTH} Given $m\ge 1$, there exist $\ks, \gamma$ depending only on $m$, for which the following holds.
\par Let $(X,d)$ be a metric tree in which each node has degree at most $m+1$.
\par Let $F:X\to\KM$ for a Banach space $Y$, and let $\lambda$ be a positive real number. Suppose that for every subset $X'\subset X$ consisting of at most $\ks$ points, the restriction $F|_{X'}$ has a Lipschitz selection $f_{X'}$ with Lipschitz  seminorm $\|f_{X'}\|_{\Lip(X',Y)}\le \lambda$.
\smallskip
\par Then $F$ has a Lipschitz selection $f$ with Lipschitz  seminorm $\|f\|_{\Lip(X,Y)}\le \gamma\lambda$.
\end{contheorem}
\begin{remark} {\em  Note that here $X'$ needn't be a subtree of $X$. Thus, in Figure 1, perhaps $X'$ contains the nodes $x$ and $y$ but not the nodes that lie between them.
\par Note also that the optimal finiteness constant $N(m,Y)$ in Theorem \reff{MAIN-FP} has been replaced in Theorem \reff{CONTH} by a sufficiently large constant $\ks$ depending only on $m$.}
\end{remark}
\par On the other hand, the work of Fefferman, Israel  and Luli \cite{FIL1} on ``$C^m$ Selection'' implies a weakened version of Theorem \reff{RM-FP}, in which $(\Mc,\rho)$ is $\RN$ with its standard Euclidean metric; the sharp finiteness constant $2^m$ in Theorem \reff{RM-FP} is replaced by $\ks$ as in Theorem \reff{CONTH}; and the constant $\gamma$ is allowed to depend on $n$ as well as on $m$. See the web posting \cite{FIL-H}.
\par To prove Theorem \reff{MAIN-FP}, we set out to adapt the arguments in \cite{FIL1} from $\RN$ to the setting of a metric tree. If we succeeded, Theorem \reff{CONTH} would follow, thus proving Theorem \reff{MAIN-FP}.
\par This attempt seemed highly unlikely to succeed; the geometry of a metric tree is of course ra\-dically different from that of $\RN$. Nevertheless, we were able to adapt \cite{FIL1} and prove Theorem \reff{CONTH}, thanks to one crucial similarity between $\RN$ and metric trees - they have finite {\it Nagata dimension}. We recall the relevant definitions (see \cite{NGT,A,LS,BDHM}).
\begin{definition}\lbl{NG-C} {\em Let $(X,d)$ be a metric space. Let $D$ be a non-negative integer and let $c$ be a positive real constant. We say that $(X,d)$ satisfies {\it Nagata $(D,c)$} if for every real number $s>0$ there exists a covering $(X_i)_{i\in I}$ of $X$ by subsets $X_i$ of diameter at most $s$, such that no ball of radius $cs$ in $X$ meets more than $D+1$ of the $X_i$. We call $D,c$ the {\it Nagata constants} of $(X,d)$.
\par The least $D$ for which $(X,d)$ satisfies Nagata $(D,c)$ for some $c>0$ is the {\it Nagata dimension} (or Assouad-Nagata dimension) of $(X,d)$.}
\end{definition}
\par Note that any {\it finite} metric space has Nagata dimension $0$. The metric space $\RN$ has Nagata dimension $n$ (\cite{LS}). The space $\ell_\infty$ has infinite Nagata dimension because $\ell_\infty$ contains $\RN$ for each $n\in\N$.  Every {\it planar} connected graph whose nodes have finite degree has Nagata dimension at most $2^{10}-1$ (see Ostrovskii, Rosenthal \cite{OR} for the precise statement and the proof).
\par Moreover, every {\it metric tree} satisfies Nagata $(1,c)$ for an absolute constant $c$. (See \cite[Lemma 3.1 and Theorem 3.2]{LS}. For the reader's convenience, in Lemma \reff{MTR-ND} we prove that one can take $c=1/16$.) This allows us to carry over arguments in \cite{FIL1} from $\RN$ to an arbitrary metric tree.
\par More precisely, we prove the following result.
\begin{theorem}\lbl{NG-FP} Given $m\ge 1$ there exists $\ks$ depending only on $m$, for which the following holds.
\par Let $(X,d)$ be a finite metric space satisfying  Nagata $(D,c)$, and let $F:X\to\CNMY$ for a Banach space $Y$. Let $\lambda$ be a positive real number. Suppose that for every $X'\subset X$ consisting of at most $\ks$ points, the restriction $F|_{X'}$ has a Lipschitz selection $f_{X'}$ with Lipschitz  seminorm $\|f_{X'}\|_{\Lip(X',Y)}\le \lambda$.
\smallskip
\par Then $F$ has a Lipschitz selection $f$ with Lipschitz  seminorm $\|f\|_{\Lip(X,Y)}\le \gamma\lambda$, where $\gamma$ depends only on $m$ and on the Nagata constants $D,c$.
\end{theorem}
\par Recall that $\CNMY$ denotes the family of all nonempty convex subsets of $Y$ of dimension at most $m$ (see \rf{CMY}).
\par As an immediate corollary, we obtain a stronger form of Theorem \reff{CONTH} in which we drop the assumption that each node has degree at most $m+1$. See Corollary \reff{A-12}. So we have proven more than we need to establish Theorem \reff{MAIN-FP}. Because we needn't assume that the nodes of our metric tree have degree at most $m+1$, we can greatly simplify the earlier reduction of Theorem \reff{MAIN-FP} to the case of metric trees.
\smallskip
\par This paper is organized as follows.
\par In Section 2 we construct ``Whitney partitions of unity'' associated to a ``lengthscale'' $r(x)>0$ defined on a metric space of finite Nagata dimension. As in H. Whitney's classic paper \cite{W1}, such partitions are used to patch together functions defined in neighborhoods of varying sizes, while maintaining the smoothness of the functions being patched.
\par In Sections 3 and 4 we associate to a Lipschitz selection problem given by $F:\Mc\to\KM$ a family of convex sets $\GL(x)\in\KM$ parametrized by $x\in\Mc$ and $\ell\ge 0$. If for every subset $\SM\subset\Mc$ consisting of at most $\ks(\ell,m)$ points there exists a Lipschitz selection $f_{\SM}$ of $F|_{\SM}$ with Lipschitz seminorm $\|f_{\SM}\|_{\Lip(\SM,Y)}\le\lambda$, then $\GL(x)$ is nonempty. That is how we use the hypothesis of Theorem \reff{NG-FP}.
\par As in \cite{FIL1,FIL-H}, we use the $\GL(x)$ in Section 4 to prove Theorem \reff{NG-FP}. This is the most technically difficult part of our proof. The idea is to measure the difficulty of a Lipschitz selection problem by examining the size and shape of the $\GL(x)$. We proceed by induction on the difficulty of the problem, reducing hard cases to easier ones by first localizing to the correct lengthscale, then patching together local Lipschitz selections by a Whitney partition of unity. By the end of  Section 4.10 we will have proven Theorem \reff{NG-FP} and deduced Corollary \reff{A-12}, the strengthened version of Theorem  \reff{CONTH} on metric trees (without any assumption of the degree of the nodes).
\par In Section 5 we return to the setting of a general metric space $(\Mc,\rho)$ and a map $F:\Mc\to\KM$. We suppose that for every $\SM\subset\Mc$ consisting of at most $\ks$ points, the restriction $F|_{\SM}$ has a Lipschitz selection with Lipschitz seminorm at most $\lambda$. Here, $\ks$ is the same constant as in Theorem \reff{NG-FP}.
\par For each $x\in\Mc$, we define a nonempty compact convex ``core''
\bel{A-CORE}
\CRF(x)\subset F(x),
\ee
with the following crucial property:
\bel{B-CORE}
\text{For every}~~x,y\in\Mc,~~\text{the {\it Hausdorff distance} from}~~\CRF(x)~~\text{to}~~\CRF(y)~~ \text{is at most}~~\gz\lambda\,\rho(x,y).
\ee
Here, $\gz$ depends only on $m$.
\par Recall that the Hausdorff distance $\dhf(A,B)$ between two nonempty compact sets $A,B\subset Y$ is defined as the least $r\ge 0$ such that for each $x\in A$ there exists $y\in B$ such that $\|x-y\|\le r$, and for each $x\in B$ there exists $y\in A$ such that $\|x-y\|\le r$.
\smallskip
\par We define $\CRF(x)$ by considering an arbitrary finite tree $T=(X,E)$ $(X=\{\text{nodes}\}$, $E=\{\text{edges}\})$ and an arbitrary map $\psi:X\to \Mc$ such that
\bel{AD-PH}
\psi(x)\ne \psi(y)~~~\text{whenever}~~~x,y\in X~~~\text{are joined by an edge.}
\ee
\par We refer to $\psi$ as an {\it admissible} mapping. (See Definition \reff{ADM}.) The map $\psi$ induces a tree metric $d$ on $X$ by setting $d(x,y)=\rho(\psi(x),\psi(y))$ whenever $x$ and $y$ are nodes in $X$ joined by an edge.
\par Moreover, we obtain a Lipschitz selection problem for the metric tree $(X,d)$ by considering the map $F\circ \psi:X\to\KM$. From Corollary \reff{A-12} (i.e., Theorem \reff{CONTH} in its strengthened form), we learn that $F\circ \psi$ has a Lipschitz selection with Lipschitz seminorm at most $\gz\,\lambda$. By considering all such Lipschitz selections for a fixed $T=(X,E)$, a node $a\in X$, and a map $\psi:X\to\Mc$ (satisfying \rf{AD-PH}) such that $\psi(a)=x$, we define a nonempty compact convex set
$$
O(x;[T,a,\psi])\subset F(x)~~~\text{for each}~~~x\in\Mc.
$$
(See Section 5.1 for the definition of the sets $O(x;[\cdot,\cdot,\cdot])$.)
\par The ``core'' $\CRF(x)$ is then defined as the intersection of the sets $O(x;[T,a,\psi])$ over all finite trees $T=(X,E)$, all nodes $a\in X$, and all $\psi:X\to\Mc$ with $\psi(a)=x$ satisfying \rf{AD-PH}. The key properties \rf{A-CORE}, \rf{B-CORE} of $\CRF$ follow easily once we know that $O(x;[T,a,\psi])$ is nonempty, and we easily deduce that key fact from Corollary \reff{A-12}.
\smallskip
\par Once we have produced a ``core'' $\CRF(x)$ satisfying
\rf{A-CORE} and \rf{B-CORE} (see Theorem \reff{HDS-M}), we can invoke a selection theorem of Shvartsman \cite{S04}, see Theorem \reff{ST-P}. This result provides the existence of a Lipschitz (with respect to the Hausdorff distance $\dhf$) map $\ST:\KM\to Y$ such that $\ST(K)\in K$ for all $K\in\KM$. Furthermore, the $\dhf$-Lipschitz seminorm of $\ST$ is bounded by a constant depending only on $m$. We refer to $\ST(K)$ as ``Steiner-type point'' of $K$. See Section 5.2 for more detail.
\par We can now apply the Steiner-type point map $\ST$ to the core $\CRF$ to establish the following weak form of Theorem \reff{MAIN-FP}.
\begin{theorem}\lbl{KSHARP} Given $m\ge 1$ there exist constants $\ks,\gone$, depending only on $m$, for which the following holds.
\par Let $(\Mc,\rho)$ be a metric space, let $F:\Mc\to\KM$ for a Banach space $Y$, and let $\lambda$ be a positive real number. Suppose that for every $\SM\subset\Mc$ consisting of at most $\ks$ points, the restriction $F|_{\SM}$ has a Lipschitz selection $f_{\SM}$ with Lipschitz  seminorm $\|f_{\SM}\|_{\Lip(\SM,Y)}\le \lambda$.
\smallskip
\par Then $F$ has a Lipschitz selection $f$ with Lipschitz  seminorm $\|f\|_{\Lip(\Mc,Y)}\le \gone\lambda$.
\end{theorem}
\smallskip
\par Note that we have here $\ks$ instead of the sharp finiteness constant $N(m,Y)$, and that $(\Mc,\rho)$ is a metric space, rather than a pseudometric space.

\smallskip
\par To prove Theorem \reff{MAIN-FP}, it remains to pass from metric spaces to pseudometric spaces, and to pass from the large finiteness constant $\ks$ to the optimal finiteness constant $N(m,Y)$.
\par We pass to pseudometric spaces in Section 6. In the context of Theorem \reff{MAIN-FP}, the task is easy. For Theorem \reff{FINITE}, the variant of Theorem \reff{MAIN-FP} in which $(\Mc,\rho)$ is finite but the sets $F(x)$ needn't be compact, it takes a bit more work.
\par Finally, we pass from $\ks$ to $N(m,Y)$ by applying a result of Shvartsman \cite[Theorem 1.2]{S02}.
\begin{theorem}\lbl{KS-N} Let $(\wtM,\trh)$ be a finite pseudometric space, let $\tF:\wtM\to\KM$ for a Banach space $Y$, and let $\lambda$ be a positive real number.
\par Suppose that for every $S\subset\wtM$ consisting of at most $N(m,Y)$ points, the restriction $\tF|_{S}$ has a Lipschitz selection $\tf_{S}$ with Lipschitz seminorm $\|\tf_{S}\|_{\Lip(S,Y)}\le \lambda$.
\par Then $\tF$ has a Lipschitz selection $\tf$ with Lipschitz seminorm $\|\tf\|_{\Lip(\wtM,Y)}\le C(\wtM)\,\lambda$, where $C(\wtM)$ depends only on $m$ and on the number of points in $\wtM$.
\end{theorem}
\par Note that Theorem \reff{KS-N} is the ``bridge'' between the existence of some finiteness constant $\ks$ in Theorem \reff{KSHARP} and the optimal finiteness constant $N(m,Y)$ in Theorem \reff{MAIN-FP}.
\medskip
\par We combine Theorem \reff{KSHARP} (for pseudometric spaces) with Theorem \reff{KS-N}, to complete the proof of Theorem \reff{MAIN-FP}, our main result. The argument is simple: Using Theorem \reff{KS-N}, we pass from $N(m,Y)$-point subsets to $\ks$-point subsets; then, using Theorem \reff{KSHARP}, we pass from $\ks$-point subsets to a full solution of our Lipschitz selection problem.
\par Finally, Section 7 states two variants of Theorem \reff{MAIN-FP}, and adds a few closing remarks.
\medskip
\par As in \cite{FIL1}, our present results lead to questions about efficient computation for Lipschitz selection problems on finite metric spaces. In connection with such issues, we ask whether the results of Har-Peled and Mendel \cite{HP-M} on the Well Separated Pairs Decomposition \cite{C-Kos} can be extended from doubling metrics to metrics of bounded Nagata dimension.
\par Readers interested in checking details of our proofs may want to consult a much more detailed version of this paper posted on the arXiv \cite{FS-2017}. We mention also that A. Brudnyi \cite{AB-PR} has advised us that he has an alternate proof of the passage from the finiteness principle for metric trees to the construction of the core.
\medskip

\SECT{2. Whitney partitions and Patching Lemma}{2}
\addtocontents{toc}{2. Whitney partitions and Patching Lemma \hfill\thepage\par\VST}

\indent\par Let $(X,d)$ be a metric space. We write $B(x,r)$ to denote the ball $\{y\in X:d(x,y)<r\}$
(strict inequality) in the metric space $(X,d)$. We also write $\diam A=\sup\,\{d(a,b):a,b\in A\}$ and
$$
\dist(A',A'')=\inf\{d(a',a''): a'\in A', a''\in A''\}
$$
to denote the diameter of a set $A\subset X$ and the distance between sets $A',A''\subset X$ respectively.
\smallskip
\medskip

\indent\par {\bf 2.1 Whitney partitions on metric spaces with finite Nagata dimension.}
\addtocontents{toc}{~~~~2.1 Whitney partitions on metric spaces with finite Nagata dimension.\hfill \thepage\par\VSU}\medskip
\par In this section, we prove the following result.
\begin{wplemma}\lbl{WPL} Let $(X,d)$ be a metric space,
and let $r(x)>0$ be a positive function on $X$. We assume the following, for constants $\cng\in(0,1]$, $\CNG\in\N\cup\{0\}$ and $\CLS\ge 1$:\vskip 1mm
\par \textbullet~({\sc Nagata $(\CNG,\cng)$}) Given $s>0$ there exists a covering of $X$ by subsets $X_i$ $(i\in I)$ of diameter at most $s$, such that every ball of radius
$\cng s$ in $X$ meets at most $\CNG+1$ of the $X_i$\,.
\smallskip
\par \textbullet~({\sc Consistency of the Lengthscale}) Let $x,y\in X$. If $d(x,y)\le r(x)+r(y)$, then
\bel{C-LSC}
\CLS^{-1}\,r(x)\le r(y)\le \CLS r(x).
\ee
\par Let $a>0$.
\medskip
\par Then there exist functions $\vf_\nu:X\to\R$, and points $x_\nu\in X$, with the following properties: \smallskip
\par \textbullet~Each $\vf_\nu\ge 0$, and each
$\vf_\nu=0$ outside $B(x_\nu,ar_\nu)$. Here and below, $r_\nu=r(x_\nu)$.
\par \textbullet~ Any given $x\in X$ satisfies $\vf_\nu(x)\ne 0$ for at most $D^*$ distinct $\nu$.
\par \textbullet~ $\sum\limits_\nu\,\vf_\nu=1$ on $X$.
\par \textbullet~For each $\nu$ and for all $x,y\in X$, we have
$$
|\vf_\nu(x)-\vf_\nu(y)|\le \CWH\,d(x,y)/r_\nu.
$$
\par Here $D^*$ and $\CWH$ are constants depending only on $\cng$, $\CNG$, $\CLS$ and $a$.
\end{wplemma}
\par {\it Proof.} We write $c,C$ to denote positive constants determined by $\cng$, $\CNG$, $\CLS$ and $a$. These symbols may denote different constants in different occurrences.
\par We introduce a large constant $A$ to be fixed later. We make the following
\begin{lassum}\lbl{WP-LA} {\em $A$ exceeds a large enough constant determined by $\cng$, $\CNG$, $\CLS$, $a$.}
\end{lassum}
\par We write $c(A), C(A)$ to denote positive constants determined by $A$, $\cng$, $\CNG$, $\CLS$, $a$. These symbols may denote different constants in different occurrences.
\par Let $P$ denote the set of all integer powers of $2$, including negative powers. For $s\in P$ let $(X(i,s))_{i\in I(s)}$ be a covering of $X$ given by the Nagata $(\CNG,\cng)$ condition. Thus,
\medskip
$$
\diam X(i,s)\le s;
$$
and, for fixed $s\in P$,
\bel{2-C}
\text{any given}~~~x\in X~~~\text{lies in at most}~~C~~\text{of the sets}~~X^{++}(i,s).
\ee
Here
$$
X^{++}(i,s)=\{y\in X: d(y,X(i,s))<\cng s/64\}~~~~(i\in I(s)).
$$
\par We also define
$$
X^{+}(i,s)=\{y\in X: d(y,X(i,s))<\cng s/128\} ~~~\text{for}~~~(i\in I(s)).
$$
\par Let
$$
\theta_{i,s}(x)=\max\{0,(1-256\,d(x,X(i,s))/(\cng s))\}
$$
for $x\in X$, $i\in I(s)$, $s\in P$.
\par Then
\bel{3a-C}
0\le\theta_{i,s}\le 1,
\ee
\bel{3b-C}
\|\theta_{i,s}\|_{\Lip(X,\R)}\le C\,s^{-1},
\ee
and
$$
\theta_{i,s}=0~~~\text{outside}~~~X^+(i,s),
$$
but
\bel{B-S9}
\theta_{i,s}=1~~~\text{on}~~~X(i,s).
\ee
\par For each $s\in P$ and $i\in I(s)$, we pick a representative point $x(i,s)\in X(i,s)$. (We may assume that the $X(i,s)$ are all nonempty.) We let \REL (relevant) denote the set of all $(i,s)$ such that
\bel{4-C}
A^{-3}r(x(i,s))\le s\le A^{-1}r(x(i,s)).
\ee
\par We establish the basic properties of the set \REL.
\begin{lemma}\lbl{WP1} Given $x_0\in X$ there exists $(i,s)\in\TRL$ such that $x_0\in X(i,s)$ and therefore $\theta_{i,s}(x_0)=1$.
\end{lemma}
\par {\it Proof.} The ``therefore'' part of the lemma follows from \rf{B-S9}.
\par Pick $s_0\in P$ such that
$$
s_0/2\le r(x_0)/A^2\le 2s_0.
$$
\par Because the $X(i,s_0)$ $(i\in I(s_0))$ cover $X$, we may fix $i_0\in I(s_0)$ such that $x_0\in X(i_0,s_0)$. The points $x_0$ and $x(i_0,s_0)$ both belong to $X(i_0,s_0)$, hence
$$
d(x_0,x(i_0,s_0))\le \diam X(i_0,s_0)\le s_0\le 2r(x_0)/A^2\,.
$$
\par The Large $A$ Assumption \reff{WP-LA} and the {\sc Consistency of the Lengthscale} together now imply that
$$
cr(x_0)\le r(x(i_0,s_0))\le Cr(x_0),
$$
and therefore
$$
cs_0\le r(x(i_0,s_0))/A^2\le Cs_0\,.
$$
Thanks to the Large $A$ Assumption \reff{WP-LA}, we therefore have \rf{4-C} for $(i_0,s_0)$. Thus, $(i_0,s_0)\in\TRL$ and $x_0\in X(i_0,s_0)$.\bx
\begin{lemma}\lbl{WP2} If $(i,s)\in\TRL$ and $x_0\in X^{++}(i,s)$, then
$$
cA^{-3}r(x_0)\le s\le CA^{-1}r(x_0),
$$
and therefore
$$
\|\theta_{i,s}\|_{\Lip(X,\R)}\le CA^{3}/r(x_0).
$$
\end{lemma}
\par {\it Proof.} Both $x_0$ and $x(i,s)$ lie in $X^{++}(i,s)$, hence
$$
d(x_0,x(i,s))\le \diam X^{++}(i,s)\le 2\cng s/64+
\diam X(i,s)\le Cs\le Cr(x(i,s))/A
$$
thanks to \rf{4-C}.
\par The Large $A$ Assumption \reff{WP-LA} and {\sc Consistency of the Lengthscale} now tell us that
$$
cr(x_0)\le r(x(i,s))\le Cr(x_0),
$$
and therefore \rf{4-C} and \rf{3b-C} imply the conclusion of Lemma \reff{WP2}.\bx
\begin{corollary}\lbl{WP-C1} Any given point $x_0\in X$ lies in $X^{++}(i,s)$ for at most  $C(A)$ distinct $(i,s)\in \TRL$. Consequently, $\theta_{i,s}(x_0)$ is nonzero for at most $C(A)$ distinct $(i,s)\in\TRL$.
\end{corollary}
\par {\it Proof.} There are at most $C(A)$ distinct $s\in P$ satisfying the conclusion of Lemma \reff{WP2}.
For each such $s$ there are at most $C$ distinct $i$ such that $x_0\in X^{++}(i,s)$; see \rf{2-C}.\bx
\begin{corollary}\lbl{WP-C2} Suppose $X^{++}(i,s)\cap X^{++}(i_0,s_0)\ne\emp$ with $(i,s),(i_0,s_0)\in\TRL$. Then 
$$
c(A)s_0\le s\le C(A)s_0.
$$
\end{corollary}
\par {\it Proof.} Pick $x_0\in X^{++}(i,s)\cap X^{++}(i_0,s_0)$. Lemma \reff{WP2} gives
$$
c(A)r(x_0)\le s\le C(A)r(x_0)~~~\text{and}~~~
c(A)r(x_0)\le s_0\le C(A)r(x_0). \BXD
$$
\begin{lemma}\lbl{WP3} Let $(i_0,s_0), (i,s)\in\TRL$. If $x\in X^+(i_0,s_0)$, then for any $y\in X$
\bel{5-C}
|\theta_{i,s}(x)-\theta_{i,s}(y)|\le C(A)\,d(x,y)/s_0.
\ee
\end{lemma}
\par {\it Proof.} We proceed by cases.\medskip
\par {\it Case 1:} $d(x,y)<\cng s_0/128$. \medskip
\par Then $x,y\in X^{++}(i_0,s_0)$. If $x$ or $y$ belongs to $X^{++}(i,s)$, then Corollary \reff{WP-C2} tells us that 
$$
c(A)s_0\le s\le C(A)s_0\,;
$$
hence, \rf{3b-C} yields the desired estimate \rf{5-C}.
\par If instead neither $x$ nor $y$ belongs to $X^{++}(i,s)$, then $\theta_{i,s}(x)=\theta_{i,s}(y)=0$, hence \rf{5-C} holds trivially.\medskip
\par {\it Case 2:} $d(x,y)\ge \cng s_0/128$. Then \rf{3a-C} gives
$$
|\theta_{i,s}(x)-\theta_{i,s}(y)|\le 1\le C\,d(x,y)/s_0.
$$
Thus, \rf{5-C} holds in all cases.\bx
\medskip
\par Now define
\bel{6-C}
\Theta(x)=\smed_{(i,s)\in\TRL}\,\theta_{i,s}(x)~~~~
\text{for all}~~~x\in X.
\ee
\par Corollary \reff{WP-C1} shows that there are at most $C(A)$ nonzero summands in \rf{6-C} for any fixed $x$. Moreover, each summand is between $0$ and $1$ (see \rf{3a-C}), and for each fixed $x$, at least one of the summands is equal to $1$ (see Lemma \reff{WP1}). Therefore,
\bel{7-C}
1\le \Theta(x)\le C(A)~~~~\text{for all}~~~~x\in X.
\ee
\begin{lemma}\lbl{WP4} Let $x,y\in X$ and $(i_0,s_0)\in\TRL$. If $x\in X^+(i_0,s_0)$, then
$$
|\Theta(x)-\Theta(y)|\le C(A)\,d(x,y)/s_0\,.
$$
\end{lemma}
\par {\it Proof.} There are at most $C(A)$ distinct $(i,s)\in\TRL$ for which $\theta_{i,s}(x)$ or $\theta_{i,s}(y)$ is nonzero. For each such $(i,s)$ we apply Lemma \reff{WP3}, then sum over $(i,s)$.\bx
\medskip
\par Now, for $(i_0,s_0)\in\TRL$, we set
\bel{8-C}
\vf_{i_0,s_0}(x)=\theta_{i_0,s_0}(x)/\Theta(x)\,.
\ee
\par This function is defined on all of $X$, and it is zero outside $X^+(i_0,s_0)$. Moreover,
\bel{9-C}
\vf_{i_0,s_0}\ge 0~~~~~\text{and}~~~~~
\smed_{(i_0,s_0)\in\TRL}\vf_{i_0,s_0}=1~~~\text{on}~~~X.
\ee
\par Note that because
$$
\diam X^+(i_0,s_0)\le Cs_0\le C\,A^{-1}r(x(i_0,s_0))
$$
(see \rf{4-C}), the function $\vf_{i_0,s_0}$ is zero outside the ball $B(x(i_0,s_0),C\,A^{-1}r(x(i_0,s_0)))$.
Thanks to our Large $A$ Assumption \reff{WP-LA}, it follows that
\bel{10-C}
\vf_{i,s}~~~\text{is identically zero outside the ball} ~~~B(x(i,s),ar(x(i,s)))\,.
\ee
\begin{lemma}\lbl{WP5} For $x,y\in X$ and $(i_0,s_0)\in\TRL$, we have
$$
|\vf_{i_0,s_0}(x)-\vf_{i_0,s_0}(y)|\le C(A)\,d(x,y)/s_0.
$$
\end{lemma}
\par {\it Proof.} Suppose first that $x\in X^+(i_0,s_0)$. Then
$$
|\vf_{i_0,s_0}(x)-\vf_{i_0,s_0}(y)|=
\left|\frac{\theta_{i_0,s_0}(x)}{\Theta(x)}-
\frac{\theta_{i_0,s_0}(y)}{\Theta(y)}\right|\nn\\
\le
\frac{|\theta_{i_0,s_0}(x)-\theta_{i_0,s_0}(y)|}{\Theta(x)}+
\theta_{i_0,s_0}(y)\,
\frac{|\Theta(x)-\Theta(y)|}{\Theta(x)\Theta(y)}\,.
$$
The first term on the right is at most $C(A)\,d(x,y)/s_0$ by \rf{3b-C} and \rf{7-C}; the second term on the right is at most $C(A)\,d(x,y)/s_0$ thanks to \rf{3a-C}, Lemma \reff{WP4} and \rf{7-C}.
Thus,
\bel{11-C}
|\vf_{i_0,s_0}(x)-\vf_{i_0,s_0}(y)|\le C(A)\,d(x,y)/s_0
~~~~\text{if}~~~x\in X^+(i_0,s_0).
\ee
Similarly, \rf{11-C} holds if $y\in X^+(i_0,s_0)$.
\par Finally, if neither $x$ nor $y$ belongs to $X^+(i_0,s_0)$, then $$\vf_{i_0,s_0}(x)=\vf_{i_0,s_0}(y)=0,$$ so \rf{11-C} is obvious.
\par Thus, \rf{11-C} holds in all cases.\bx
\begin{corollary}\lbl{WP-C3} For $x,y\in X$ and $(i_0,s_0)\in\TRL$, we have
$$
|\vf_{i_0,s_0}(x)-\vf_{i_0,s_0}(y)|\le C(A)\,d(x,y)/r(x(i_0,s_0)).
$$
\end{corollary}
\par {\it Proof.} Immediate from Lemma \reff{WP5} and inequalities \rf{4-C}.\bx
\bigskip
\par We can now finish the proof of the Whitney Partition Lemma \reff{WPL}. We pick $A$ to be a constant determined by $\cng$, $\CNG$, $\CLS$, $a$, taken large enough to satisfy the Large $A$ Assumption \reff{WP-LA}. We then take our functions $\vf_\nu$ to be the $\vf_{(i,s)}$ $((i,s)\in\TRL)$, and we take our $x_\nu$ to be the points $x(i,s)$  $((i,s)\in\TRL)$. We set $r_\nu=r({x_\nu})$.
\smallskip
\par The following hold: \smallskip
\par \textbullet~Each $\vf_\nu\ge 0$, and each
$\vf_\nu=0$ outside $B(x_\nu,ar_\nu)$; see \rf{9-C} and \rf{10-C}.
\smallskip
\par \textbullet~ Any given $x\in X$ satisfies $\vf_\nu(x)\ne 0$ for at most $C$ distinct $\nu$.
This follows from Corollary \reff{WP-C1}, definition \rf{8-C}, and the fact that $A$ is now determined by $\cng$, $\CNG$, $\CLS$, $a$.\smallskip
\par \textbullet~ $\sum\limits_\nu\,\vf_\nu=1$ on $X$; see \rf{9-C}.
\smallskip
\par \textbullet~For each $\nu$ and for all $x,y\in X$, we have
$$
|\vf_\nu(x)-\vf_\nu(y)|\le C\,d(x,y)/r_\nu;
$$
see Corollary \reff{WP-C3}, and note that $A$ is now determined by $\cng$, $\CNG$, $\CLS$ and $a$.
\smallskip
\par The proof of the Whitney Partition Lemma \reff{WPL} is complete.\bx
\smallskip
\begin{remark} {\em Later on there will be another Large $A$ Assumption different from that in this section.}
\end{remark}

\indent\par {\bf 2.2 Patching Lemma.}
\addtocontents{toc}{~~~~2.2 Patching Lemma.\hfill \thepage\par\VSU}\medskip
\begin{ptchlemma}\lbl{PTHM} Let $(X,d)$ be a metric space, and let $Y$ be a Banach space. For each $\nu$ in some index set, assume we are given the following objects:\smallskip
\par \textbullet~ A point $x_\nu\in X$ and a positive number $r_\nu>0$ (a ``lengthscale'').
\smallskip
\par \textbullet~ A function $\theta_\nu:X\to \R$\,.
\smallskip
\par \textbullet~ A vector $\eta_\nu\in Y$ and a vector-valued function $F_\nu:X\to Y$.
\medskip
\par We make the following assumptions: We are given positive constants $\CLS\ge 1$, $\CWH$, $\CETA$, $\CAGR$, $\CLIP$, $\DS$, such that the following conditions are satisfied for each $\mu,\nu$
\smallskip
\par \textbullet~ ({\sc Consistency of the Lengthscale})
\bel{C-LS}
\CLS^{-1}\le r_\nu/r_\mu\le \CLS~~~
\text{whenever}~~~d(x_\mu,x_\nu)\le r_\mu+r_\nu.
\ee
\par ({\sc Whitney Partition Assumptions})
\smallskip
\par \textbullet~ $\theta_\nu\ge 0$ on $X$ and
$\theta_\nu= 0$ outside $B(x_\nu,a\,r_\nu)$, where
\bel{A-SMALL}
a=(4\,\CLS)^{-1}.
\ee
\par \textbullet~
$|\theta_\nu(x)-\theta_\nu(y)|\le \CWH\cdot d(x,y)/r_\nu$~
for~ $x,y\in X$.
\smallskip
\par \textbullet~ Any given $x\in X$ satisfies $\theta_\nu(x)\ne 0$~ for at most $\DS$ distinct $\nu$.
\smallskip
\par \textbullet~ $\smed\limits_\nu\,\theta_\nu=1$ on $X$.
\smallskip
\par \textbullet~ ({\sc Consistency of the $\eta_\nu$})~
$\|\eta_\mu-\eta_\nu\|\le \CETA\cdot [r_\nu+r_\nu+d(x_\mu,x_\nu)]$.
\smallskip
\par \textbullet~ ({\sc Agreement of $F_\nu$ with $\eta_\nu$})~
$\|F_\nu(x)-\eta_\nu\|\le \CAGR\, r_\nu$~ for~
$x\in B(x_\nu,r_\nu)$.
\smallskip
\par \textbullet~ ({\sc Lipschitz continuity of $F_\nu$})~
$\|F_\nu(x)-F_\nu(y)\|\le \CLIP\cdot d(x,y)$~ for~
$x,y\in B(x_\nu,r_\nu)$.
\bigskip
\par Define
$$
F(x)=\sum_\nu\theta_\nu(x)\,F_\nu(x)~~~\text{for}~~x\in X.
$$
Then $F$ satisfies
$$
\|F(x)-F(y)\|\le C\,d(x,y)~~~\text{for}~~~x,y\in X,
$$
where $C$ is determined by $\CLS$, $\CWH$, $\CETA$, $\CAGR$, $\CLIP$, $\DS$.
\end{ptchlemma}
\medskip
\par To start the proof of the {\sc Patching Lemma} \reff{PTHM}, we define a set of relevant $\nu$ by setting
$$
\RLV(x)=\{\nu: \theta_\nu(x)\ne 0\}, ~~~x\in X.
$$

Then $1\le\#(\RLV(x))\le D^*$, and
\bel{D-NU}
d(x,x_\nu)\le a\,r_\nu~~~\text{for}~~~v\in\RLV(x).
\ee
\par We also recall that $\CLS\ge 1$ and $a=(4\,\CLS)^{-1}$ so that
\bel{CLS-A}
\CLS\cdot a= 1/4~~~\text{and}~~~a\le 1/4.
\ee
\vskip 2mm
\par We will use the following result.
\begin{lemma}\lbl{STR} Let $\nu,\nu_0\in\RLV(x)$, $\mu_0\in\RLV(y)$, and suppose that $d(x,y)\le a\cdot[r_{\nu_0}+r_{\mu_0}]$. Then
$$
x,y\in B(x_\nu,r_\nu)\cap B(x_{\nu_0},r_{\nu_0})\cap B(x_{\mu_0},r_{\mu_0})
$$
and the ratios
$$
r_{\nu_0}/r_{\mu_0},~r_{\mu_0}/r_{\nu_0},~
r_{\nu}/r_{\nu_0},~ r_{\nu_0}/r_{\nu},~
r_{\nu}/r_{\mu_0},~ r_{\mu_0}/r_{\nu}
$$
are at most $\CLS$.
\end{lemma}
\par {\it Proof.} We have the following inequalities:
\smallskip
\par $(\bigstar 1)$~
$d(x_\nu,x_{\nu_0})\le  d(x_\nu,x)+d(x,x_{\nu_0})\le a\,r_\nu+a\,r_{\nu_0}$,
\smallskip
\par $(\bigstar 2)$~
$d(x_{\nu_0},x_{\mu_0})\le  d(x_{\nu_0},x)+d(x,y)+d(y,x_{\mu_0})\le a\,r_{\nu_0}+[a\,r_{\nu_0}+a\,r_{\mu_0}]+a\,r_{\mu_0}$,
\smallskip
\par $(\bigstar 3)$~
$d(x_{\nu},x_{\mu_0})\le  d(x_{\nu},x)+d(x,y)+d(y,x_{\mu_0})\le a\,r_{\nu}+[a\,r_{\nu_0}+a\,r_{\mu_0}]+a\,r_{\mu_0}$.
\medskip
\par From $(\bigstar 1)$, $(\bigstar 2)$, \rf{CLS-A}, and {\sc Consistency of the Lengthscale} \rf{C-LS}, we have
$$
r_{\nu}/r_{\nu_0},~ r_{\nu_0}/r_{\nu},~
r_{\nu_0}/r_{\mu_0},~r_{\mu_0}/r_{\nu_0}\le \CLS.
$$
Therefore, $(\bigstar 3)$ and \rf{CLS-A} imply that
$$
d(x_{\nu},x_{\mu_0})
\le a\,r_{\nu}+\CLS a\,r_{\nu}+2a\,r_{\mu_0}
\le r_{\nu}+r_{\mu_0},
$$
and, consequently, another application of {\sc Consistency of the Lengthscale} \rf{C-LS} gives
$$
r_{\nu}/r_{\mu_0},~ r_{\mu_0}/r_{\nu}\le \CLS.
$$
\par Next, note that, by \rf{D-NU} and \rf{CLS-A},
$$
d(x,x_{\nu})\le a\,r_{\nu}<r_{\nu}
$$
and
$$
d(y,x_{\nu})\le d(y,x)+ d(x,x_{\nu})\le  [a\,r_{\nu_0}+a\,r_{\mu_0}]+a\,r_{\nu}\le
(3\CLS\,a)r_\nu< r_\nu.
$$
Hence,
$$
x,y\in B(x_\nu,r_\nu).
$$
\smallskip
\par Similarly,
$$
d(x,x_{\nu_0})\le  a\,r_{\nu_0}<r_{\nu_0}
$$
and
$$
d(y,x_{\nu_0})\le d(y,x)+ d(x,x_{\nu_0})\le  [a\,r_{\mu_0}+a\,r_{\nu_0}]+a\,r_{\nu_0}\le
(3\CLS\,a)r_{\nu_0}< r_{\nu_0}.
$$
Hence,
$$
x,y\in B(x_{\nu_0},r_{\nu_0}).
$$
\par Finally,
$$
d(y,x_{\mu_0})\le  a\,r_{\mu_0}< r_{\mu_0}
$$
and
$$
d(x,x_{\mu_0})\le d(x,y)+ d(y,x_{\mu_0})\le  [a\,r_{\mu_0}+a\,r_{\nu_0}]+a\,r_{\mu_0}\le
(3\CLS\,a)r_{\mu_0}< r_{\mu_0}.
$$
Hence,
$$
x,y\in B(x_{\mu_0},r_{\mu_0}).
$$
\par The proof of the lemma is complete.\bx

\bigskip
\par {\it Proof of the Patching Lemma \reff{PTHM}.}
\smallskip
\par We write $c,C,C'$, etc. to denote positive constants determined by $\CLS$, $\CWH$, $\CETA$, $\CAGR$, $\CLIP$, $\DS$. These symbols may denote different constants in different occurrences.
\par Let $x,y\in X$ be given. We must show that
$$
\|F(x)-F(y)\|\le C\,d(x,y).
$$
\par Fix $\mu_0,\nu_0$, with $\nu_0\in\RLV(x)$ and $\mu_0\in\RLV(y)$. We distinguish two cases.
\medskip
\par {\it CASE 1:} Suppose
$$
d(x,y)\le a\cdot[r_{\nu_0}+r_{\mu_0}]~~~\text{with}~~~
a=(4\,\CLS)^{-1}.
$$
Then Lemma \reff{STR} yields
\bel{E-L2}
x,y\in B(x_\nu,r_\nu)\cap B(x_{\nu_0},r_{\nu_0})\cap B(x_{\mu_0},r_{\mu_0})
\ee
for all $\nu\in \RLV(x)\cup\RLV(y)$. (If $\nu\in\RLV(y)$,
we apply Lemma \reff{STR} with $y,x,\mu_0,\nu_0$ in place of $x,y,\nu_0,\mu_0$.) Also, for such $\nu$, Lemma \reff{STR} gives
\bel{RMN}
c\,r_{\nu_0}\le r_{\nu}\le C\,r_{\nu_0}
~~~\text{and}~~~c\,r_{\nu_0}\le r_{\mu_0}\le C\,r_{\nu_0}. \ee
\par For $v\in\RLV(x)$, we have
\bel{Y-BL}
\|F_\nu(y)-\eta_{\nu_0}\|\le \|F_\nu(y)-\eta_{\nu}\|+
\|\eta_\nu-\eta_{\nu_0}\|\le C\,r_\nu+C\,[r_\nu+r_{\nu_0}+d(x_\nu,x_{\nu_0})].
\ee
(Here, we may apply {\sc Consistency
of the $\eta_\nu$} and {\sc Agreement of $F_\nu$ with $\eta_\nu$}, because $y\in B(x_\nu,r_{\nu})$.) Also, by \rf{E-L2},
$$
d(x_\nu,x_{\nu_0})\le d(x_\nu,x)+d(x,x_{\nu_0})\le r_{\nu}+r_{\nu_0}~~~\text{for}~~~
\nu\in\RLV(x).
$$
\par The above estimates and \rf{RMN} tell us that
$$
\|F_\nu(y)-\eta_{\nu_0}\|\le C\,r_{\nu_0}~~~\text{if}~~~\nu\in\RLV(x).
$$
\par Similarly, suppose $v\in\RLV(y)$. Then \rf{Y-BL} holds. (We may apply {\sc Agreement of $F_\nu$ with $\eta_\nu$}, because $y\in B(x_\nu,r_{\nu})$.)
Also, by \rf{E-L2},
$$
d(x_\nu,x_{\nu_0})\le d(x_\nu,y)+d(y,x_{\nu_0})\le r_{\nu}+\,r_{\nu_0}
~~\text{for all}~~\nu\in\RLV(y).
$$
\par The above estimates and \rf{RMN} tell us that
$$
\|F_\nu(y)-\eta_{\nu_0}\|\le C\,r_{\nu_0}~~~\text{for all}~~~\nu\in\RLV(y).
$$
Thus,
$$
\|F_\nu(y)-\eta_{\nu_0}\|\le C\,r_{\nu_0}~~~\text{for all}~~~\nu\in\RLV(x)\cup\RLV(y).
$$
\par We now write
$$
F(x)-F(y)=\smed_{\nu\in\RLV(x)\cup\RLV(y)}
\theta_\nu(x)\cdot[F_\nu(x)-F_\nu(y)]+
\smed_{\nu\in\RLV(x)\cup\RLV(y)}
[\theta_\nu(x)-\theta_\nu(y)]\cdot[F_\nu(y)-\eta_{\nu_0}]
\equiv I+II.
$$
Here we have used the {\sc Whitney Partition Assumption} that the sum of all $\theta_\nu$ equals $1$.
\par It follows from Lipschitz continuity of $F_\nu$ that
$$
\|I\|\le \smed_{\nu\in\RLV(x)\cup\RLV(y)}
\theta_\nu(x)\cdot[C\,d(x,y)]=C\,d(x,y).
$$
\par Each summand in $II$ satisfies
$$
|\theta_\nu(x)-\theta_\nu(y)|\le C\,d(x,y)/r_\nu
~~~\text{and}~~~
\|F_\nu(y)-\eta_{\nu_0}\|\le C\,r_{\nu_0}\le C\,\CLS\,r_\nu\,,
$$
see \rf{C-LSC}. Hence
$$
\|[\theta_\nu(x)-\theta_\nu(y)]
\cdot[F_\nu(y)-\eta_{\nu_0}]\|
\le C\,d(x,y).
$$
\par Because there are at most $2\DS$ summands in $II$, it follows that
$$
\|II\|\le C\,d(x,y).
$$
Combining our estimates for terms $I$ and $II$, we find that
$$
\|F(x)-F(y)\|\le C\,d(x,y)~~~\text{in CASE 1.}
$$
\par {\it CASE 2:} Suppose
$$
d(x,y)> a\cdot[r_{\nu_0}+r_{\mu_0}]~~~\text{with}~~~
a=(4\,\CLS)^{-1}.
$$
For $\nu\in\RLV(x)$, we have
$$
d(x_\nu,x_{\nu_0})\le  d(x_\nu,x)+d(x,x_{\nu_0}) \le a\cdot r_{\nu}+a\cdot r_{\nu_0},
$$
hence, by {\sc Consistency of the Lengthscale} (see \rf{C-LSC}),
$$
c\,r_{\nu_0}\le r_{\nu}\le C\,r_{\nu_0}
$$
and
$$
\|F_\nu(x)-\eta_{\nu_0}\|\le \|F_\nu(x)-\eta_{\nu}\|+
\|\eta_\nu-\eta_{\nu_0}\|\le C\,r_\nu+[C\,r_\nu+C\,r_{\nu_0}+Cd(x_\nu,x_{\nu_0})]
\le
C\,r_{\nu_0}.
$$
Here we use \rf{D-NU} and \rf{CLS-A}, and the fact that
$x\in B(x_\nu,r_\nu)$.
\par Consequently,
$$
\|F(x)-\eta_{\nu_0}\|=\left\|\smed_{v\in\RLV(x)}
\theta_\nu(x)\cdot[F_\nu(x)-\eta_{\nu_0}]\right\|\le
C\,r_{\nu_0}
\smed_{v\in\RLV(x)}\theta_\nu(x)=C\,r_{\nu_0}.
$$
\par Similarly,
$$
\|F(y)-\eta_{\mu_0}\|\le C\,r_{\mu_0}.
$$
Therefore,
\be
\|F(x)-F(y)\|&\le& C\,r_{\nu_0}+C\,r_{\mu_0}+\|\eta_{\nu_0}-\eta_{\mu_0}\|
\le
C'\,r_{\nu_0}+C'\,r_{\mu_0}+C'd(x_{\nu_0},x_{\mu_0})\nn\\
&\le&
C'\,r_{\nu_0}+C'\,r_{\mu_0}+C'[d(x_{\nu_0},x)+d(x,y)+
d(y,x_{\mu_0})]\nn\\
&\le&
C''\,r_{\nu_0}+C''\,r_{\mu_0}+C''d(x,y).\nn
\ee
\par Moreover, because we are in CASE 2, we have
$$
r_{\nu_0}+r_{\mu_0}\le\tfrac1a\,d(x,y)=4\,\CLS\,d(x,y).
$$
\par It now follows that
$$
\|F(x)-F(y)\|\le C'''d(x,y)~~~\text{in\, CASE}~2.
$$
\par Thus, the conclusion of the {\sc Patching Lemma} holds in all cases.\bx

\SECT{3. Sets $\Gamma_\ell$\,, labels and bases}{3}
\addtocontents{toc}{3. Sets $\Gamma_\ell$, labels and bases \hfill\thepage\par\VST}

\indent\par {\bf 3.1 Main properties of $\Gamma_\ell$\,.}
\addtocontents{toc}{~~~~3.1 Main properties of $\Gamma_\ell$. \hfill \thepage\par}\medskip
\indent\par We recall that $(Y,\|\cdot\|)$ denotes a Banach space. Given a convex set $S\subset\BS$ we let $\affspan(S)$ denote the affine hull of $S$, i.e., the smallest (with respect to inclusion) affine subspace of $\BS$ containing  $S$. We define the affine dimension $\dim S$ of $S$ as the dimension of its affine hull, i.e.,
$$
\dim S=\dim\affspan (S).
$$
\par Given $y\in\BS$ and $r>0$ we let
$$
\BL(y,r)=\{z\in \BS:~\|z-y\|\le r\}
$$
denote a {\it closed} ball in $\BS$ with center $y$ and radius $r$. By $\BY=\BL(0,1)$ we denote the unit ball in $\BS$.
\par Given non-empty sets $A, B\subset Y$ we let  $A+B=\{a+b: a\in A, b\in B\}$ denote the Minkowski sum of these sets. Given a positive real number $\lambda$ by $\lambda A$ we denote the set $\lambda A=\{\lambda a: a\in A\}$.
\smallskip
\par We call a pseudometric space $(\MS,\dm)$ finite if $\Mc$ is finite, but we say that the pseudometric $\rho$ is finite if $\rho(x,y)$ is finite for every $x,y\in\Mc$.
\par Let $(\Mc,\rho)$ be a {\it finite} pseudometric space
with a finite pseudometric $\rho$. Let us fix a constant $\lambda>0$, an integer $m\ge 0$, and a set-valued mapping $F:\Mc\to\CNMY$. Recall that $\CNMY$ denotes the family of all nonempty convex subsets of $Y$ of dimension at most $m$.
\smallskip
\par In this section we introduce a family of convex sets $\GL(x)\subset Y$ parametrized by $x\in \Mc$ and a non-negative integer $\ell$. To do so, we first define integers $k_0,k_1,k_2,...$ by the formula
\bel{KEL}
k_\ell=(m+2)^\ell~~~(\ell\ge 0).
\ee
\begin{definition}\lbl{GML-D} {\em Let $x\in \Mc$ and let $S\subset \Mc$. A point $\xi\in Y$ belongs to the set $\Gamma(x,S)$ if there exists a mapping $f:S\cup\{x\}\to Y$ such that:
\par (i) $f(x)=\xi$ and $f(z)\in F(z)$ for all $z\in S\cup\{x\}$;
\par (ii) For every $z,w\in S\cup\{x\}$ the following inequality
$$
\|f(z)-f(w)\|\le \lambda\,\rho(z,w)
$$
holds.}
\end{definition}
\par We then define
\bel{BC-3}
\GL(x)=\bigcap_{\substack{S\subset \Mc\\\#S\le k_\ell}}
\Gamma(x,S)~~~\text{for}~~~x\in \Mc,~\ell\ge 0.
\ee
\par For instance, given $x\in\Mc$ let us present an explicit formula for $\Gamma_0(x)$. By \rf{BC-3} for $\ell=0$,
$$
\Gamma_0(x)=\bigcap_{S\subset\Mc,\,\,\#S\le 1}
\Gamma(x,S).
$$
Clearly, by Definition \reff{GML-D},
$$
\Gamma(x,\{z\})=F(x)\,\capsm
\left(F(z)+\lambda\,\rho(x,z)\BY\right)~~\text{for every}~~z\in\Mc\,,
$$
and $\Gamma(x,\emp)=F(x)$, so that
\bel{GM-ZR}
\Gamma_0(x)=\bigcap_{z\in\Mc}
\left(F(z)+\lambda\,\rho(x,z)\BY\right)\,.
\ee
\begin{remark} {\em (i) Of course, the sets $\GL(x)$ also depend on the set-valued mapping $F$, the constant $\lambda$ and $m$. However, we use $\Gamma$'s only in this section, Sections 3-4 and Section 6.2 where these objects, i.e., $F$, $\lambda$ and $m$, are clear from the context. Therefore we omit $F$, $\lambda$ and $m$ in the notation of $\Gamma$'s.
\smallskip
\par (ii) As in the statement of Theorem \reff{RM-FP}, we may restrict attention to $S$ containing exactly $k_\ell$ points in \rf{BC-3}, provided $\Mc$ contains at least $k_\ell$ points.}
\end{remark}
\smallskip
\par The above $\Gamma'\text{s}$ are (possibly empty) {\it convex subsets} of $Y$. Note that
\bel{BC-4a}
\Gamma(x,S)\subset F(x)~~~\text{for all}~~~x\in\Mc
~~~\text{and}~~~S\subset\Mc.
\ee
Hence,
\bel{AFF-4}
\Gamma(x,S)\subset \affspan(F(x)) ~~~~~~x\in\Mc,~S\subset\Mc.
\ee
\par From \rf{BC-4a} and \rf{BC-3} we obtain
\bel{GL-F2}
\GL(x)\subset F(x)~~~~\text{for}~~~x\in \Mc,~ \ell\ge 0.
\ee
Also, obviously,
\bel{BC-5}
\GL(x)\subset \GLO(x)~~~~\text{for}~~~x\in \Mc,~ \ell\ge 1.
\ee
\par We describe main properties of the sets $\GL$ in Lemma \reff{G-AB} below. The proof of this lemma relies on
Helly's intersection theorem \cite{DGK}, a classical result from the Combinatorial Geometry of convex sets.
\begin{hellytheorem}\lbl{HT-2} Let $\Kc$ be a finite family of nonempty convex subsets of $Y$ lying in an affine subspace of $Y$ of dimension $m$. Suppose that every subfamily of $\Kc$ consisting of at most $m+1$ elements has a common point. Then there exists a point common to all of the family $\Kc$.
\end{hellytheorem}
\begin{lemma}\lbl{G-AB} Let $\ell\ge 0$. Suppose that the restriction $F|_{\Mc'}$ of $F$ to an arbitrary subset $\Mc'\subset\Mc$ consisting of at most $k_{\ell+1}$ points has a Lipschitz selection $f_{\Mc'}:\Mc'\to Y$ with  $\|f_{\Mc'}\|_{\Lip(\Mc',Y)}\le\lambda$. Then for all $x\in \Mc$ \smallskip
\par (a)\, $\GL(x)\ne\emp$\,;\medskip
\par (b)\, $\GL(x)\subset \GLO(y)+\lambda\,\rho(x,y)\,\BY$~ for all $y\in\Mc$\,, provided $\ell\ge 1$.
\end{lemma}
\par {\it Proof.} Thanks to \rf{BC-3}, \rf{AFF-4} and Helly's Theorem \reff{HT-2}, conclusion (a) will follow if we can show that
\bel{BC-7}
\Gamma(x,S_1)\,\capsm...\capsm\,\Gamma(x,S_{m+1})\ne\emp
\ee
for every $S_1,...,S_{m+1}\subset \Mc$ such that $\#S_i\le k_\ell$ (each $i$). (We note that, by \rf{AFF-4}, each set $\Gamma(x,S)$ is a subset of the affine space $\affspan(F(x))$ of dimension at most $m$. We also use the fact that there are only finitely many $S\subset \Mc$ because $\Mc$ is finite.)\medskip
\par However, $S_1\cup...\cup S_{m+1}\cup\{x\}\subset \Mc$ has cardinality at most
$$
(m+1)\cdot k_\ell+1\le k_{\ell+1}.
$$
\par The lemma's hypothesis therefore produces a function $\tff:S_1\cup...\cup S_{m+1}\cup\{x\}\to Y$ such that
$\tff(z)\in \FF(z)$ for all $z\in S_1\cup...\cup S_{m+1}\cup\{x\}$, and
$$
\|\tff(z)-\tff(w)\|\le \lambda\,\rho(z,w)~~~\text{for all}~~~z,w\in S_1\cup...\cup S_{m+1}\cup\{x\}.
$$
Then $\tff(x)$ belongs to $\Gamma(x, S_i)$ for $i=1,...,m+1$, proving \rf{BC-7} and thus also proving {\it(a)}. \medskip
\par To prove {\it(b)}, let $x,y \in \Mc$, and let $\xi\in \Gamma_\ell(x)$ with $\ell\ge 1$. We must show that there exists $\eta\in \Gamma_{\ell-1}(y)$ such that
$\|\xi-\eta\|\le\lambda\cdot \rho(x,y)$. To produce such an $\eta$, we proceed as follows.
\par Given a set $S\subset \Mc$ we introduce a set $\GH(x,y,\xi,S)$ consisting of all points $\eta\in Y$ such that there exists a mapping $f:S\cup\{x,y\}\to Y$ satisfying the following conditions:\smallskip
\par (i) $f(x)=\xi$, $f(y)=\eta$, and $f(z)\in F(z)$ for all $z\in S\cup\{x,y\}$;
\par (ii) For every $z,w\in S\cup\{x,y\}$ the following inequality
$$
\|f(z)-f(w)\|\le \lambda\,\rho(z,w)
$$
holds.
\smallskip
\par Clearly, $\GH(x,y,\xi,S)$ is a {\it convex} subset of $F(y)$. Let us show that
\bel{BC-8}
\bigcap_{\substack{S\subset \Mc\\\#S\le k_{\ell-1}}}
\GH(x,y,\xi,S)\,\ne\emp\,.
\ee
Thanks to Helly's Theorem \reff{HT-2}, \rf{BC-8} will follow if we can show that
\bel{BC-9}
\GH(x,y,\xi,S_1)\cap...\cap \GH(x,y,\xi,S_{\MPL})\,\ne\emp
\ee
for all $S_1,...,S_{\MPL}\subset \Mc$ with $\#S_i\le k_{\ell-1}$ (each $i$).
\par We set $\tS=S_1\cup...\cup S_{\MPL}\cup\{y\}$. Then $\tS\subset \Mc$ with
$$
\#\tS\le (\MPL)\cdot k_{\ell-1}+1\le k_\ell.
$$
Because $\xi\in \GL(x)\subset\Gamma(x,\tS)$ (see \rf{BC-3}), there exists $\tff:S_1\cup...\cup S_{\MPL}\cup\{x,y\}\to Y$ such that
$$
\tff(x)=\xi,~ \tff(z)\in \FF(z)~~~\text{for all}~~~
z\in S_1\cup...\cup S_{\MPL}\cup\{x,y\},
$$
and
$$
\|\tff(z)-\tff(w)\|\le \lambda\,\rho(z,w)~~~\text{for}~~~z,w\in S_1\cup...\cup S_{\MPL}\cup\{x,y\}.
$$
We then have $\tff(y)\in\GH(x,y,\xi,S_i)$ for $i=1,...,\MPL$, proving \rf{BC-9} and therefore also proving \rf{BC-8}.
\par Let
$$
\eta\in \bigcap_{\substack{S\subset \Mc\\\#S\le k_{\ell-1}}}
\GH(x,y,\xi,S)\,.
$$
Taking $S=\emp$, we obtain a function $\ff:\{x,y\}\to Y$ with $\ff(x)=\xi$, $\ff(y)=\eta$ and
$$
\|\ff(z)-\ff(w)\|\le \lambda\,\rho(z,w)~~~\text{for}~~~z,w\in \{x,y\}.
$$
Therefore,
\bel{BC-10}
\|\eta-\xi\|\le \lambda\,\rho(z,w).
\ee
Moreover, because $\GH(x,y,\xi,S)\subset\Gamma(y,S)$ for any $S\subset \Mc$ (see Definition \reff{GML-D}), we have
\bel{BC-11}
\eta\in \bigcap_{\substack{S\subset \Mc\\\#S\le k_{\ell-1}}}
\Gamma(y,S)=\GLO(y)\,.
\ee
\par Our results \rf{BC-10}, \rf{BC-11} complete the proof of {\it(b)}.\bx

\bigskip
\indent\par {\bf 3.2 Statement of the Finiteness Theorem for bounded Nagata dimension.}
\addtocontents{toc}{~~~~3.2 Statement of the Finiteness Theorem for bounded Nagata dimension.\hfill \thepage\par}\medskip
\par We place ourselves in the following setting. \smallskip
\par \textbullet~We fix a positive integer $m$.\smallskip
\par \textbullet~$(X,d)$ is a finite metric space satisfying Nagata $(\CNG,\cng)$ (see Definition \reff{NG-C}). \smallskip
\par \textbullet~$\BS$ is a Banach space. We write $\|\cdot\|$ for the norm in $Y$, and $\|\cdot\|_{Y^*}$ for the norm in the dual space $Y^*$. We write $\ip{e,y}$ to denote the natural pairing between vectors $y\in Y$ and dual vectors $e\in Y^*$.\smallskip
\par \textbullet~For each $x\in X$ we are given a convex set
$$
\FF(x)\subset \Aff_{\FF}(x)\subset Y,
$$
where
$$
\Aff_F(x)~~~\text{is an affine subspace of}~~Y,~~\text{of dimension at most}~~m.
$$
Say, $\Aff_F(x)$ is a translate of the vector subspace $\Vect_{\FF}(x)\subset Y$.\smallskip
\par \textbullet~We make the following assumption for a large enough $\ks$ determined by $m$. \medskip
\begin{fnassumption}\lbl{FNA} Given $S\subset X$ with $\#S\le \ks$, there exists $f^S:S\to Y$ with Lipschitz seminorm at most $1$, such that $f^S(x)\in F(x)$ for all $x\in S$.
\end{fnassumption}
\smallskip
\par The above assumption implies the existence of a Lipschitz selection with a controlled Lipschitz seminorm. More precisely, we have the following result.
\begin{theorem}\lbl{PFT-FT}(Finiteness Theorem for bounded Nagata dimension) Let $(X,\dt)$ be a finite metric space
satisfying Nagata $(\CNG,\cng)$.
\par Given $m\in\N$ there exist a constant $\ks\in\N$ depending only on $m$, and a constant $\gnd>0$ depending only on $m$, $\cng$, $\CNG$, for which the following holds: Let $Y$ be a Banach space. For each $x\in X$, let $F(x)\subset Y$ be a convex set of (affine) dimension at most $m$.\smallskip
\par Suppose that for each $S\subset X$ with $\#S\le \ks$ there exists $f^S:S\to Y$ with Lipschitz seminorm at most $1$, such that $f^S(x)\in F(x)$ for all $x\in S$.
\par Then there exists $f:X\to Y$ with Lipschitz seminorm at most $\gnd$, such that $f(x)\in F(x)$ for all $x\in X$.
\end{theorem}
\par By applying Theorem \reff{PFT-FT} to the metric space $(X,\lambda\,d)$ we establish Theorem \reff{NG-FP}.
\smallskip
\par We place ourselves in the above setting until the end of the proof of Theorem \reff{PFT-FT} in the end of  Section 4.9.
\bigskip\medskip
\indent\par {\bf 3.3 Labels and bases.}
\addtocontents{toc}{~~~~3.3 Labels and bases.\hfill \thepage\par\VSU}\medskip
\par A ``label'' is a finite sequence $\Ac=(e_1,e_2,...,e_s)$ of functionals $e_a\in Y^*$, $a=1,...,s$, with  $s\le m$. Here, $m$ is as in the hypothesis of Theorem \reff{PFT-FT}.
\par We write $\#\Ac$ to denote the number $s$ of functionals $e_a$ appearing in $\Ac$. We allow the case $\#\Ac=0$, in which case $\Ac$ is the empty sequence $\Ac=(~\,)$.
\par Let $\Gamma\subset Y$ be a convex set, let $\Ac=(e_1,e_2,...,e_s)$ be a label, and let $r, C_B$ be positive real numbers. Finally, let $\zeta\in Y$.
\begin{definition}\lbl{LB-DB} {\em An {\it $(\Ac,r,C_B)$-basis} for $\Gamma$ at $\zeta$ is a sequence of $s$ vectors $v_1,...,v_s\in Y$, with the following properties:\medskip
\par (B0)~ $\zeta\in\Gamma$\,.
\medskip
\par (B1)~ $\ip{e_a,v_b}=\dl_{ab}$ (Kronecker delta) for $a,b=1,...,s$\,.
\medskip
\par (B2)~ $\|v_a\|\le C_B$ and $\|e_a\|_{Y^*}\le C_B$ for $a=1,...,s$\,.
\medskip
\par (B3)~ $\zeta+\frac{r}{C_B}v_a$~ and~ $\zeta-\frac{r}{C_B}v_a$ belong to $\Gamma$ for $a=1,...,s$\,.
}
\end{definition}
\smallskip
\par If $s\ge 1$, then of course (B3) implies (B0).
\par Let us note several elementary properties of $(\Ac,r,C_B)$-bases.
\begin{remark}\lbl{LBL-R}{\em (i) If $s=0$ then (B1), (B2), (B3) hold vacuously, so the assertion that $\Gamma$ has an $\left((~\,),r,C_B\right)$-basis at $\zeta$ means simply that $\zeta\in\Gamma$;
\smallskip
\par (ii) If $r'\le r$ and $C'_B\ge C_B$, then any $(\Ac,r,C_B)$-basis for $\Gamma$ at $\zeta$ is also an
$(\Ac,r',C'_B)$-basis for $\Gamma$ at $\zeta$;\smallskip
\par (iii) If $K\ge 1$, then any $(\Ac,r,C_B)$-basis for $\Gamma$ at $\zeta$ is also an $(\Ac,Kr,KC_B)$-basis for $\Gamma$ at $\zeta$;\smallskip
\par (iv) If $\Gamma\subset \Gamma'$, then every  $(\Ac,r,C_B)$-basis for $\Gamma$ at $\zeta$ is also an $(\Ac,r,C_B)$-basis for $\Gamma'$ at $\zeta$.}
\end{remark}
\begin{lemma}\lbl{LB-L1} (``Adding a vector'')\, Suppose $\Gamma\subset Y$ (convex) has an $(\Ac,r,C_B)$-basis at $\xi$, where $\Ac=(e_1,e_2,...,e_s)$ and $s\le m-1$.
\par Let $\eta\in\Gm$, and suppose that
$$
\|\eta-\xi\|\ge r
$$
and
$$
\ip{e_a,\eta-\xi}=0~~~\text{for}~~~a=1,...,s.
$$
\par Then there exist $\zeta\in\Gm$ and $e_{s+1}\in Y^*$ with the following properties:\medskip
\par \textbullet~ $\|\zeta-\xi\|=\tfrac12 r$.
\medskip
\par \textbullet~ $\ip{e_a,\zeta-\xi}=0$~ for $a=1,...,s$ (not necessarily for $a=s+1$).
\medskip
\par \textbullet~ $\Gm$ has an
$(\Ac^+,r,C'_B)$-basis at $\zeta$, where $\Ac^+=(e_1,...,e_s,e_{s+1})$ and $C'_B$ is determined by $C_B$ and $m$.
\end{lemma}
\par {\it Proof.} In this proof, we write $C$ to denote a positive constant determined by $C_B$ and $m$. This symbol may denote different constants in different occurrences.
\par Let $(v_1,...,v_s)$ be an $(\Ac,r,C_B)$-basis for $\Gm$ at $\xi$. Thus, $\xi\in\Gamma$,
\bel{LB-2}
\ip{e_a,v_b}=\dl_{ab}~~~\text{for}~~~a,b=1,...,s,
\ee
\bel{LB-3}
\|e_a\|_{Y^*}\le C_B,~~\|v_a\|\le C_B~~~\text{for}~~~ a=1,...,s,
\ee
\bel{LB-4}
\xi+\tfrac{r}{C_B}v_a,~~\xi-\tfrac{r}{C_B}\,v_a\in\Gamma~~~
\text{for}~~~a=1,...,s\,.
\ee
\smallskip
\par Let
$$
\zeta=\tau\,\eta+(1-\tau)\,\xi~~~\text{with}~~~
\tau=\tfrac12\,r\,\|\xi-\eta\|^{-1}\in(0,\tfrac12].
$$
(Note that, by the lemma's hypothesis, $\|\xi-\eta\|$ is non-zero so that $\tau$ and $\zeta$ are well defined.)
\par Our hypotheses on $\xi$ and $\eta$ tell us that
\bel{LB-5}
\zeta\in\Gm,~~\|\zeta-\xi\|=\tfrac12 r,~~ \ip{e_a,\zeta-\xi}=0
~~~\text{for}~~~a=1,...,s.
\ee
Because $\eta\in\Gm$, $\Gm$ is convex, and $\tau\in(0,\tfrac12]$,~ \rf{LB-4} implies
\bel{LB-6}
\zeta+\tfrac12\tfrac{r}{C_B}\,v_a,~ \zeta-\tfrac12\tfrac{r}{C_B}\,v_a\in\Gamma~~~\text{for}~~~ a=1,...,s\,.
\ee
\par Let
\bel{LB-7}
v_{s+1}=\frac{\zeta-\xi}{\|\zeta-\xi\|}\,.
\ee
(The denominator is nonzero, by \rf{LB-5}.) Then
$$
\zeta+\|\zeta-\xi\|\,v_{s+1}=\zeta+(\zeta-\xi)=2\zeta-\xi=
2\tau\eta+(1-2\tau)\xi\in\Gm
$$
because $\xi,\eta\in\Gm$, $\Gamma$ is convex and $\tau\in(0,\tfrac12]$.
\par Also,
$$
\zeta-\|\zeta-\xi\|\,v_{s+1}=\zeta-(\zeta-\xi)=\xi\in\Gm
\,.
$$
Recall that $\|\zeta-\xi\|=\frac12\,r$, hence the above remarks and \rf{LB-6} together yield
\bel{LB-8}
\zeta+\tfrac{r}{C}v_a,~\zeta-\tfrac{r}{C}v_a\in\Gm~~~
\text{for}~~~a=1,...,s+1,
~~~\text{and}~~~C=2\max\{C_B,1\}.
\ee
\par Also, because $\ip{e_a,\zeta-\xi}=0$ for $a=1,...,s$, (see \rf{LB-5}), the definition of $v_{s+1}$, together with \rf{LB-2}, tells us that
\bel{LB-9}
\ip{e_a,v_b}=\dl_{ab}~~~\text{for}~~~a=1,...,s~~~\text{and}
~~~b=1,...,s+1.
\ee
\par We prepare to define a functional $e_{s+1}\in Y^*$. To do so, we first prove the estimate
\bel{LB-10}
\sum_{a=1}^{s+1}\,|\lambda_a|\le C\,\left\|\sum_{a=1}^{s+1}\lambda_a\,v_a\right\|
~~~\text{for all}~~~\lambda_1,...,\lambda_{s+1}\in\R.
\ee
\par To see this, we first note that for any $b=1,...,s$,\, \rf{LB-9} yields the estimate
\bel{LB-11}
|\lambda_b|=\left|\ip{e_b,\sum_{a=1}^{s+1}\lambda_a\,v_a}
\right|\le
\|e_b\|_{Y^*}\cdot \left\|\sum_{a=1}^{s+1}\lambda_a\,v_a\right\|
\le C_B\,\left\|\sum_{a=1}^{s+1}\lambda_a\,v_a\right\|\,.
\ee
Consequently,
\be
|\lambda_{s+1}|&=&\left\|\lambda_{s+1}v_{s+1}\right\|
\le \left\|\sum_{a=1}^{s+1}\lambda_a\,v_a\right\|+
\sum_{a=1}^{s}|\lambda_a|\,\|v_a\|
\nn\\
&\le&
\left\|\sum_{a=1}^{s+1}\lambda_a\,v_a\right\|+
C_B\sum_{a=1}^{s}|\lambda_a|
\le (1+m\,C_B^2)\left\|\sum_{a=1}^{s+1}\lambda_a\,v_a
\right\|\,.\nn
\ee
(Recall that $0\le s\le m$.) Together with \rf{LB-11}, this completes the proof of \rf{LB-10}.\medskip
\par By \rf{LB-10} and the Hahn-Banach theorem, the linear functional
$$
\sum_{a=1}^{s+1}\lambda_a\,v_a\to\lambda_{s+1}
$$
on the span of $v_1,...,v_{s+1}$ extends to a linear functional $e_{s+1}\in Y^*$, with
\bel{LB-12}
\|e_{s+1}\|_{Y^*}\le C
\ee
and
\bel{LB-13}
\ip{e_{s+1},v_a}=\dl_{s+1,a}~~~\text{for}~~~a=1,...,s+1.
\ee
\par From \rf{LB-3}, \rf{LB-5}, \rf{LB-7}, \rf{LB-9}, \rf{LB-12}, \rf{LB-13} we have
\bel{LB-14}
\zeta\in\Gm,
\ee
\bel{LB-15}
\|e_a\|_{Y^*},~ \|v_a\|\le C~~~\text{for}~~~ a=1,...,s+1,
\ee
\bel{LB-16}
\ip{e_a,v_b}=\dl_{ab}~~~\text{for}~~~ a,b=1,...,s+1.
\ee
From  \rf{LB-8}, \rf{LB-14}, \rf{LB-15}, \rf{LB-16}, we see that $v_1,...,v_{s+1}$ form an
$((e_1,...,e_{s+1}),r,C)$-basis for $\Gm$ at $\zeta$.
\par Together with \rf{LB-5}, this completes the proof of Lemma \reff{LB-L1}.\bx
\begin{lemma}\lbl{LB-L2} (``Transporting a Basis'')\,
Given $m\in\N$ and $C_B>0$ there exists a constant $\ve_0\in(0,1]$ depending only on $m$, $C_B$, for which the following holds:
\par Suppose $\Gamma\subset Y$ (convex) has an $(\Ac,r,C_B)$-basis at $\xi_0$, where $\Ac=(e_1,e_2,...,e_s)$ and $s\le m$. Suppose $\Gm'\subset Y$ (convex) satisfies:\medskip
\par (*)~ Given any $\xi\in\Gm$ there exists $\eta\in\Gm'$ such that $\|\xi-\eta\|\le\ve_0r$.
\medskip
\par Then there exists $\eta_0\in\Gm'$ with the following properties:\medskip
\par \textbullet~ $\|\eta_0-\xi_0\|\le C\,r\,.$
\medskip
\par \textbullet~ $\ip{e_a,\eta_0-\xi_0}=0$~ for $a=1,...,s$.
\medskip
\par \textbullet~ $\Gm'$ has an $(\Ac,r,C)$-basis at $\eta_0$.
\medskip
\par Here, C is determined by $C_B$ and $m$.
\end{lemma}
\par {\it Proof.} In the trivial case $s=0$ (see Remark \reff{LBL-R} (i)), Lemma \reff{LB-L2} holds because it simply asserts that there exists $\eta_0\in\Gm'$ such that $\|\eta_0-\xi_0\|\le C\,r$, which is immediate from (*). We suppose $s\ge 1$.
\par We take
\bel{EP0-A}
\ve_0~~\text{to be less than a small enough positive constant determined by}~~C_B~~\text{and}~~m.
\ee
We can take $\ve_0$ to be, say, $\frac12$ times that small positive constant.
\smallskip
\par We write $c_1,c_2,c_3,C$ to denote positive constants determined by $C_B$ and $m$. These symbols may denote different constants in different occurrences.
\par Let $(v_1,...,v_s)$ be an $(\Ac,r,C_B)$-basis for $\Gm$ at $\xi_0$. Thus, $\xi_0\in\Gamma$,
\bel{LB-18}
\ip{e_a,v_b}=\dl_{ab}~~~\text{for}~~~a,b=1,...,s,
\ee
\bel{LB-19}
\|e_a\|_{Y^*}\le C_B,~~\|v_a\|\le C_B~~~\text{for}~~~ a=1,...,s,
\ee
and
\bel{LB-20}
\xi_0+c_1\sigma r\,v_a\in\Gamma~~~
\text{for}~~~a=1,...,s,~~\sigma=\pm 1~~~\text{and}~~~c_1=1/C_B\,.
\ee
\smallskip
\par Applying our hypothesis (*) to the vectors in \rf{LB-20}, we obtain vectors
$$
\zeta_{a,\sigma}\in Y~~~~(a=1,...,s,~~\sigma=\pm 1)
$$
such that
\bel{LB-21}
\xi_0+c_1\sigma r\,v_a+\zeta_{a,\sigma}\in\Gamma'~~~
\text{for}~~~a=1,...,s,~~\sigma=\pm 1,
\ee
and
\bel{LB-22}
\|\zeta_{a,\sigma}\|\le \ve_0\,r~~~
\text{for}~~~a=1,...,s,~~\sigma=\pm 1\,.
\ee
\par We define vectors
\bel{LB-23}
\eta_{00}=\frac{1}{2s}\,\sum_{a=1}^s\,\sum_{\sigma=\pm 1}
\,(\xi_0+c_1\sigma r v_a+\zeta_{a,\sigma})=
\xi_0+\frac{1}{2s}\,\sum_{a=1}^s\,\sum_{\sigma=\pm 1}
\zeta_{a,\sigma}
\ee
and
\bel{LB-24}
\tv_a=\frac{[\xi_0+c_1rv_a+\zeta_{a,1}]
-[\xi_0-c_1 r v_a+\zeta_{a,-1}]}{2c_1r}=
v_a+\left(\frac{\zeta_{a,1}-\zeta_{a,-1}}{2c_1r}\right)
\ee
for $ a=1,...,s$.
\par From \rf{LB-21} and the first equality in \rf{LB-23}, we have $\eta_{00}\in \Gm'$. From \rf{LB-22} and the second equality in \rf{LB-23}, we have
\bel{LB-26}
\|\eta_{00}-\xi_0\|\le \ve_0 r.
\ee
From \rf{LB-22} and the second equality in \rf{LB-24}, we have
\bel{LB-27}
\|\tv_a-v_a\|\le C\,\ve_0~~~\text{for}~~~a=1,...,s.
\ee
Also, for $b=1,...,s$ and $\hsg=\pm 1$, the first equalities in \rf{LB-23}, \rf{LB-24} give
$$
\eta_{00}+\frac{1}{s}c_1r\hsg\,\tv_b
=\frac{1}{2s}\,\sum_{a=1}^s\,\sum_{\sigma=\pm 1}
\,(\xi_0+c_1\sigma r v_a+\zeta_{a,\sigma})+
\frac{\hsg}{2s}
[(\xi_0+c_1rv_b+\zeta_{b,1})-
(\xi_0-c_1 r v_b+\zeta_{b,-1})],
$$
which exhibits $\eta_{00}+\frac{1}{s}c_1r\hsg\,\tv_b$ as a convex combination of the vectors in \rf{LB-21}. Consequently,
$$
\eta_{00}+c_2r\,\tv_b,~~\eta_{00}-c_2r\,\tv_b\in\Gamma'
~~~\text{for}~~~b=1,...,s,
$$
which implies that
\bel{LB-28}
\eta_{00}+c_2r\sum_{a=1}^s\,\tau_a\tv_a\in\Gm'
~~~\text{for any}~~~\tau_1,...,\tau_s\in\R~~~\text{with}~~  \sum_{a=1}^s\,|\tau_a|\le 1.
\ee
Here we use the following remark on convex sets: Suppose $\xi+\eta_i,\xi-\eta_i$, $(i=1,...,I)$ belong to a convex set $\Gamma$. Then
$$
\xi+\sum_{i=1}^I\,\tau_i\eta_i\in\Gamma~~~\text{for all}~~~\tau_1,...,\tau_I\in\R~~~\text{with}~~
\sum_{i=1}^I\,|\tau_i|\le 1.
$$
\smallskip
\par From \rf{LB-18}, \rf{LB-19}, \rf{LB-27}, we have
\bel{REVISION-1}
|\ip{e_a,\tv_b}-\dl_{ab}|\le C\ve_0 ~~~\text{for}~~~a,b=1,...,s.
\ee
\par We let $A$ denote the $s\times s$ matrix $A=(\ip{e_a,\tv_b})_{a,b=1}^s$. Let $I=(\delta_{ab})_{a,b=1}^s$ be the identity matrix. Given an $s\times s$ matrix $T$, we let $\|T\|_{op}$ denote the operator norm of $T$ as an operator from $\ell^2_s$ into $\ell^2_s$. Clearly, $\|T\|_{op}$ is equivalent (with constants depending only on $s$) to $\max\{|t_{ab}|:1\le a,b\le s\}$ provided $T=(t_{ab})_{a,b=1}^s$.
\par Hence, by \rf{REVISION-1},
\bel{MN-S}
\|A-I\|_{op}\le C\ve_0\,.
\ee
\par We recall the standard fact from matrix algebra which states that an $s\times s$ matrix $T$ is invertible and the inequality $\|T^{-1}-I\|_{op}\le \|T-I\|_{op}/(1-\|T-I\|_{op})$ is satisfied provided $\|T-I\|_{op}<1$. Therefore, by \rf{MN-S}, for $\ve_0$ small enough, the matrix $A$ is invertible, and the following inequality
\bel{A-KL}
\|A^{-1}-I\|_{op}\le 2\,\|A-I\|_{op}
\ee
holds.
\par Let $(A^{\bf T})^{-1}=(M_{gb})_{g,b=1,...,s}$ where $A^{\bf T}$ denotes the transpose of $A$. Then
\bel{LB-29}
\ip{e_a,\sum_{b=1}^s\,M_{gb}\,\tv_b}=\dl_{ag}
~~~\text{for}~~~a,g=1,...,s.
\ee
Moreover, by \rf{MN-S} and \rf{A-KL},
\bel{LB-30}
|M_{gb}-\dl_{gb}|\le C\,\ve_0
~~~\text{for}~~~g,b=1,...,s.
\ee
\par We set
\bel{LB-31}
\hv_g=\sum_{b=1}^s\,M_{gb}\,\tv_b
~~~\text{for}~~~g=1,...,s.
\ee
\par Then \rf{LB-19}, \rf{LB-27}, \rf{LB-30}, \rf{LB-31} yield
\bel{LB-32}
\|\hv_g\|\le C~~~\text{for}~~~g=1,...,s,
\ee
while \rf{LB-29}, \rf{LB-31} give
\bel{LB-33}
\ip{e_a,\hv_g}=\dl_{ag}~~~\text{for}~~~a,g=1,...,s.
\ee
\par Moreover, \rf{LB-28}, \rf{LB-30}, \rf{LB-31} together imply that
\bel{LB-34}
\eta_{00}+c_3r\sum_{g=1}^s\,\tau_g\hv_g\in\Gm'
~~~\text{for all}~~~\tau_1,...,\tau_s~~~\text{such that each}~~  |\tau_g|\le 1.
\ee
\par To see this, we simply write the linear combination
of the $\hv_g$ in \rf{LB-34} as a linear combination of the $\tv_b$ using \rf{LB-31}, and then recall \rf{LB-28}.
\par From \rf{LB-19}, \rf{LB-26} we have
\bel{LB-35}
|\ip{e_a,\eta_{00}-\xi_0}|\le C\ve_0\, r ~~~\text{for}~~~a=1,...,s.
\ee
\par We set
\bel{LB-36}
\eta_0=\eta_{00}-
\sum_{g=1}^s\,\ip{e_g,\eta_{00}-\xi_0}\,\hv_g,
\ee
so that by \rf{LB-33},
\bel{LB-37}
\ip{e_a,\eta_0-\xi_0}=\ip{e_a,\eta_{00}-\xi_0}-
\sum_{g=1}^s\,\ip{e_g,\eta_{00}-\xi_0}\ip{e_a,\hv_g}=0
~~~\text{for}~~~a=1,...,s.
\ee
Also, by \rf{LB-26}, \rf{LB-32}, \rf{LB-35},
\bel{LB-38}
\|\eta_0-\xi_0\|\le \|\eta_{00}-\xi_0\|+
\sum_{g=1}^s\,|\ip{e_g,\eta_{00}-\xi_0}|\cdot\|\hv_g\|
\le C\ve_0 r\,.
\ee
\par From \rf{LB-35} and our {\it small $\ve_0$ assumption} \rf{EP0-A}, we have
$$
|\ip{e_a,\eta_{00}-\xi_0}|\le \tfrac12 c_3 r
~~~\text{for}~~~a=1,...,s,
$$
with $c_3$ as in \rf{LB-34}.
\par Therefore \rf{LB-34} and \rf{LB-36} tell us that
$$
\eta_0+c_3r
\sum_{g=1}^s\,\tau_g\hv_g\in\Gm' ~~~\text{for any}~~~
\tau_1,...,\tau_s~~~\text{such that}~~~
|\tau_g|\le\tfrac12~~~\text{for each}~~g.
$$
In particular,
\bel{LB-39}
\eta_0\in\Gm'
\ee
and
$$
\eta_0+\tfrac12 c_3r\,\hv_g,~\eta_0-\tfrac12 c_3r\,\hv_g\in\Gm'~~~\text{for}~~~g=1,...,s.
$$
\par Also, recalling \rf{LB-19}, \rf{LB-32}, \rf{LB-33}, we note that
$$
\|e_a\|_{Y^*},~\|\hv_a\|\le C~~~\text{for}~~~a=1,...,s
$$
and
\bel{LB-42}
\ip{e_a,\hv_g}=\dl_{ag}~~~\text{for}~~~a,g=1,...,s.
\ee
\par Our results \rf{LB-39},...,\rf{LB-42} tell us that $\hv_1,...,\hv_s$ form an $(\Ac,r,C)$-basis for $\Gm'$ at $\eta_0$, with $\Ac=(e_1,...,e_s)$. That's the third bullet point in the statement of Lemma \reff{LB-L2}. The other two bullet points are immediate from our results  \rf{LB-38} and \rf{LB-37}.
\par The proof of Lemma \reff{LB-L2} is complete.\bx
\SECT{4. The Main Lemma}{4}
\addtocontents{toc}{4. The Main Lemma\hfill\thepage \par\VST}

\indent\par {\bf 4.1 Statement of the Main Lemma.}
\addtocontents{toc}{~~~~4.1 Statement of the Main Lemma.\hfill \thepage\par}\medskip
\par Recall that $(X,d)$ is a (finite) metric space
satisfying Nagata $(\CNG,\cng)$.
\par For any label $\Ac=(e_1,...,e_s)$, we define
\bel{ELL-2}
\ell(\Ac)=2+3\cdot(m-\#\Ac)=2+3\cdot(m-s).
\ee
\par Note that
$$
\ell(\Ac)\ge\ell(\Ac^+)+3~~~\text{whenever}~~~ \#\Ac^+>\#\Ac.
$$
\par We now choose the constant $\ks$ in our Finiteness Assumption \reff{FNA}. We take
\bel{KS-DF}
\ks=k_{\ell^{\#}+1}=(m+2)^{\ell^{\#}+1}
\ee
as in equation \rf{KEL}, with
\bel{LS-DF}
\ell^{\#}=2+3m.
\ee
\par In this setting we define a family $\GL(x)$ of basic convex sets as in Section 3.1. More specifically, let $(\Mc,\rho)=(X,\dt)$, $\lambda=1$ and let
$F:X\to\Conv_m(Y)$ be the set-valued mapping from Theorem \reff{PFT-FT}. We apply Definition \reff{GML-D} and formulae \rf{KEL}, \rf{BC-3} to these objects and obtain a family
$$
\{\GL(x):x\in X,\ell=0,1,...\}
$$
of convex subsets of $Y$.
\par Finally, we apply Lemma \reff{G-AB} to the setting of this section. The Finiteness Assumption \reff{FNA} enables us to replace the hypothesis of this lemma with the requirement $\ks\ge k_{\ell+1}$, which together with definition \rf{ELL-2} of $\ell(\Ac)$ leads us to the following statement.
\begin{lemma}\lbl{BCS-L1P} Let $\Ac$ be a label. Then
\smallskip
\par (A) $\GL(x)\ne\emp$~ for any $x\in X$ and any $\ell\le\ell(\Ac)$.
\medskip
\par (B) Let $1\le\ell\le\ell(\Ac)$, let $x,y\in X$, and let\, $\xi\in\GL(x)$. Then there exists $\eta\in\GLO(y)$ such that
$$
\|\xi-\eta\|\le d(x,y).
$$
\end{lemma}
\par In Sections 4.2-4.9 we will prove the following result.
\begin{mainlemma}\lbl{SML-ML} Let $x_0\in X$, $\xi_0\in Y$, $r_0>0$, $C_B\ge 1$ be given, and let $\Ac$ be a label.
\par Suppose that $\Gm_{\ell(\Ac)}(x_0)$ has an $(\Ac,\ve^{-1}r_0,C_B)$-basis at $\xi_0$, where $\ve>0$ is less than a small enough constant $\emin>0$ determined by $m$, $C_B$, $\cng$, $\CNG$.
\par Then there exists $f:\BXR\to Y$ with the following properties:\medskip
\bel{SML-A1}
\|\ff(z)-\ff(w)\|\le C(\ve)\,d(z,w)~~~\text{for all}~~~z,w\in \BXR,
\ee
\bel{SML-A2}
\|\ff(z)-\xi_0\|\le C(\ve)\,r_0~~~\text{for all}~~~z\in \BXR,
\ee
\bel{SML-A3}
\ff(z)\in\Gm_0(z)~~~\text{for all}~~~z\in \BXR.
\ee
Here $C(\ve)$ is determined by $\ve$, $m$, $C_B$, $\cng$, $\CNG$.
\end{mainlemma}
\medskip
\par We will prove the Main Lemma \reff{SML-ML} by downward induction on $\#\Ac$, starting with the case ${\#\Ac=m}$, and ending with the case $\#\Ac=0$.

\indent\par {\bf 4.2 Proof of the Main Lemma in the base case.}
\addtocontents{toc}{~~~~4.2 Proof of the Main Lemma in the base case.\hfill \thepage\par}\medskip
\par In this section, we assume the hypothesis of the Main Lemma \reff{SML-ML} in the base case ${\Ac=(e_1,...,e_m)}$. Thus, in this case $\#\Ac=m$ and $\ell(\Ac)=2$, (see \rf{ELL-2}).
\par We recall that for each $x\in X$ we have $\GL(x)\subset \FF(x)\subset \Aff_{\FF}(x)$ (all $\ell\ge0$), where $\Aff_{\FF}(x)$ is a translate of the vector space $\Vect_{\FF}(x)$ of dimension $\le m$. We write $C$ to denote a positive constant determined by $m$, $C_B$, $\cng$, $\CNG$. This symbol may denote different constants in different occurrences.
\begin{lemma}\lbl{SML-L1} For each $z\in \BXR$, there exists
\bel{SML-1}
\eta^z\in\Gm_{1}(z)
\ee
such that
\bel{SML-2}
\|\,\eta^z-\xi_0\,\|\le C\,\ve^{-1}r_0,
\ee
\bel{SML-3}
\ip{e_a,\eta^z-\xi_0}=0~~~\text{for}~~~a=1,...,m,
\ee
\bel{SML-4}
\Gm_{1}(z)~~~\text{has an}~~~(\Ac,\ve^{-1}r_0,C)\text{-basis at}~~~\eta^z.
\ee
\end{lemma}
\par {\it Proof.} We apply Lemma \reff{LB-L2}, taking $\Gm$ to be $\Gm_{2}(x_0)$, $\Gm'$ to be
$\Gm_{1}(z)$, and $r$ to be $\ve^{-1}r_0$. To apply that lemma, we must check the key hypothesis (*), which asserts in the present case that
\bel{SML-5}
\text{Given}~~\xi\in \Gm_{2}(x_0)~~\text{there exists}~~\eta\in\Gm_{1}(z)~~\text{such that}~~\|\xi-\eta\|\le\ve_0\cdot(\ve^{-1}r_0),
\ee
where $\ve_0$ is a small enough constant determined by $C_B$ and $m$.\medskip
\par To check \rf{SML-5}, we recall Lemma \reff{BCS-L1P} (B). Given $\xi\in\Gm_{2}(x_0)$ there exists $\eta\in\Gm_{1}(z)$ such that
$$
\|\xi-\eta\|\le d(z,x_0)\le r_0~~(\text{because}~ z\in\BXR)<\ve_0\cdot(\ve^{-1}r_0);
$$
here, the last inequality holds thanks to our assumption that $\ve$ is less than a small enough constant determined by $m$, $C_B$, $\cng$, $\CNG$.
\par Thus, \rf{SML-5} holds, and we may apply Lemma \reff{LB-L2}. That lemma   provides a vector $\eta^z$ satisfying \rf{SML-1},...,\rf{SML-4}, completing the proof of Lemma \reff{SML-L1}.\bx \smallskip
\par For each $z\in\BXR$, we fix a vector $\eta^z$ as in Lemma \reff{SML-L1}. Repeating the idea of the proof of Lemma \reff{SML-L1}, we establish the following result.
\begin{lemma}\lbl{SML-L2} Given $z,w\in \BXR$, there exists a vector
\bel{SML-6}
\eta^{z,w}\in\Gm_{0}(w)
\ee
such that
\bel{SML-7}
\|\eta^{z,w}-\eta^{z}\|\le C\,\ve^{-1}d(z,w)
\ee
and
\bel{SML-8}
\ip{e_a,\eta^{z,w}-\eta^{z}}=0~~~\text{for}~~~a=1,...,m.
\ee
\end{lemma}
\par {\it Proof.} If $z=w$, we can just take $\eta^{z,w}=\eta^{z}$. Suppose $z\ne w$. Because $z,w\in\BXR$, we have $0<d(z,w)\le 2r_0$. Therefore, \rf{SML-4} and Remark \reff{LBL-R} (ii) tell us that
\bel{SML-9}
\Gm_{1}(z)~~~\text{has an}~~~(\Ac,\tfrac12\ve^{-1}d(z,w), C)\text{-basis at}~~\eta^z.
\ee
\par We prepare to apply Lemma \reff{LB-L2}, this time taking
$$
\Gm=\Gm_{1}(z),~~~\Gm'=\Gm_{0}(w), ~~~r=\tfrac12 \ve^{-1}d(z,w).
$$
We must verify the key hypothesis (*), which asserts in the present case that:
\medskip
\par Given any $\xi\in\Gm_{1}(z)$ there exists $\eta\in\Gm_{0}(w)$ such that
\bel{SML-10}
\|\xi-\eta\|\le\ve_0\cdot(\tfrac12 \ve^{-1}d(z,w)),
\ee
where $\ve_0$ arises from the constant $C$ in \rf{SML-9}
as in Lemma \reff{LB-L2}. In particular, $\ve_0$ depends only on $m$ and $C_B$. Therefore, our assumption that $\ve$ is less than a small enough constant determined by $m$, $C_B$, $\cng$, $\CNG$ tells us that
$$
d(z,w)<\ve_0\cdot(\tfrac12 \ve^{-1}d(z,w)).
$$
\par Consequently, Lemma \reff{BCS-L1P} (B) produces for each $\xi\in \Gm_{1}(z)$ an $\eta\in \Gm_{0}(w)$ such that
$$
\|\xi-\eta\|\le d(z,w)<\ve_0\cdot(\tfrac12 \ve^{-1}d(z,w)),
$$
which proves \rf{SML-10}.
\par Therefore, we may apply Lemma \reff{LB-L2}. That lemma provides a vector $\eta^{z,w}$ satisfying  \rf{SML-6}, \rf{SML-7}, \rf{SML-8}, and additional properties that we don't need here.
\par The proof of Lemma \reff{SML-L2} is complete.\bx
\begin{lemma}\lbl{SML-L3} Let $w\in\BXR$. Then any vector $v\in\Vect_{\FF}(w)$ satisfying $\ip{e_a,v}=0$ for $a=1,...,m$ must be the zero vector.
\end{lemma}
\par {\it Proof.} Applying \rf{SML-4}, we obtain an $(\Ac,\ve^{-1}r_0,C)$-basis $(v_1,...,v_m)$ for $\Gm_{1}(w)$ at $\eta^w$. From the definition of an $(\Ac,\ve^{-1}r_0,C)$-basis, see Definition \reff{LB-DB}, we have
\bel{SML-11}
\ip{e_a,v_b}=\dl_{ab}~~~\text{for}~~~a,b=1,...,m,
\ee
and
$$
\eta^w+\tfrac{1}{C}\,\ve^{-1}r_0v_a,~\eta^w-
\tfrac{1}{C}\,\ve^{-1}r_0v_a\in
\Gm_{1}(w)\subset \FF(w)\subset\Aff_{\FF}(w)
~~~\text{for}~~~a=1,...,m,
$$
from which we deduce that
\bel{SML-12}
v_a\in \Vect_{\FF}(w)~~~\text{for}~~~a=1,...,m.
\ee
From \rf{SML-11}, \rf{SML-12} we see that
$$
v_1,...,v_m\in\Vect_{\FF}(w)
$$
are linearly independent. However, $\Vect_{\FF}(w)$ has dimension at most $m$. Therefore, $v_1,...,v_m$ form a basis for $\Vect_{\FF}(w)$. Lemma \reff{SML-L3} now follows at once from \rf{SML-11}.\bx
\medskip
\par Now let $z,w\in\BXR$. From Lemmas \reff{SML-L1} and \reff{SML-L2} we have
$$
\eta^w,\, \eta^{z,w}\in \Gm_{0}(w)\subset \FF(w)\subset\Aff_{\FF}(w),
$$
and consequently
\bel{SML-13}
\eta^w-\eta^{z,w}\in\Vect_{\FF}(w).
\ee
\par On the other hand, \rf{SML-3} and \rf{SML-8} tell us that
$$
\ip{e_a,\eta^w-\xi_0}=0,~~\ip{e_a,\eta^z-\xi_0}=0,~~
\ip{e_a,\eta^z-\eta^{z,w}}=0~~~\text{for}~~~a=1,...,m.
$$
Therefore,
\bel{SML-14}
\ip{e_a,\eta^w-\eta^{z,w}}=0~~~\text{for}~~~a=1,...,m.
\ee
\par From \rf{SML-13}, \rf{SML-14} and Lemma \reff{SML-L3}, we conclude that $\eta^{z,w}=\eta^w$. Therefore, from
\rf{SML-7}, we obtain the estimate
\bel{SML-15}
\|\eta^z-\eta^w\|\le C\ve^{-1}\,d(z,w)
~~~\text{for}~~~z,w\in\BXR.
\ee
\par We now define
$$
\ff(z)=\eta^z~~~~\text{for}~~~z\in\BXR.
$$
Then \rf{SML-1}, \rf{SML-2}, \rf{SML-15} tell us that
\bel{SML-17}
\ff(z)\in\Gm_0(z)~~~\text{for all}~~~z\in\BXR,
\ee
\bel{SML-18}
\|\ff(z)-\xi_0\|\le C\ve^{-1}r_0~~~\text{for}~~~z\in\BXR,
\ee
and
\bel{SML-19}
\|\ff(z)-\ff(w)\|\le C\ve^{-1}\,d(z,w)~~~\text{for}~~~z,w\in\BXR.
\ee
\par Our results \rf{SML-17}, \rf{SML-18}, \rf{SML-19} immediately imply the conclusions of the Main Lemma \reff{SML-ML}.
\par This completes the proof of the Main Lemma \reff{SML-ML} in the base case $\#\Ac=m$.\bx
\bigskip
\vskip 3mm
\indent\par {\bf 4.3 Setup for the induction step.}
\addtocontents{toc}{~~~~4.3 Setup for the induction step. \hfill \thepage\par}\medskip
\par Fix a label $\Ac=(e_1,...,e_s)$ with $0\le s\le m-1$.
We assume the
\begin{indhypothesis}\lbl{SIS-IH} {\em Let $x_0^+\in X$, $\xi_0^+\in Y$, $r_0^+>0$, $C_B^+\ge 1$ be given, and let $\Ac^+$ be a label such that $\#\Ac^+>\#\Ac$.
\par Then the Main Lemma \reff{SML-ML} holds, with  $x_0^+$, $\xi_0^+$, $r_0^+$, $C_B^+$, $\Ac^+$, in place of $x_0$, $\xi_0$, $r_0$, $C_B$, $\Ac$, respectively.}
\end{indhypothesis}
\par We assume the
\begin{hmla}\lbl{SIS-HMLA} {\em $x_0\in X$, $\xi_0\in Y$, $r_0>0$, $C_B\ge 1$, $\Gm_{\ell(\Ac)}(x_0)$ has an $(\Ac,\ve^{-1}r_0, C_B)$-basis at $\xi_0$.}
\end{hmla}
\par We introduce a positive constant $A$, and we make the following assumptions.
\begin{lassumption}\lbl{SIS-LAA} {\em $A$ exceeds a large enough constant determined by $m$, $C_B$, $\cng$, $\CNG$.}
\end{lassumption}
\begin{smepsassumption}\lbl{SIS-SEA} {\em $\ve$ is less than a small enough constant determined by $A$, $m$, $C_B$, $\cng$, $\CNG$.}
\end{smepsassumption}
\smallskip
\par We write $C$ to denote a positive constant determined by $m$, $C_B$, $\cng$, $\CNG$; we write $C(\ve,A)$ and $C'(\ve,A)$ to denote positive constants determined by $\ve$, $m$, $A$, $C_B$, $\cng$, $\CNG$. These symbols may denote different constants in different occurrences.
\par Under the above assumptions, we will prove that there exists $\ff:\BXR\to Y$ satisfying
\bel{SIS-A1*}
\|\ff(z)-\ff(w)\|\le C(\ve,A)\,d(z,w)~~~\text{for all} ~~~z,w\in\BXR,
\ee
\bel{SIS-A2*}
\|\ff(z)-\xi_0\|\le C(\ve,A)\,r_0~~~\text{for all}~~~z\in\BXR,
\ee
\bel{SIS-A3*}
\ff(z)\in\Gm_0(z)~~~\text{for all}~~~z\in\BXR.
\ee
These conclusions differ from the conclusions \rf{SML-A1},
\rf{SML-A2}, \rf{SML-A3} of the Main Lemma \reff{SML-ML} only in that here, $C(\ve)$ is replaced by $C(\ve,A)$.
\par Once we have proven the existence of such an $f$ under the above assumptions, we then pick $A$ to be a constant determined by $m$, $C_B$, $\cng$, $\CNG$, taken large enough to satisfy the Large $A$ Assumption \reff{SIS-LAA}.
\par Once we do so, our present Small $\ve$ Assumption \reff{SIS-SEA} will follow from the small $\ve$ assumption made in the Main Lemma \reff{SML-ML}. Moreover, the conclusions \rf{SIS-A1*}, \rf{SIS-A2*}, \rf{SIS-A3*} will then imply conclusions \rf{SML-A1}, \rf{SML-A2}, \rf{SML-A3}. Consequently, we will have proven
the Main Lemma \reff{SML-ML} for $\Ac$. That will complete our downward induction on $\#\Ac$, thereby proving
the Main Lemma \reff{SML-ML} for all labels.
\bigskip
\par To recapitulate:\smallskip
\par We assume the Inductive Hypothesis \reff{SIS-IH} and the Hypotheses of the Main Lemma for the Label $\Ac$ \reff{SIS-HMLA}, and we make the Large $A$ Assumption \reff{SIS-LAA} and the Small $\ve$ Assumption \reff{SIS-SEA}.
\par Under the above assumptions, our task is to prove that there exists  $\ff:\BXR\to Y$ satisfying \rf{SIS-A1*}, \rf{SIS-A2*}, \rf{SIS-A3*}. Once we do that, the Main Lemma \reff{SML-ML} will follow.\medskip
\par We keep the assumptions and notation of this section in force until the end of the proof of the Main Lemma \reff{SML-ML}.
\bigskip\bigskip
\indent\par {\bf 4.4 A family of useful vectors.}
\addtocontents{toc}{~~~~4.4 A family of useful vectors. \hfill \thepage\par}\medskip
\par Recall that $\GLA(x_0)$ has an $(\Ac,\ve^{-1}r_0,C_B)$-basis at $\xi_0$.
\par Let $z\in B(x_0,\pl r_0)$. Then, thanks to our Small $\ve$ Assumption \reff{SIS-SEA}, we have
\bel{FUV-1}
d(z,x_0)\le \pl r_0<\ve_0\cdot(\ve^{-1}r_0),
\ee
where $\ve_0$ arises from $C_B,m$ as in Lemma \reff{LB-L2}.
\medskip
\par We apply that lemma, taking $\Gm=\GLA(x_0)$, $\Gm'=\GLAO(z)$, and $r= \ve^{-1}\,r_0$, and using \rf{FUV-1} and Lemma \reff{BCS-L1P} (B) to verify the key hypothesis (*) in Lemma \reff{LB-L2}. Thus, we obtain a vector $\eta^z\in Y$, with the following properties:
\bel{FUV-UV1}
\GLAO(z)~~~\text{has an}~~~(\Ac,\ve^{-1}r_0,C)\text{-basis
at}~~~\eta^z,
\ee
\bel{FUV-UV2}
\|\eta^z-\xi_0\|\le C\ve^{-1}r_0,
\ee
and
\bel{FUV-UV3}
\ip{e_a,\eta^z-\xi_0}=0~~~\text{for}~~~a=1,...,s.
\ee
We fix such a vector $\eta^z$ for each
$z\in B(x_0,\pl r_0)$.

\indent\par {\bf 4.5 The basic lengthscales.}
\addtocontents{toc}{~~~~4.5 The basic lengthscales. \hfill \thepage\par}\medskip
\begin{definition} {\em Let $x\in B(x_0,5r_0)$, and let $r>0$. We say that $(x,r)$ is OK if both conditions (OK1) and (OK2) below are satisfied.\medskip
\par (OK1)~ $d(x_0,x)+5r\le 5r_0$.\smallskip
\par (OK2)~ Either condition (OK2A) or condition (OK2B) below is satisfied.\medskip
\par~~~~(OK2A)~ $\#B(x,5r)\le 1$ (i.e., $B(x,5r)$ is the singleton $\{x\}$).
\smallskip
\par~~~~(OK2B)~ For some label $\Ac^+$ with $\#\Ac^+>\#\Ac$, the following holds:\medskip
\par\hspace*{19mm} For each $w\in B(x,5r)$ there exists a vector $\zeta^w\in Y$ satisfying conditions (OK2Bi),
\par\hspace*{19mm} (OK2Bii), (OK2Biii) below:\medskip
\par\hspace*{19mm} (OK2Bi)~ $\Gm_{\ell(\Ac)-3}(w)$ has an $(\Ac^+,\ve^{-1}r,A)$-basis at $\zeta^w$.\smallskip
\par\hspace*{19mm} (OK2Bii)~ $\|\zeta^w-\xi_0\|
\le A\ve^{-1}r_0$.\smallskip
\par\hspace*{19mm} (OK2Biii)~ $\ip{e_a,\zeta^w-\xi_0}=0$ for $a=1,...,s$.
}
\end{definition}
\medskip
\par Of course (OK1) guarantees that $B(x,5r)\subset B(x_0,5r_0)$.\medskip
\par Note that $(x,r)$ cannot be OK if $r>r_0$, because then (OK1) cannot hold. On the other hand, if $x\in B(x_0,5r_0)$, then $d(x_0,x)<5r_0$, hence (OK1) holds for
small enough $r$, and (OK2) holds as well (because $B(x,5r)=\{x\}$ for small enough $r$; recall that $(X,d)$ is a finite metric space). Thus, for fixed $x\in B(x_0,5r_0)$, we find that $(x,r)$ is OK if $r$ is small enough, but not if $r$ is too big.
\par For each $x\in B(x_0,5r_0)$ we may therefore
\bel{BLSC}
\text{fix a {\it basic lengthscale}}~~~r(x)>0,
\ee
such that
\bel{BL-1}
(x,r(x))~~\text{is OK, but}~~(x,2r(x))~~\text{is not OK.}
\ee
\par Indeed, we may just take $r(x)$ to be any $r'$ such that $(x,r')$ is OK and
$$
r'>\tfrac12\sup\,\{r:(x,r)~~\text{is OK}\}.
$$
\par We let $\RELX$ (relevant $X$) denote the set of all $x\in B(x_0,5r_0)$ such that
\bel{BL-RELX}
B(x,r(x))\cap B(x_0,r_0)\ne\emp.
\ee
\par Clearly,
\bel{BXR-R}
B(x_0,r_0)\subset \RELX.
\ee
\par From \rf{BL-1} and (OK1), we have
$$
d(x_0,x)+5r(x)\le 5r_0~~~\text{for each}~~~x\in B(x_0,5r_0).
$$
\begin{lemma}\lbl{BL-L1} Let $z_1,z_2\in B(x_0,5r_0)$. If
\bel{BL-A2}
d(z_1,z_2)\le r(z_1)+r(z_2),
\ee
then
$$
\tfrac14r(z_1)\le r(z_2)\le 4r(z_1).
$$
\end{lemma}
\par {\it Proof.} Suppose not. After possibly interchanging $z_1$ and $z_2$, we have
\bel{BL-A1}
r(z_1)<\tfrac14 r(z_2).
\ee
\par Now $(z_2,r(z_2))$ is OK (see \rf{BL-1}). Therefore it satisfies (OK1), i.e.,
$$
d(x_0,z_2)+5r(z_2)\le 5r_0.
$$
Therefore, by \rf{BL-A2},
\be
d(x_0,z_1)+5\cdot (2r(z_1))
&\le& d(x_0,z_2)+d(z_1,z_2) +10r(z_1)
\le d(x_0,z_2)+r(z_1)+r(z_2) +10r(z_1)
\nn\\
&\le& d(x_0,z_2)+\tfrac{11}{4}r(z_2)+r(z_2)
<d(x_0,z_2)+5r(z_2)\le 5r_0,\nn
\ee
i.e., $(z_1,2r(z_1))$ satisfies (OK1).
\par Moreover,
\bel{BL-A3}
B(z_1,10r(z_1))\subset B(z_2,5r(z_2)).
\ee
Indeed, if $w\in B(z_1,10r(z_1))$, then \rf{BL-A1} and \rf{BL-A2} give
$$
d(w,z_2)\le  d(w,z_1)+d(z_1,z_2)
\le 10r(z_1)+r(z_1)+r(z_2)
\le \tfrac{11}{4}r(z_2)+r(z_2)<5r(z_2),
$$
proving \rf{BL-A3}.
\par Because $(z_2,r(z_2))$ is OK, it satisfies (OK2A) or (OK2B). If $(z_2,r(z_2))$ satisfies (OK2A), then so does
$(z_1,2r(z_1))$, thanks to  \rf{BL-A3}. In that case,
$(z_1,2r(z_1))$ satisfies (OK1) and (OK2A), hence
$(z_1,2r(z_1))$ is OK, contradicting  \rf{BL-1}.
\smallskip
\par On the other hand, suppose $(z_2,r(z_2))$ satisfies (OK2B). Fix $\Ac^+$ with $\#\Ac^+>\#\Ac$ such that for every $w\in B(z_2,5r(z_2))$ there exists $\zeta^w$ satisfying
\medskip
\par \textbullet~ $\Gm_{\ell(\Ac)-3}(w)$ has an $(\Ac^+,\ve^{-1}r(z_2),A)$-basis at $\zeta^w$.
\medskip
\par \textbullet~ $\|\zeta^w-\xi_0\|\le A\ve^{-1}r_0$.
\medskip
\par \textbullet~ $\ip{e_a,\zeta^w-\xi_0}=0$~~ for $a=1,...,s$.
\medskip
\par Thanks to \rf{BL-A3} there exists such a\, $\zeta^w$ for every $w\in B(z_1,5\cdot(2r(z_1)))$. 
\par Note that, by \rf{BL-A1} and Remark \reff{LBL-R} (ii), the $(\Ac^+,\ve^{-1}r(z_2),A)$-basis in the first bullet point above is also an $(\Ac^+,\ve^{-1}\cdot(2r(z_1)),A)$-basis.
\par It follows that $(z_1,2r(z_1))$ satisfies (OK2B). We have seen that $(z_1,2r(z_1))$ satisfies (OK1), so again
$(z_1,2r(z_1))$ is OK, contradicting \rf{BL-1}.
\par Thus, in all cases, our assumption that Lemma \reff{BL-L1} fails leads to a contradiction.\bx
\bigskip
\indent\par {\bf 4.6 Consistency of the useful vectors.}
\addtocontents{toc}{~~~~4.6 Consistency of the useful vectors. \hfill \thepage\par}\medskip
\par Recall the useful vectors $\eta^z$ $(z\in B(x_0,\pl r_0))$, see \rf{FUV-UV1}, \rf{FUV-UV2}, \rf{FUV-UV3}, and the set $\RELX$, see \rf{BL-RELX}. In this section we establish the following result.
\begin{lemma}\lbl{CUV-L1} Let $z_1,z_2\in \RELX$. Then
$$
\|\eta^{z_1}-\eta^{z_2}\|\le C\ve^{-1}[r(z_1)+r(z_2)+d(z_1,z_2)].
$$
\end{lemma}
\par {\it Proof.} If
$$
r(z_1)+r(z_2)+d(z_1,z_2)\ge r_0/\pl,
$$
then the lemma follows from \rf{FUV-UV2} applied to $z=z_1$ and to $z=z_2$.
\par Suppose
\bel{MN-910}
r(z_1)+r(z_2)+d(z_1,z_2)< r_0/\pl.
\ee
\par Because $z_1\in\RELX$, we have $d(z_1,x_0)\le r_0+r(z_1)$, hence
$$
d(z_1,x_0)+5\cdot(2r(z_1))\le r_0+11r(z_1)<5r_0.
$$
Thus $(z_1,2r(z_1))$ satisfies (OK1), and, in particular,  $B(z_1,10r(z_1))\subset B(x_0,5r_0)$.
\par Recall from \rf{FUV-UV1} that $\GLAO(z_2)$ has an $(\Ac,\ve^{-1}r_0,C)$-basis at $\eta^{z_2}$. By \rf{MN-910} and Remark \reff{LBL-R} (ii), it follows that
\bel{CUV-1}
\GLAO(z_2)~~~\text{has an}~~~(\Ac,\ve^{-1}[r(z_1)+r(z_2)+d(z_1,z_2)],C)
\text{-basis}~~\text{at}~~
\eta^{z_2}.
\ee
\par Our Small $\ve$ Assumption \reff{SIS-SEA} shows that
$$
d(z_1,z_2)\le \ve_0\cdot
\ve^{-1}[r(z_1)+r(z_2)+d(z_1,z_2)],
$$
for the $\ve_0$ arising from Lemma \reff{LB-L2}, where we use the constant $C$ in \rf{CUV-1} as the constant $C_B$ in Lemma \reff{LB-L2}.
\par Therefore, by Lemma \reff{LB-L2} and Lemma \reff{BCS-L1P} (B), with
$$
\Gm=\GLAO(z_2),~~\Gm'=\Gm_{\ell(\Ac)-2}(z_1),~~
r=\ve^{-1}[r(z_1)+r(z_2)+d(z_1,z_2)],
$$
we obtain a vector $\zeta\in\Gm_{\ell(\Ac)-2}(z_1)$
such that
\bel{444}
\|\zeta-\eta^{z_2}\|\le C\ve^{-1}[r(z_1)+r(z_2)+d(z_1,z_2)]
\ee
and
$$
\ip{e_a,\zeta-\eta^{z_2}}=0~~~\text{for}~~~a=1,...,s,
$$
hence
\bel{CUV-5}
\ip{e_a,\zeta-\eta^{z_1}}=0~~~\text{for}~~~a=1,...,s.~~~~
\text{(\,See \rf{FUV-UV3}.\,)}
\ee
\par We will prove that
$$
\|\zeta-\eta^{z_1}\|\le \ve^{-1}r(z_1);
$$
\rf{444} will then imply the conclusion of Lemma \reff{CUV-L1}.
\par Suppose instead that
\bel{CUV-6}
\|\zeta-\eta^{z_1}\|> \ve^{-1}r(z_1).
\ee
\par We will derive a contradiction.
By \rf{FUV-UV1}, Remark \reff{LBL-R} (iv), and because $r(z_1)<r_0/\pl$ (see \rf{MN-910}), we know that
\bel{CUV-7}
\Gm_{\ell(\Ac)-2}(z_1)~~~\text{has an}~~~ (\Ac,\ve^{-1}r(z_1),C)\text{-basis at}~~\eta^{z_1}.
\ee
Our results \rf{CUV-5}, \rf{CUV-7} and our assumption \rf{CUV-6} are the hypotheses of Lemma \reff{LB-L1} (``Adding a vector''). Applying that lemma, we obtain a vector $\hz\in\Gm_{\ell(\Ac)-2}(z_1)$, with the following properties:
\bel{CUV-9}
\|\hz-\eta^{z_1}\|=\tfrac12 \ve^{-1}r(z_1),
\ee
$$
\ip{e_a,\hz-\eta^{z_1}}=0~~~\text{for}~~~a=1,...,s;
$$
also
\bel{CUV-11}
\Gm_{\ell(\Ac)-2}(z_1)~~~\text{has an}~~~ (\Ac^+,\ve^{-1}r(z_1),C)\text{-basis at}~~\hz,
\ee
for a label of the form $\Ac^+=(e_1,...,e_s,e_{s+1})$; and 
\bel{CUV-10}
\ip{e_a,\hz-\xi_0}=0~~~\text{for}~~~a=1,...,s.
\ee
See \rf{FUV-UV3}.
\par In particular, $\#\Ac^+=\#\Ac+1$.
\smallskip
\par From \rf{CUV-11} and Remark \reff{LBL-R} (iii) we have, with a larger constant $C$,
\bel{CUV-13}
\Gm_{\ell(\Ac)-2}(z_1)~~~\text{has an}~~~ (\Ac^+,\ve^{-1}\cdot(2r(z_1)),C)\text{-basis at}~~\hz.
\ee
\par Now let $w\in B(z_1,5\cdot(2r(z_1)))$. Let $\ve_0$ arise from Lemma \reff{LB-L2} where we use $C$ from \rf{CUV-13} as the constant $C_B$ in Lemma \reff{LB-L2}. We have
$$
d(z_1,w)< 10r(z_1)<\ve_0\cdot(\ve^{-1}\cdot(2r(z_1))),
$$
thanks to our Small $\ve$ Assumption \reff{SIS-SEA}. Therefore, Lemma \reff{BCS-L1P} (B) allows us to verify the key hypothe\-sis (*) in  Lemma \reff{LB-L2}, with $\Gm=\Gm_{\ell(\Ac)-2}(z_1)$, $\Gm'=\Gm_{\ell(\Ac)-3}(w)$, $r=\ve^{-1}\cdot(2r(z_1))$.
\par Applying Lemma \reff{LB-L2}, we obtain a vector
$\zeta^w\in\Gm_{\ell(\Ac)-3}(w)$ with the following properties:
\bel{CUV-15}
\|\zeta^w-\hz\|\le C\ve^{-1}\cdot(2r(z_1)),
\ee
$$
\ip{e_a,\zeta^w-\hz}=0~~~\text{for}~~~a=1,...,s+1;
$$
hence by \rf{CUV-10},
\bel{CUV-17}
\ip{e_a,\zeta^w-\xi_0}=0~~~\text{for}~~~a=1,...,s.
\ee
Also,
\bel{CUV-18}
\Gm_{\ell(\Ac)-3}(w)~~~\text{has an}~~~ (\Ac^+,\ve^{-1}\cdot(2r(z_1)),C)\text{-basis at}~~\zeta^w.
\ee
\par We have
$$
\|\zeta^w-\xi_0\|\le \|\zeta^w-\hz\|+\|\hz-\eta^{z_1}\|+
\|\eta^{z_1}-\xi_0\|
\le  C\ve^{-1}r(z_1)+\tfrac12\ve^{-1}r(z_1)+C\ve^{-1}r_0
$$
by \rf{CUV-15}, \rf{CUV-9} and \rf{FUV-UV2}.
\par Recalling that $r(z_1)<r_0/\pl$, we conclude that
\bel{CUV-19}
\|\zeta^w-\xi_0\|\le C\ve^{-1}\cdot r_0.
\ee
\par Thus, for every $w\in B(z_1,5\cdot(2r(z_1)))$, our vector $\zeta^w$ satisfies \rf{CUV-17}, \rf{CUV-18}, \rf{CUV-19}. Comparing \rf{CUV-18}, \rf{CUV-19}, \rf{CUV-17} with (OK2Bi), (OK2Bii), (OK2Biii), and recalling our Large $A$ Assumption \reff{SIS-LAA},
we conclude that (OK2B) holds for $(z_1,2r(z_1))$. We have already seen that (OK1) holds for $(z_1,2r(z_1))$. Thus $(z_1,2r(z_1))$ is OK, contradicting the defining property \rf{BL-1} of $r(z_1)$.
\medskip
\par This contradiction proves that \rf{CUV-6} cannot hold, completing the proof of Lemma \reff{CUV-L1}.\bx
\bigskip\bigskip
\indent\par {\bf 4.7  Additional useful vectors.}
\addtocontents{toc}{~~~~4.7  Additional useful vectors. \hfill \thepage\par}\medskip
\begin{lemma}\lbl{AUV-L1} Let $x\in B(x_0,5r_0)$, and suppose that $\#B(x,5r(x))\ge 2$.
\par Then there exist a vector $\zeta^x\in Y$ and a label $\Ac^+$ with the following properties:\medskip
\bel{AUV-1}
\#\Ac^+>\#\Ac\,,
\ee
\bel{AUV-2}
\Gm_{\ell(\Ac)-3}(x)~~~\text{has an}~~~ (\Ac^+,\ve^{-1}r(x),A)\text{-basis at}~~\zeta^x,
\ee
\bel{AUV-3}
\|\zeta^x-\eta^x\|\le \ve^{-1}r(x),
\ee
\bel{AUV-4}
\ip{e_a,\zeta^x-\eta^x}=0~~~\text{for}~~~a=1,...,s.
\ee
\end{lemma}
\par {\it Proof.} Recall that $(x,r(x))$ is OK. We are assuming that (OK2A) fails for $(x,r(x))$, hence (OK2B)
holds. Fix $\Ac^+$ as in (OK2B), and let $\zeta^x$ be as in (OK2B) with $w=x$. Then \rf{AUV-1}, \rf{AUV-2}, \rf{AUV-4}
hold, thanks to (OK2B); however, \rf{AUV-3} may fail in case $r(x)$ is much smaller than $r_0$. If \rf{AUV-3} holds, we are done. 
\par Suppose instead that \rf{AUV-3} fails, i.e.,
\bel{AUV-5}
\|\zeta^x-\eta^x\|> \ve^{-1}r(x).
\ee
We recall from \rf{FUV-UV1} that $\Gm_{\ell(\Ac)-1}(x)$ has an  $(\Ac,\ve^{-1}r_0,C)$-basis at $\eta^x$. We have also $r(x)\le r_0$ because $(x,r(x))$ is OK; and $$\Gm_{\ell(\Ac)-1}(x)\subset\Gm_{\ell(\Ac)-3}(x).$$ Therefore, by Remark \reff{LBL-R} (iv),
\bel{AUV-6}
\Gm_{\ell(\Ac)-3}(x)~~~\text{has an}~~~ (\Ac,\ve^{-1}r(x),C)\text{-basis at}~~\eta^x.
\ee
\par From \rf{AUV-4}, \rf{AUV-5}, \rf{AUV-6} and Lemma \reff{LB-L1} (``Adding a vector''), we obtain a vector $\hz\in Y$ and a label $\hAc$ with the following properties:
\bel{AUV-7}
\#\hAc>\#\Ac\,,
\ee
\bel{AUV-8}
\|\hz-\eta^x\|= \tfrac12 \ve^{-1}r(x),
\ee
\bel{AUV-9}
\ip{e_a,\hz-\eta^x}=0~~~\text{for}~~~a=1,...,s,
\ee
\bel{AUV-10}
\Gm_{\ell(\Ac)-3}(x)~~~\text{has an}~~~ (\hAc,\ve^{-1}r(x),C')\text{-basis at}~~\hz.
\ee
\par Comparing \rf{AUV-7},...,\rf{AUV-10} with \rf{AUV-1},...,\rf{AUV-4}, and recalling our Large $A$ Assumption \reff{SIS-LAA}, we see that $\hz$ and $\hAc$ have all the properties asserted for $\zeta^x$ and $\Ac^+$ in the statement of Lemma \reff{AUV-L1}.
\medskip
\par Thus, Lemma \reff{AUV-L1} holds in all cases.\bx
\bigskip\bigskip
\indent\par {\bf 4.8  Local selections.}
\addtocontents{toc}{~~~~4.8  Local selections. \hfill \thepage\par}\medskip
\begin{lemma}\lbl{LS-L1} Given $x\in\RELX$, there exists $\ff:B(x,r(x))\to Y$ with the following properties:
\medskip
\par (I)~ $\|\ff(z)-\ff(w)\|\le C(\ve,A)\,d(z,w)$~~ for~ $z,w\in B(x,r(x))$.
\medskip
\par (II)~ $\ff(z)\in\Gm_0(z)$~~ for~ $z\in B(x,r(x))$.
\medskip
\par (III)~ $\|\ff(z)-\eta^x\|\le C(\ve,A)\,r(x)$~~ for~ $z\in B(x,r(x))$.
\medskip
\par (IV)~ $\|\ff(z)-\xi_0\|\le C(\ve,A)\,r_0$~~ for~ $z\in B(x,r(x))$.
\medskip
\end{lemma}
\par {\it Proof.} We proceed by cases.\medskip
\par {\it Case 1.} Suppose $\#B(x,5r(x))>1$.
\par Then Lemma \reff{AUV-L1} applies. Let $\Ac^+$, $\zeta^x$ be as in that lemma. Thus,
\bel{LS-1}
\#\Ac^+>\#\Ac\,,
\ee
\bel{LS-2}
\|\zeta^x-\eta^x\|\le \ve^{-1}r(x)
\ee
and
$$
\Gm_{\ell(\Ac)-3}(x)~~~\text{has an}~~~ (\Ac^+,\ve^{-1}r(x),A)\text{-basis at}~~\zeta^x\,;
$$
hence, by Remark \reff{LBL-R} (iv),
\bel{LS-4}
\Gm_{\ell(\Ac^+)}(x)~~~\text{has an}~~~ (\Ac^+,\ve^{-1}r(x),A)\text{-basis at}~~\zeta^x,
\ee
because $\ell(\Ac)-3\ge \ell(\Ac^+)$ whenever $\#\Ac^+>\#\Ac$.
\par We recall from our Small $\ve$ Assumption \reff{SIS-SEA} that
\bel{LS-5}
\ve~~\text{is less than a small enough constant determined by}~A,\cng,\CNG,m.
\ee
\par Thanks to \rf{LS-4}, \rf{LS-5}, the Hypotheses of the Main Lemma \reff{SIS-HMLA} are satisfied, with $\Ac^+$, $x$, $\zeta^x$, $r(x)$, $A$, in place of
$\Ac$, $x_0$, $\xi_0$, $r_0$, $C_B$, respectively. Moreover, thanks to \rf{LS-1} and the Inductive Hypothesis \reff{SIS-IH}, we are assuming the validity of the Main Lemma \reff{SML-ML} for $\Ac^+,...,A$.
\par Therefore, we obtain a function $f:B(x,r(x))\to Y$ satisfying (I), (II) and the inequality
$$
\|\ff(z)-\zeta^x\|\le C(\ve,A)\,r(x),~~~~ z\in B(x,r(x)).
$$
This inequality together with \rf{LS-2} implies
(III).
\par Moreover, (IV) follows from (III) because, for $z\in B(x,r(x))\subset B(x_0,5r_0)$, we have
$$
\|\ff(z)-\xi_0\|\le \|\ff(z)-\eta^x\|+\|\eta^x-\xi_0\|
\le C(\ve,A)r(x)+C\ve^{-1}r_0\le C'(\ve,A)r_0;
$$
here we use \rf{FUV-UV2} and the fact that $(x,r(x))$ satisfies (OK1).
\par This completes the proof of Lemma \reff{LS-L1} in Case 1.\bigskip
\par {\it Case 2.} Suppose $\#B(x,5r(x))\le 1$.
\par Then, $B(x,5r(x))=\{x\}$ and $\eta^x\in\GLAO(x)\subset\Gamma_0(x)$. Hence the function
$f(x)=\eta^x$ satisfies (I),(II),(III), and also (IV) thanks to \rf{FUV-UV2}.
\par Thus, Lemma \reff{LS-L1} holds in all cases.\bx
\bigskip\bigskip\bigskip
\par {\bf 4.9 Proof of the Main Lemma: the final step.}
\addtocontents{toc}{~~~~4.9 Proof of the Main Lemma: the final step. \hfill \thepage\par\VSU}\medskip
\indent\par Let $\Bc_0$ be the metric space
$$
\Bc_0=\left(B(x_0,r_0),d|_{B(x_0,r_0)\times B(x_0,r_0)}\right),
$$
i.e., the ball $B(x_0,r_0)$ supplied with the metric $d$.
\par For the rest of Section 4.9, we work in the metric space $\Bc_0$. Given $x\in B(x_0,r_0)$ and $r>0$, we write $\BO(x,r)$ to denote the ball in $\Bc_0$ with center $x$ and radius $r$; thus $\BO(x,r)= B(x,r)\cap B(x_0,r_0)$.
\par Note that, since $(X,d)$ satisfies Nagata $(\CNG,\cng)$, the metric space $\Bc_0$ satisfies Nagata $(\CNG,\cng)$ as well.
See Definition \reff{NG-C}.
\par Let $r:X\to\R_+$ be {\it the basic lengthscale} constructed in Section 4.5 (see \rf{BLSC}), and let
\bel{A-CL}
\CLS=4 ~~~\text{and}~~~a=(4\,\CLS)^{-1}.
\ee
Note that, by Lemma \reff{BL-L1}, {\sc Consistency of the Lengthscale} (see \rf{C-LSC}) holds for the lengthscale $r(x)$ on $B(x_0,r_0)$ with the constant $\CLS$ given by \rf{A-CL}.
\par We apply the Whitney Partition Lemma \reff{WPL} to the metric space $\Bc_0$, the lengthscale
$$\{r(x):x\in B(x_0,r_0)\}$$
and the constants $\CLS$, $a$ determined by \rf{A-CL}, and obtain a partition of unity ${\{\theta_\nu:
B(x_0,r_0)\to\R_+\}}$ and points 
\bel{XV-BR}
x_\nu\in B(x_0,r_0)
\ee
with the following properties.
\medskip
\par \textbullet~ Each $\theta_\nu\ge 0$ and for each $\nu$, $\theta_\nu=0$ outside $\BO(x_\nu,a r_\nu)$; here
$a$ is determined by \rf{A-CL}, and $r_\nu=r(x_\nu)$.
\smallskip
\par \textbullet~ Any given $x$ satisfies $\theta_\nu(x)\ne 0$ for at most $\DS$ distinct $\nu$, where $\DS$ depends only on $\cng$, $\CNG$.
\smallskip
\par \textbullet~ $\smed\limits_\nu\,\theta_\nu(x)=1$ for all $x\in B(x_0,r_0)$.
\smallskip
\par \textbullet~ Each $\theta_\nu$ satisfies
$$
|\theta_\nu(x)-\theta_\nu(y)|\le \frac{C}{r_\nu}\,d(x,y)
$$
for all $x,y\in B(x_0,r_0)$; here again $r_\nu=r(x_\nu)$.
\smallskip
\par From Lemma \reff{BL-L1}, we know that
\smallskip
\par \textbullet~ For each $\mu,\nu$, if $d(x_\mu,x_\nu)\le r_\mu+r_\nu$, then $\frac14 r_\nu\le r_\mu\le 4r_\nu$.
\smallskip
\par Moreover, by \rf{BXR-R} and \rf{XV-BR},
\bel{RELX-N}
x_\nu\in\RELX~~~\text{for each}~~\nu,
\ee
so that, by Lemma \reff{LS-L1}, there exists a function $\hf_\nu:B(x_\nu,r_\nu)\to Y$ satisfying the following conditions
\smallskip
\par \textbullet~ $\|\hf_\nu(z)-\hf_\nu(w)\|\le C(\ve,A)\,d(z,w)$~~ for~ $z,w\in B(x_\nu,r_\nu)$.
\smallskip
\par \textbullet~ $\hf_\nu(z)\in\Gm_0(z)$~~ for~ $z\in B(x_\nu,r_\nu)$.
\smallskip
\par \textbullet~ $\|\hf_\nu(z)-\eta_\nu\|\le C(\ve,A)\,r_\nu$~~ for~ $z\in B(x_\nu,r_\nu)$, where $\eta_\nu\equiv \eta^{x_\nu}$.
\smallskip
\par \textbullet~ $\|\hf_\nu(z)-\xi_0\|\le C(\ve,A)\,r_0$~~ for~ $z\in B(x_\nu,r_\nu)$.
\medskip
\par Let $f_\nu=\hf_\nu|_{\BO(x_\nu,r_\nu)}$. We extend $f_\nu$ from $\BO(x_\nu,r_\nu)=B(x_\nu,r_\nu)\cap B(x_0,r_0)$ to all of $B(x_0,r_0)$ by setting $f_\nu=0$ outside $\BO(x_\nu,r_\nu)$.
\smallskip
\par Since each $x_\nu\in\RELX$ (see \rf{RELX-N}), from Lemma \reff{CUV-L1}, we have
\smallskip
\par \textbullet~ $\|\eta_\nu-\eta_\mu\|\le C(\ve,A)\cdot
[r_\nu+r_\mu+d(x_\nu,x_\mu)]$ for each $\mu,\nu$.
\medskip
\par The above conditions on the $\theta_\nu$, $\eta_\nu$, $\hf_\nu$, $f_\nu$, $r_\nu$ and $a$ (cf. \rf{A-SMALL} with \rf{A-CL}) allow us to apply the Patching Lemma
\reff{PTHM} on $\Bc_0$. We conclude that
$$
f(x)=\smed_\nu\,\theta_\nu(x)\,f_\nu(x)~~~
(\text{all}~~x\in B(x_0,r_0))
$$
satisfies
$$
\|f(x)-f(y)\|\le C(\ve,A)\,d(x,y)~~~\text{for}~~
x,y\in B(x_0,r_0).
$$
\par Moreover, for fixed $x\in B(x_0,r_0)$, we know that $f(x)$ is a convex combination of finitely many values $f_\nu(x)$ with $\BO(x_\nu,ar_\nu)\ni x$; for those $\nu$ we have $f_\nu(x)\in \Gamma_0(x)$ and $\|f_\nu(x)-\xi_0\|\le C(\ve,A)\,r_0$. Therefore, $f(x)\in \Gamma_0(x)$ and
$\|f(x)-\xi_0\|\le C(\ve,A)\,r_0$ for all $x\in B(x_0,r_0)$.
\medskip
\par Thus, $f$ satisfies \rf{SIS-A1*}, \rf{SIS-A2*} and \rf{SIS-A3*}, completing the proof of the Main Lemma \reff{SML-ML}.\bx
\medskip
\medskip
\par {\it Proof of the Finiteness Theorem \reff{PFT-FT} for bounded Nagata dimension.} Let $x_0\in X$, $r_0=\diam X+1$, $C_B=1$, and $\Ac=(~)$. Let $\ve=\tfrac12\,\emin$ where $\emin$ is as in the Main Lemma \reff{SML-ML} for $m$, $C_B=1$, $\cng$ and $\CNG$. Thus, $\ve$ depends only on $m$, $\cng$ and $\CNG$.
\par By Lemma \reff{BCS-L1P} (A), $\Gm_{\ell(\Ac)}(x_0)\ne\emp$ so that there exists $\xi_0\in \Gm_{\ell(\Ac)}(x_0)$. Since $\#\Ac=0$, the set
$\Gm_{\ell(\Ac)}(x_0)$ has an $(\Ac,\ve^{-1}r_0,C_B)$-basis at $\xi_0$. See Remark \reff{LBL-R}, (i).
\par Hence, by the Main Lemma \reff{SML-ML}, there exists a mapping $f:\BXR\to Y$ such that
$$
\|f(z)-f(w)\|\le C\,d(z,w)~~~\text{for all}~~~z,w\in \BXR,
$$
and
$$
f(z)\in\Gm_0(z)~~~\text{for all}~~~z\in \BXR.
$$
Here $C$ is a constant determined by $\ve$, $m$, $C_B$, $\cng$, $\CNG$. Thus, $C$ depends only on $m$, $\cng$, $\CNG$. 
\par Clearly, $\BXR=X$. Furthermore, $\Gm_0(z)\subset F(z)$ for every $z\in X$ (see \rf{GL-F2}), so that $f(z)\in F(z)$, $z\in X$. Thus, $f$ is a Lipschitz selection of $F$ on $X$ with Lipschitz seminorm at most a certain constant depending only on $m$, $\cng$, $\CNG$.
\par The proof of Theorem \reff{PFT-FT} is complete.\bx
\medskip
\smallskip
\par Recall that Theorem \reff{PFT-FT} immediately implies Theorem \reff{NG-FP}.
\bigskip\smallskip

\indent\par {\bf 4.10 The Finiteness Principle on metric trees.}
\addtocontents{toc}{~~~~4.10 The Finiteness Principle on metric trees. \hfill \thepage\par}\medskip

\smallskip
\par Let us consider an important example of a metric space with finite Nagata dimension.
\par Let $T=(X,E)$ be a {\it finite tree}. Here $X$ denotes the set of nodes and $E$ denotes the set of edges of $T$. We write $x\je y$ to indicate that nodes $x,y\in X$, $x\ne y$, are joined by an edge; we denote that edge by $[xy]$.
\par Suppose we assign a positive number $\Delta(e)$ to each edge $e\in E$. Then we obtain a notion of distance $d(x,y)$ for any $x,y\in X$, as follows.
\par We set
\bel{D-XX}
d(x,x)=0~~~\text{for every}~~~x\in X.
\ee
\par Because $T$ is a tree, any two distinct nodes $x,y\in X$ are joined by one and only ``path''
$$
x=x_0\je x_1\je...\je x_L=y~~~\text{with all the}~~x_i ~~\text{distinct}.
$$
\par We define
\bel{D-TR}
d(x,y)=\smed_{i=1}^L\, \Delta([x_{i-1}x_i]).
\ee
We call the resulting metric space $(X,d)$ a {\it metric tree}.
\medskip
\par For the reader's convenience we prove the following slight variant of a result from \cite{LS}.
\begin{lemma}\lbl{MTR-ND} Every metric tree satisfies Nagata $(1,c)$ with $c=1/16$. (See Definition \reff{NG-C}).
\end{lemma}
\par {\it Proof.} Given a metric tree $(X,d)$, we fix an origin $0\in X$ and make the following definition:
\bigskip
\par Every point $x\in X$ is joined to the origin by one and only one ``path''
$$
0=x_0\je x_1\je...\je x_L=x,
~~~\text{with all the}~~x_i~~\text{\it distinct.}
$$
\par We call $x_0,x_1,..., x_L$ the {\it ancestors} of $x$. We define the {\it distinguished ancestor} of $x$, denoted  $\DA(x)$, to be $x_i$ for the smallest $i\in \{0,...,L\}$ for which
\bel{ANC}
d(0,x_i)>\lfloor d(0,x)\rfloor-1,
\ee
where $\lfloor \cdot\rfloor$ denotes the greatest integer function. (Note that there is at least one $x_i$ satisfying \rf{ANC}, namely $x_L=x$. Thus, every $x\in X$ has a distinguished ancestor.)
\medskip
\par We note two simple properties of $\DA(x)$, namely,
\medskip
\par (1) $d(x,\DA(x))\le 2$;
\medskip
\par (2) $\DA(x)$ is an ancestor of any ancestor $y$ of $x$ that satisfies  $d(0,y)>\lfloor d(0,x)\rfloor-1$.
\medskip
\par We now exhibit a Nagata covering of $X$ for the lengthscale $s=4$.\smallskip
\par For $q=0,1$ and $z\in X$, let
$$
X_q(z)=\{x\in X:
z=\DA(x)~~\text{and}~~\lfloor d(0,x)\rfloor\equiv q\mod 2\}.
$$
Clearly, the $X_q(z)$ cover $X$. Moreover, (1) tells us that each $X_q(z)$ has diameter at most $4$.
\par We assert the following\smallskip
\par {\sc Claim:} If $z\ne z'$ and $q=q'$, then the distance from $X_q(z)$ to $X_{q'}(z')$ is at least $1/2$.
\medskip
\par The {\sc Claim} immediately implies that any given ball $B\subset X$ of radius $1/4$ meets at most one of the $X_0(z)$ and at most one of the $X_1(z)$, hence at most two of the $X_q(z)$.\smallskip
\par Let us establish the {\sc Claim}; if it were false, then we could find
$$
z\ne z',~q\in\{0,1\},~x\in X_q(z),~x'\in X_{q}(z')~~~
\text{with}~~~d(x,x')\le 1/2.
$$
\par We will derive a contradiction from these conditions as follows.
\par Because $d(x,x')\le 1/2$, we have
$$
|\,\lfloor d(0,x)\rfloor-\lfloor d(0,x')\rfloor\,|\le 1.
$$
On the other hand, $\lfloor d(0,x)\rfloor\equiv \lfloor d(0,x')\rfloor \mod 2$. Hence, $\lfloor d(0,x)\rfloor= \lfloor d(0,x')\rfloor$.
\par Next, let $\tz$ be the closest common ancestor of $x,x'$. Because $d(x,x')\le 1/2$, we have $d(x,\tz)\le 1/2$ and $d(x',\tz)\le 1/2$, and therefore the ancestor $\tz$ of $x$ satisfies
$$
d(0,\tz)>\lfloor d(0,x)\rfloor-1.
$$
Hence, (2) implies that $z$ is an ancestor of $\tz$. Similarly, $z'$ is an ancestor of $\tz$.
\par It follows that either $z$ is an ancestor of $z'$, or $z'$ is an ancestor of $z$. Without loss of generality, we may suppose that $z$ is an ancestor of $z'$. Consequently, $z$ is an ancestor of $x'$; moreover,
$$
d(0,z)>\lfloor d(0,x)\rfloor-1=\lfloor d(0,x')\rfloor-1.
$$
\par Thanks to (2), we now know that $z'$ is an ancestor of $z$. Thus, each of the points $z,z'$ is an ancestor of the other, and therefore $z=z'$, contradicting an assumption that the {\sc Claim} is false.
\smallskip
\par We have produced a covering of an arbitrary metric tree by subsets $X_i$ of diameter at most $4$, such that no ball of radius $1/4$ intersects more than two of the $X_i$.
\par Applying the above result to the metric tree $(X,\ltfrac{4}{s}\,d)$ for given $s>0$, we produce a covering of $X$ by $X_i$ such that, with respect to $d$, each $X_i$ has diameter at most $s$, and no ball of radius $s/16$ meets more than two of the $X_i$. Thus, we have verified the Nagata condition for metric trees.\bx
\smallskip
\par Let us apply Theorem \ref{NG-FP} to {\it metric trees}. Thus, we obtain the following
\begin{corollary}\lbl{A-12} Let $m\in\N$, let $(X,d)$ be a metric tree and let $\lambda$ be a positive constant.
\par Let ${F:X\to\CNMY}$ be a set-valued mapping such that, for every subset $X'\subset X$ with $\#X'\le \ks$, the restriction $F|_{X'}$ has a Lipschitz selection $f_{X'}:X'\to\BS$ with $\|f_{X'}\|_{\Lip(X',\BS)}\le \lambda$.
\par Then $F$ has a Lipschitz selection $f:X\to\BS$ with $\|f\|_{\Lip(X,\BS)}\le \gamma_0\,\lambda$.
\par Here $\ks=\ks(m)$ is the constant from Theorem \reff{NG-FP}, and $\gamma_0=\gamma_0(m)$ is a constant depending only on $m$.
\end{corollary}

\SECT{5. Metric trees and Lipschitz selections with respect to the Hausdorff distance}{5}
\addtocontents{toc}{5. Metric trees and Lipschitz selections with respect to the Hausdorff distance
\hfill\thepage\par\VST}

\indent\par We recall that $(Y,\|\cdot\|)$ denotes a Banach space, and $\KM$ denotes the family of all nonempty compact convex subsets $K\subset Y$ of dimension at most $m$. Recall also that $\dhf(A,B)$ denotes the Hausdorff distance between $A,B\in\KM$.
\par In this section we work with finite trees $T=(X,E)$, where $X$ denotes the set of nodes and $E$ denotes the set of edges of $T$. As in Section 4.10, we write $u\je v$ to indicate that $u,v\in X$ are distinct nodes joined by an edge in $T$. \smallskip
\par We supply $X$ with a metric $d$ defined by formulae \rf{D-XX} and \rf{D-TR}, and we refer to the metric space $(X,d)$ as a {\it metric tree (with respect to the tree $T=(X,E)$)}.
\begin{remark} {\em Sometimes we will be looking simultaneously at two different pseudometrics, say $\rho$ and $\trh$, on a pseudometric space, say on $\Mc$. In this case we will speak of a $\rho$-Lipschitz selection and $\rho$-Lipschitz seminorm, or a $\trh$-Lipschitz selection and $\trh$-Lipschitz seminorm to make clear which pseudometric we are using. Furthermore, sometimes given a mapping $f:\Mc\to Y$ we will write $\|f\|_{\Lip((\Mc,\rho),Y)}$ to denote the Lipschitz seminorm of $f$ with respect to the pseudometric $\rho$.}
\end{remark}

\bigskip

\indent\par {\bf 5.1 The ``core'' of a set-valued mapping and the Finiteness Principle.}
\addtocontents{toc}{~~~~5.1 The ``core'' of a set-valued mapping and the Finiteness Principle. \hfill \thepage\par}\medskip
\par Until the end of Section 6 we write $\ks$ and $\gz$ to denote the constants from  Corollary \reff{A-12}. Recall that these constants depend only on $m$.
\smallskip
\par In this and the next subsection we prove the following result.
\begin{theorem}\lbl{HDS-M} Let $(\Mc,\rho)$ be a metric space, and let $F:\Mc\to\KM$ for a Banach space $Y$. Let $\lambda$ be a positive real number.
\par Suppose that for every subset $\Mc'\subset\MS$ consisting of at most $\ks$ points, the restriction $F|_{\Mc'}$ has a Lipschitz selection $f_{\Mc'}$ with Lipschitz seminorm $\|f_{\Mc'}\|_{\Lip(\Mc',\BS)}\le \lambda$.
\par Then there exists a mapping $G:\MS\to\KM$ satisfying the following conditions:\smallskip
\par (i). $G(x)\subset F(x)$ for every $x\in\MS$;
\medskip
\par (ii). For every $x,y\in \MS$ the following inequality
$$
\dhf(G(x),G(y))\le \gz\,\lambda\,\dm(x,y)
$$
holds.
\end{theorem}
\smallskip
\par Let $(\MS,\dm)$ be a metric space and let $F:\MS\to\KM$ be a set-valued mapping. We suppose that the following assumption is satisfied.
\begin{assumption}\lbl{A-F} For every subset $\MS'\subset\MS$ consisting of at most $\ks$ points, the restriction $F|_{\MS'}$ of $F$ to $\MS'$ has a $\dm$-Lipschitz selection $f_{\MS'}:\MS'\to\BS$ with  $\|f_{\MS'}\|_{\Lip((\MS',\,\dm),\BS)}\le \lambda$.
\end{assumption}
\par Our aim is to prove the existence of a mapping $G:\Mc\to\KMY$ satisfying conditions (i) and (ii) of Theorem \reff{HDS-M}. We refer to $G$ as a ``core'' of the set-valued mapping $F$.

\smallskip
\par Let $T=(X,E)$ be an arbitrary finite tree. We introduce the following
\begin{definition} \lbl{ADM}{\em A mapping $\psi:X\to\MS$ is said to be {\it admissible} with respect to $T$ if for every two distinct nodes $u,v\in X$ with $u\je v$ (i.e., $u$ is joined by an edge to $v$), we have $\psi(u)\ne \psi(v)$.}
\end{definition}
\medskip
\par Let $\psi:X\to\MS$ be an admissible mapping. Then $\psi$ gives rise a {\it tree metric} $d_{T,\psi}:X\times X\to \R_+$ defined by 
\bel{DTREE}
d_{T,\psi}(u,v)=\dm(\psi(u),\psi(v))~~~\text{for every}~~ u,v\in X,~ u\je v\,.
\ee
See \rf{D-TR}.
\par Clearly, by the triangle inequality,
\bel{RH-D}
\dm(\psi(u),\psi(v))\le d_{T,\psi}(u,v)~~~
\text{for every}~~ u,v\in X\,.
\ee
\par Now define a set-valued mapping $F_{T,\psi}:X\to\Kc_m(\BS)$ by the formula
$$
F_{T,\psi}(u)=F(\psi(u)),~~~u\in X.
$$
\begin{lemma}\lbl{S-T} The set-valued mapping $F_{T,\psi}=F\circ \psi$ has a $d_{T,\psi}$-Lipschitz selection $f:X\to\BS$ such that
\bel{LS-FT}
\|f\|_{\Lip((X,\,d_{T,\psi}),\BS)}\le \gz\,\lambda.
\ee
\end{lemma}
\par {\it Proof.}  Let $X'\subset X$ be an arbitrary subset of $X$ with $\# X'\le \ks$, and  let $\MS'=\psi(X')$. Then
$$
\#\MS'\le \#X'\le \ks
$$
so that, by Assumption \reff{A-F}, the restriction $F|_{\MS'}$ has a $\dm$-Lipschitz selection $f_{\MS'}:\MS'\to \BS$ with  $\|f_{\MS'}\|_{\Lip((\MS',\,\dm),\BS)}\le \lambda$.
\par Let $g_{X'}:X'\to Y$ be defined by
$$
g_{X'}(u)=f_{\MS'}(\psi(u)),~~~u\in X'.
$$
Then $g_{X'}$ is a {\it selection} of the restriction $F_{T,\psi}|_{X'}$, i.e., $g_{X'}(u)\in F_{T,\psi}(u)$ for all $u\in X'$. Furthermore, for every $u,v\in X'$
$$
\|g_{X'}(u)-g_{X'}(v)\|=
\|f_{\MS'}(\psi(u))-f_{\MS'}(\psi(v))\|\le
\lambda\,\dm(\psi(u),\psi(v))
$$
so that, by \rf{RH-D},
$$
\|g_{X'}(u)-g_{X'}(v)\|\le \lambda\, d_{T,\psi}(u,v)
$$
proving that the $d_{T,\psi}$-Lipschitz seminorm of $g_{X'}$ is bounded by $\lambda$.
\par Hence, by Corollary \reff{A-12}, the set-valued mapping $F_{T,\psi}$ has a $d_{T,\psi}$-Lipschitz selection $f:X\to \BS$ satisfying inequality \rf{LS-FT}.\bx
\medskip
\par We will need the following two definitions.
\begin{definition} \lbl{L-TAP}{\em Let $x\in\MS$. The  family $\APT(x)$ consists of all triples $L=[T,a,\psi]$ where
\smallskip
\par \textbullet~ $T=(X,E)$ is a finite tree with the family of nodes $X$ and the family of edges $E$;
\smallskip
\par \textbullet~ $a\in X$ is a node of $T$;
\smallskip
\par \textbullet~ $\psi:X\to\MS$ is an admissible mapping with respect to $T$ such that $\psi(a)=x$.
\smallskip
\par We refer to each triple $L=[T,a,\psi]\in\APT(x)$ as
{\it an admissibly placed tree rooted at $a$}. We call $\APT(x)$ the family of all {\it Admissibly Placed Trees} associated with $x$.}
\end{definition}
\begin{definition}\lbl{OR-D} {\em \par Let $x\in\MS$. Given a finite tree $T=(X,E)$ and a triple $L=[T,a,\psi]\in \APT(x)$ we let $O(x;L)$ denote the subset of $\BS$ defined by
$$
O(x;L)=\{f(a):f~\,\text{is a}~~d_{T,\psi}\text{-Lipschitz selection of}~~ F_{T,\psi}~\,\text{with} ~~\|f\|_{\Lip((X,d_{T,\psi}),\BS)}\le\gz\,\lambda\}\,.
$$
}
\end{definition}
\medskip
\par We recall that a convex subset of $Y$ has dimension at most $m$ if it is contained in an affine subspace of
$Y$ of dimension at most $m$.
\begin{lemma}\lbl{OR-PR} Let $x\in\MS$ and let  $L=[T,a,\psi]\in \APT(x)$. Then $O(x;L)$ is a nonempty compact convex subset of $F(x)$ of dimension at most $m$.
\end{lemma}
\par {\it Proof.} By Lemma \reff{S-T}, the mapping $F_{T,\psi}=F\circ \psi$ has a $d_{T,\psi}$-Lipschitz selection $f:X\to\BS$ with $\|f\|_{\Lip((X,d_{T,\psi}),\BS)}\le\gz\,\lambda$. Therefore, by Definition \reff{OR-D}, $f(a)\in O(x;L)$ proving that $O(x;L)\ne\emp$.
\smallskip
\par The convexity of $O(x;L)$ directly follows from the convexity of sets $F(y)$ $(y\in \Mc)$ and Definition \reff{OR-D}. Furthermore, if $f:X\to Y$ is a selection of $F_{T,\psi}=F\circ \psi$, then $f(a)\in F(\psi(a))=F(x)$
(recall that $x=\psi(a)$, see Definition \reff{L-TAP}). 
\par Hence, $O(x;L)\subset F(x)$. This also proves that $\dim O(x;L)\le \dim F(x)\le m$.
\smallskip
\par Let us prove that $O(x;L)$ is compact whenever each set $F(y), y\in \Mc$, is. Since $O(x;L)\subset F(x)$ and $F(x)$ is a compact set, $O(x;L)$ is a bounded set. We prove that  $O(x;L)$ is closed.
\par Let $h\in\BS$, and a let $h_n\in O(x;L), n=1,2,...$ be a sequence of points converging to $h$:
\bel{HN-L1}
h=\lim_{n\to\infty}h_n\,.
\ee
We will prove that $h\in O(x;L)$.
\par By Definition \reff{OR-D}, there exists a sequence of mappings $f_n\in \Lip((X,d_{T,\psi}),\BS)$ such that
\bel{FN-OR1}
f_n(u)\in F(\psi(u)) ~~~\text{and}~~~
\|f_n\|_{\Lip((X,d_{T,\psi}),\BS)}\le\gz\,\lambda
\ee
for every $u\in X$ and $n\in\N$, and
\bel{HN-FNA1}
h_n=f_n(a),~~~n=1,2,...\,.
\ee
\par Note that $(X,d_{T,\psi})$ is a {\it finite} metric space, and each set $F(\psi(u)), u\in X$, is a finite dimensional compact subset of $Y$. Therefore, there exists a subsequence $n_k\in\N$, $k=1,2,...$, such that $(f_{n_k}(u))_{k=1}^\infty$ converges in $Y$ for every $u\in X$. Let
\bel{TF-DOR1}
\tf(u)=\lim_{k\to\infty}f_{n_k}(u), ~~~u\in X.
\ee
\par Then, by \rf{HN-L1} and \rf{HN-FNA1},
\bel{H-FA1}
h=\lim_{k\to\infty}h_{n_k}
=\lim_{k\to\infty}f_{n_k}(a)=\tf(a).
\ee
\par Since each set $F(\psi(u))$, $u\in X$, is closed, by \rf{FN-OR1} and \rf{TF-DOR1}, $\tf(u)\in F(\psi(u))$ for every $u\in X$, proving that $\tf$ is a {\it selection} of the set-valued mapping $F_{T,\psi}=F\circ \psi$ on $X$. Since each mapping $f_n:X\to\BS$ is $d_{T,\psi}$-Lipschitz with  $\|f_n\|_{\Lip((X,d_{T,\psi}),\BS)}\le\gz\,\lambda$, by \rf{TF-DOR1}, $\tf$ is $d_{T,\psi}$-Lipschitz as well, with $\|\tf\|_{\Lip((X,d_{T,\psi}),\BS)}\le\gz\,\lambda$.
\par Thus, by \rf{H-FA1} and Definition \reff{OR-D}, $h\in O(x;L)$ proving the lemma.\bx
\medskip
\par Given $x\in\MS$ let
\bel{GX-D}
G(x)=\bigcap_{L\,\in\,\APT(x)}\,O(x;L)\,.
\ee
Clearly, by Lemma \reff{OR-PR}, for every $x\in\MS$ the set 
\bel{G-KKS}
G(x)~~\text{\it is a convex compact subset of}~~F(x).
\ee
\par In the next section, we will prove that $G(x)\ne\emp$ for each $x\in\MS$ and that
\bel{GX-HD}
\dhf(G(x),G(y))\le \gz\,\lambda\,\dm(x,y)~~~
\text{for every}~~x,y\in\MS\,.
\ee
Recall that $\dhf$ denotes the Hausdorff distance between subsets of $\BS$.
\bigskip\vskip 3mm
\indent\par {\bf 5.2 Lipschitz continuity of the ``core'' with respect to the Hausdorff distance.}
\addtocontents{toc}{~~~~5.2 Lipschitz continuity of the ``core'' with respect to the Hausdorff distance.\hfill \thepage\par\VSU}
\begin{lemma}\lbl{L1} For every $x\in\MS$\,, the set $G(x)\ne\emp$.
\end{lemma}
\par {\it Proof.} We must show that
$$
\bigcap_{L\,\in\,\APT(x)}\,O(x;L)\ne\emp.
$$
See \rf{GX-D}. By Lemma \reff{OR-PR}, each $O(x;L)$ is a nonempty compact subset of the compact set $F(x)$. Therefore, it is enough to show that
\bel{O-LI2}
O(x;L_1)\cap...\cap O(x;L_N)\,\ne\emp
\ee
for every finite subcollection $\{L_1,...,L_N\}\subset\APT(x)$.
\par Let $L_1,...,L_N\in \APT(x)$ with $L_i=[T_i,a_i,\psi_i]$,~$i=1,...,N$, where each $T_i=(X_i,E_i)$ is a finite tree.
\medskip
\par We introduce a procedure for gluing the finite trees $T_i=(X_i,E_i)$,~$i=1,...,N$, together. Recall that $X_i$ here denotes the set of nodes of $T_i$, and $E_i$ denotes the set of edges of $T_i$. By passing to isomorphic copies of the $T_i$, we may assume that the sets $X_i$ are pairwise disjoint. Then we form a finite tree $T^+=(X^+, E^+)$ from $T_1,...,T_N$ by identifying together all the nodes $a_1,...,a_N$.  We spell out details below. \smallskip
\par For each $i$, we write $J_i$ to denote the set of all the neighbors of $a_i$ in $T_i$. Also, we write $X'_i$ to denote the set $X_i\setminus\{a_i\}$, and we write $E'_i$ to denote all the edges in $T_i$ that join together points of $X'_i$ (i.e. not including $a_i$ as an endpoint).\smallskip
\par We introduce a new node $a^+$ distinct from all the nodes of all the $T_i$.\smallskip
\par The finite tree $T^+=(X^+,E^+)$ is then defined as follows. The nodes in $X^+$ are all the nodes in all the $X'_i$, together with the single node $a^+$. The edges in $E^+$ are all the edges belonging to any of the $E'_i$, together with edges joining $a^+$ to all the nodes in all the $J_i$. One checks easily that $T^+$ is a finite tree. We say that $T^+$ arises by ``gluing together the $T_i$ by identifying the $a_i$''.\smallskip
\par Note that $T^+$ contains an isomorphic copy of each $T_i$ as a subtree; the relevant isomorphism $\varphi_i$ carries the node $a_i$ of $T_i$ to the node $a^+$ of $T^+$, and $\varphi_i$ is the identity on all other nodes of $T_i$. Each edge $[ab]$ of the tree $T_i$ is carried to the edge $[\varphi_i(a)\,\varphi_i(b)]$ of $T^+$.
\smallskip
\par This concludes our discussion of the gluing of trees $T_i$.
\medskip
\par  We define a map $\psi^+:X^+\to\Mc$ by setting
\bel{WP-1}
\psi^+(a^+)=x
\ee
and
\bel{WP-2}
\psi^+(b)=\psi_i(b)~~~\text{for all}~~~b\in X'_i=X_i\setminus\{a_i\},~i=1,...,N.
\ee
\par One checks that $\psi^+$ is an admissible map, and $\psi^+(a^+)=x$. Thus , $L^+=[T^+,a^+,\psi^+]$ belongs to $\APT(x)$. Consequently, by Lemma \reff{S-T}, there exists a $d_{T^+,\psi^+}$-Lipschitz selection $f^+$ of $F\circ \psi^+$ with $d_{T^+,\psi^+}$-Lipschitz  seminorm $\le \gz\,\lambda$. (We recall that the metric $d_{T^+,\psi^+}$ is determined by formula \rf{DTREE}.)
\smallskip
\par The map
$$
f_i(b)=\left \{
\begin{array}{ll}
f^+(b),& \text{if}~~b\in X_i\setminus\{a_i\},\vspace{2mm}\\
f^+(a^+),& \text{if}~~b=a_i,
\end{array}
\right.
$$
is a $d_{T_i,\psi_i}$-Lipschitz selection of $F\circ \psi_i$ with $d_{T_i,\psi_i}$-Lipschitz seminorm $\le\gz\,\lambda$, therefore $$f^+(a^+)\in O(x;L_i)~~~\text{for each}~~~i=1,...,N.$$
\par Thus, \rf{O-LI2} holds, completing the proof of Lemma \reff{L1}.\bx
\smallskip
\par We know that the affine dimension of each set $F(x)$ is at most $m$. Since $G(x)\subset F(x)$, the same is true for each set $G(x)$, $x\in \MS$. This observation, Lemma \reff{L1} and statement \rf{G-KKS} imply that $G$ maps the metric space $\MS$ into the family $\KM$.
\smallskip
\par We are in a position to prove inequality \rf{GX-HD}.
\begin{lemma}\lbl{HD-Z} For every $x,y\in\MS$ the following inequality
$$
\dhf(G(x),G(y))\le \gz\,\lambda\,\dm(x,y)
$$
holds.
\end{lemma}
\par {\it Proof.} We may suppose $x\ne y$, else the desired conclusion is obvious. Let us prove that
\bel{R-IN}
I=G(x)+\gz\,\lambda\,\dm(x,y)\,\BY\supset G(y)\,.
\ee
Recall that by $\BY=\BL(0,1)$ we denote the closed unit ball in $\BS$.
\par If we can prove that, then by interchanging the roles of $x$ and $y$ we obtain also
$$
G(y)+\gz\,\lambda\,\dm(x,y)\,\BY\supset G(x)\,.
$$
These two inclusions tell us that $\dhf(G(x),G(y))\le \gz\,\lambda\,\dm(x,y)$, proving the lemma.\smallskip
\par Let us prove \rf{R-IN}. By definition,
$$
I=\left[\,\bigcap_{L\,\in\,\APT(x)}\,O(x;L)\right]+
\gz\,\lambda\,\dm(x,y)\,\BY\,.
$$
See \rf{GX-D}. We will check that
\bel{Z2}
\left[\,\bigcap_{L\,\in\,\APT(x)}\,O(x;L)\right]+
\gz\,\lambda\,\dm(x,y)\,\BY=\bigcap\,\left\{\,
\left[O(x;L_1)\cap...\cap O(x;L_N)\right]+\gz\,\lambda\,\dm(x,y)\,\BY\right\},
\ee
where the first intersection of the right-hand side is taken over all finite sequences $L_1,...,L_N$ of elements of $\APT(x)$.
\par Indeed, the left-hand side of \rf{Z2} is obviously contained in the right-hand side. Conversely, let $\xi$ belong to the right-hand side of \rf{Z2}. Then any finite subcollection of the compact sets
$$
K_L=\{\eta\in\BY: \xi-\gz\,\lambda\,\dm(x,y)\,\eta\in O(x;L)\}
$$
has nonempty intersection. (The above sets are compact because $O(x;L)$ is compact.)
\par Therefore,
$$
\bigcap_{L\in\APT(x)}\,K_L\ne\emp,
$$
proving that $\xi$ belongs to the left-hand side of \rf{Z2}. The proof of \rf{Z2} is complete.\smallskip
\par Thanks to \rf{Z2}, our desired inclusion \rf{R-IN} will follow if we can show that
\bel{Z3}
\left[O(x;L_1)\cap...\cap O(x;L_N)\right]+\gz\,\lambda\,\dm(x,y)\,\BY\supset G(y)
\ee
for any $L_1,...,L_N\in\APT(x)$. Then the proof of Lemma \reff{HD-Z} is reduced to the task of proving \rf{Z3}.
\smallskip
\par Let $L_i=[T_i,a_i,\psi_i]$ where $T_i=(X_i,E_i)$. Then $a_i$ is a node of the tree $T_i$, $i=1,...,N$. We introduce a new node $a^+$ and form the tree $T^+=(X^+,E^+)$ as in the proof of Lemma \reff{L1}. Thus $T^+$ arises by {\it gluing together the trees $T_i$ by identifying the $a_i$}.
\par We also introduce an admissible map $\psi^+:X^+\to\Mc$ as in the proof of Lemma \reff{L1}, see \rf{WP-1} and \rf{WP-2}.
\smallskip
\par We now introduce a new node $\ta$ not present in $T^+$. We define a new tree $\tT=(\tX,\tE)$ as follows.
\smallskip
\par \textbullet~ The nodes in $\tX$ are the nodes in $X^+$, together with the new node $\ta$.
\smallskip
\par \textbullet~ The edges in $\tE$ are the edges in $E^+$, together with a single edge joining $\ta$ to $a^+$.
\medskip
\par We define a map $\tpsi:\tT\to\Mc$ by setting
$$
\tpsi=\psi^+~~~\text{on}~~~T^+,~~\tpsi(\ta)=y.
$$
Then one checks that $\tT=(\tX,\tE)$ is a tree and $\tpsi(\ta)=y$. Furthermore, since $\psi^+$ is admissible on $T^+$ and $x\ne y$, the mapping $\tpsi$ is admissible on $\tT$.
\smallskip
\par Let $\tL=[\tT,\ta,\tpsi]$, and let $\eta\in G(y)$. Then, by definition \rf{GX-D}, $\eta\in O(y;\tL)$ so that there exists a $d_{\tT,\tpsi}$-Lipschitz selection $\tf$ of $F\circ \tpsi$, with $d_{\tT,\tpsi}$-Lipschitz seminorm $\le\gz\,\lambda$, satisfying $\tf(\ta)=\eta$. See Definition \reff{OR-D}. (We also recall that the metric $d_{\tT,\tpsi}$ is defined by formulae \rf{DTREE} and \rf{D-TR}.)
\smallskip
\par Restricting this $\tf$ to $T^+$ and arguing as in the proof of Lemma \reff{L1}, we see that                                                      
$$
\tf(a^+)\in
O(x;L_1)\cap...\cap O(x;L_N).
$$
\par On the other hand, our Lipschitz bound for $\tf$ gives
$$
\|\tf(a^+)-\eta\|=\|\tf(a^+)-\tf(\ta)\|
\le\gz\,\lambda\,\rho(\tpsi(a^+),\tpsi(\ta))
=\gz\,\lambda\,\rho(x,y).
$$
Then,
$$
\eta\in \left[O(x;L_1)\cap...\cap O(x;L_N)\right]+\gz\,\lambda\,\dm(x,y)\,\BY
$$
proving \rf{Z3}.\bx
\medskip
\par The proof of Theorem \reff{HDS-M} is complete.\bx
\medskip
\par We turn to the final step of the proof of Theorem \reff{KSHARP}. The following selection theorem is a special case of \cite[Theorem 1.2]{S04}.
\begin{theorem}\lbl{ST-P} Let $Y$ be a Banach space, and let $m\ge 1$. Then there exists a map $\ST:\KM\to Y$ such that
\smallskip
\par $(\alpha)$~~ $\ST(K)\in K$ for all $K\in\KM$
\smallskip
\par\noindent and
\par $(\beta)$~~ $\|\ST(K)-\ST(K')\|\le C(m)\cdot \dhf(K,K')$~~ for all $K,K'\in\KM$.
\smallskip
\par Here $C(m)$ depends only on $m$.
\end{theorem}
\smallskip
\par We refer to $\ST(K)$ as the {\it ``Steiner-type point'' of $K$}. In the special case $Y=\RM$, we can take $\ST(K)$ to be the {\it Steiner point} of $K$. Recall that the Steiner point of $K$ may be defined as the limit as $R\to \infty$ of the barycenter of $K+B(R)$, where ``$+$'' denotes Minkowski sum, and $B(R)$ is the standard Euclidean ball of radius $R$ about $0$. For general $Y$, there is no simple description of the ``Steiner-type point'' $\ST(K)$ in \cite{S04}.
\medskip
\par To construct the Lipschitz selection $f$ and establish Theorem \reff{KSHARP}, we just set
$$
f(x)=\ST(\CRF(x))~~~\text{for}~~~x\in\Mc,
$$
where $\CRF$ is the core defined by \rf{GX-D}. Since $G(x)\in\KM$ for each $x\in\Mc$, the function $f$ is well defined on $\Mc$.
\par By part (i) of Theorem \reff{HDS-M} and part $(\alpha)$ of Theorem \reff{ST-P},
$$
f(x)=\ST(\CRF(x))\in \CRF(x)\subset F(x)~~~\text{for}~~~x\in\Mc.
$$
On the other hand, part (ii) of Theorem \reff{HDS-M} and part $(\beta)$ of Theorem \reff{ST-P} imply that
$$
\|f(x)-f(y)\|=\|\ST(\CRF(x))-\ST(\CRF(y))\|\le C(m)\cdot\dhf(\CRF(x),\CRF(y))
\le C(m)\cdot \gz\,\lambda\,\rho(x,y)
$$
for all $x,y\in\Mc$. Thus, $f$ is a Lipschitz selection of $F$ with Lipschitz seminorm at most $C(m)\cdot \gz\,\lambda$. Recalling that $C(m)$ and $\gz$ depend only on $m$, we conclude that Theorem \reff{KSHARP} holds.\bx
\SECT{6. Pseudometric spaces}{6}
\addtocontents{toc}{6. Pseudometric spaces \hfill\thepage\par\VST}

\indent\par In this section we prove Theorem \reff{MAIN-FP}, the Finiteness Principle for Lipschitz Selections, and Theorem \reff{FINITE}, a variant of Theorem \reff{MAIN-FP} for finite pseudometric spaces.
\par Until the end of Section 6 we write $\gone$ to denote the constant $\gone=\gone(m)$ from Theorem \reff{KSHARP}.
Everywhere in Section 6 we write $\ks=\ks(m)$ and $\gz=\gz(m)$ to denote the constants from  Corollary \reff{A-12}.
\par Let $(\Mc,\dm)$ be a pseudometric space. Recall that we say that the {\it pseudometric $\rho$ is finite} if
\bel{FIN-A}
\rho(x,y)<\infty~~~\text{for all}~~x,y\in\Mc\,.
\ee
\par On the other hand, we say that $(\Mc,\rho)$ is a {\it finite pseudometric space} if $\Mc$ contains only finitely many points.
\par Given a set-valued mapping $F:\Mc\to\CNMY$, by a {\it selection} of $F$ (not necessarily Lipschitz) we mean simply a map $f:\Mc\to Y$ such that $f(x)\in F(x)$ for all $x\in\Mc$.

\indent\par {\bf 6.1 The final step of the proof of the Finiteness Principle.}
\addtocontents{toc}{~~~~6.1 The final step of the proof of the Finiteness Principle.\hfill \thepage\par}
\smallskip
\par In this section we prove an analog of Theorem \reff{KSHARP} for pseudometric spaces.
\begin{proposition}\lbl{PM-SP} Let $(\Mc,\dm)$ be a pseudometric space satisfying \rf{FIN-A}, and let $\lambda>0$. Let $F:\MS\to\KM$ be a set-valued mapping such that for every subset $\Mc'\subset\MS$ consisting of at most $\ks$ points, the restriction $F|_{\Mc'}$ of $F$ to $\Mc'$ has a Lipschitz selection $f_{\Mc'}:\Mc'\to\BS$ with  $\|f_{\Mc'}\|_{\Lip(\Mc',\BS)}\le \lambda$.
\par Then $F$ has a Lipschitz selection $f:\Mc\to\BS$ with $\|f\|_{\Lip(\Mc,\BS)}\le \gone\lambda$.
\end{proposition}
\par {\it Proof.} A selection of $F$ may be regarded as a point of the Cartesian product
$$
\Fc=\prod_{x\in\Mc}\,F(x)\,.
$$
We endow $\Fc$ with the product topology. Then, by Tychonoff's theorem, $\Fc$ is compact because each $F(x)$ is compact.
\par For $\ve>0$ and $x,y\in\Mc$, let
$$
\rho_\ve(x,y)=\left \{
\begin{array}{ll}
\rho(x,y)+\ve,& \text{if}~~x\ne y,\vspace*{2mm}\\
0,& \text{if}~~x=y.
\end{array}
\right.
$$
\par Then $(\Mc,\rho_\ve)$ is a metric space. For any $\Mc'\subset\Mc$ with $\#\Mc'\le \ks$ there exists a selection of $F|_{\Mc'}$ with $\rho$-Lipschitz seminorm $\le \lambda$, hence with $\rho_\ve$-Lipschitz seminorm $\le \lambda$. By Theorem \reff{KSHARP}, $F$ has a selection with $\rho_\ve$-Lipschitz seminorm $\le\gone\lambda$.
\par Let $\Sel(\ve)$ be the set of all selections of $F$ with $\rho_\ve$-Lipschitz seminorm at most $\gone\lambda$. Then $\Sel(\ve)$ is a closed subset of $\Fc$. We have just seen that $\Sel(\ve)$ is nonempty. Because
$$
\Sel(\ve)\subset \Sel(\ve')~~~\text{for}~~~\ve<\ve',
$$
it follows that
$$
\Sel(\ve_1)\cap\Sel(\ve_2)\cap...\cap\Sel(\ve_N)\ne\emp
$$
for any $\ve_1,\ve_2,...,\ve_N>0$.
\par Because $\Fc$ is compact and each $\Sel(\ve)$ is closed in $\Fc$, it follows that
$$
\bigcap_{\ve>0}\Sel(\ve)\ne\emp\,.
$$
Furthermore, any $f\in\capbig\{\Sel(\ve):\ve>0\}$ is a selection of $F$ with $\rho$-Lipschitz seminorm $\le\gone\lambda$.
\par The proof of Proposition \reff{PM-SP} is complete.\bx
\medskip
\par {\it Proof of Theorem \reff{MAIN-FP}.} Suppose that $\rho$ is a finite pseudometric, i.e., condition \rf{FIN-A} holds.
\par Let $\Mc'$ be an arbitrary subset of $\Mc$ consisting of at most $\ks$ points. (Note that, by definitions \rf{NMY-1}, \rf{NMY-2} and \rf{KS-DF}, \rf{LS-DF}, the constant $\ks>N(m,Y)$ for every $m\in\N$.) Then, by the theorem's hypothesis, for every set $S\subset\Mc'$ with $\#S\le\FN$, the restriction $F|_S$  has a Lipschitz selection $f_{S}:S\to\BS$ with $\|f_S\|_{\Lip(S,\BS)}\le \lambda$. Hence, by Theorem \reff{KS-N}, the restriction $F|_{\Mc'}$ of $F$ to $\Mc'$ has a Lipschitz selection $f_{\Mc'}:\Mc'\to\BS$ whose seminorm satisfies $\|f_{\Mc'}\|_{\Lip(\Mc',\BS)}\le \gamma\lambda$ where $\gamma$ is a constant depending only on $m$ and $\#\Mc'$. Since $\#\Mc'\le \ks$ and $\ks$ depends only on $m$, the constant $\gamma$ depends only on $m$ as well.
\smallskip
\par Hence, by Proposition \reff{PM-SP}, $F$ has a Lipschitz selection $f:\Mc\to Y$ with
$\|f\|_{\Lip(\Mc,Y)}\le\gone\gamma\lambda$.
Recall that $\gone$ is a constant depending only on $m$.
\par This completes the proof of Theorem \reff{MAIN-FP} for the case of a {\it finite} pseudometric $\rho$.
\smallskip
\par To pass to the general case in which $\rho(x,y)$ may take the value $+\infty$ is an easy exercise. We define an equivalence relation on $\Mc$ by calling $x$ and $y$ equivalent when $\rho(x,y)$ is finite. On each equivalence class we produce a Lipschitz selection of $F$, with controlled Lipschitz seminorm, by invoking the known case of Theorem \reff{MAIN-FP} in which all distances are finite. By combining those Lipschitz selections into a single function defined on the union of all the equivalence classes, we obtain the desired Lipschitz selection of $F$.
Details are spelled out in \cite{FS-2017}.
\par The proof of Theorem \reff{MAIN-FP} is complete. \bx
\medskip


\indent\par {\bf 6.2 Finite pseudometric spaces.}
\addtocontents{toc}{~~~~6.2 Finite pseudometric spaces.\hfill \thepage\par}
\smallskip

\par In this section we prove a variant of our main result, Theorem \reff{MAIN-FP}, related to the case of {\it finite} pseudometric spaces. As we have noted in the Introduction, for the case of the trivial distance function $\dm\equiv 0$ defined on a finite pseudometric space, Theorem \reff{FINITE} below agrees with the classical Helly's Theorem \cite{DGK} (up to the values of $N(m,Y)$ and the optimal finiteness constant for $\dm\equiv 0$ (see \rf{DM-NMY})).
\begin{theorem}\lbl{FINITE} Let $(\Mc,\dm)$ be a finite pseudometric space, and let $F:\Mc\to\CNMY$ be a set-valued mapping from $\Mc$ into the family $\CNMY$ of all convex subsets of $Y$ of affine dimension at most $m$. Let $\lambda$ be a positive real number.
\par Suppose that for every subset $\Mc'\subset\MS$ consisting of at most $\FN$ points, the restriction $F|_{\Mc'}$ of $F$ to $\Mc'$ has a Lipschitz selection $f_{\Mc'}$ with Lipschitz  seminorm $\|f_{\Mc'}\|_{\Lip(\Mc',\BS)}\le \lambda$.
\par Then $F$ has a Lipschitz selection $f$ with Lipschitz  seminorm $\|f\|_{\Lip(\MS,\BS)}\le \gamma\,\lambda$.
\par Here, $\gamma$ depends only on $m$.
\end{theorem}
\par Our proof of this result relies on an analog of Proposition \reff{PM-SP} for a {\it finite} pseudometric space $(\Mc,\rho)$ and a set-valued mapping $F:\Mc\to\CNMY$. See Proposition \reff{FN-PM} below.
\par We will need three auxiliary lemmas.
\begin{lemma}\lbl{FM-B} Let $\lambda>0$ and let $(\Mc,\dm)$ be a finite metric space. Let $F$ be a set-valued mapping on $\Mc$ which to every $x\in\Mc$ assigns a nonempty convex bounded subset of $Y$ of dimension at most $m$.
\par Suppose that for every subset $\Mc'\subset\MS$ with $\#\Mc'\le \ks$, the restriction $F|_{\Mc'}$ of $F$ to $\Mc'$ has a Lipschitz selection $f_{\Mc'}:\Mc'\to\BS$ with  $\|f_{\Mc'}\|_{\Lip(\Mc',\BS)}\le \lambda$.
\par Then $F$ has a Lipschitz selection $f:\Mc\to\BS$ with $\|f\|_{\Lip(\Mc,\BS)}\le 2\gone\lambda$. Here, $\gone$ is as in Proposition \reff{PM-SP}.
\end{lemma}
\par {\it Proof.} We introduce a new set-valued mapping on $\Mc$ defined by
$$
\tF(x)=(F(x))^{\cl}~~~~\text{for all}~~~x\in\Mc.
$$
Here the sign $\cl$ denotes the closure of a set in $Y$.
\par Since the sets $F(x)$, $x\in\Mc$, are finite dimensional and bounded, each set $\tF(x)$ is compact so that $\tF:\Mc\to\KM$. Furthermore, since $F(x)\subset\tF(x)$ on $\Mc$, the mapping $\tF$ satisfies the hypothesis of Proposition \reff{PM-SP}.
\par By this proposition, there exists a mapping $\tf:\Mc\to Y$ such that
\bel{TF-CL}
\tf(x)\in\tF(x)=(F(x))^{\cl}~~~\text{for all}~~x\in\Mc,
\ee
and
\bel{FW-11}
\|\tf(x)-\tf(y)\|\le \gone\,\lambda\,\rho(x,y) ~~~\text{for all}~~x,y\in\Mc.
\ee
\par Since $\Mc$ is a {\it finite} metric space, the following quantity
\bel{DL-M}
\delta=\gone\,\lambda\,\min_{x,y\in\Mc,\,x\ne y}\rho(x,y)
\ee
is positive. Therefore, by \rf{TF-CL}, for each $x\in\Mc$ there exists a point $f(x)\in F(x)$ such that
$$
\|f(x)-\tf(x)\|\le \delta/2\,.
$$
\par Thus $f:\Mc\to Y$ is a selection of $F$ on $\Mc$. Let us estimate its Lipschitz seminorm. For every $x,y\in\Mc$ (distinct), by \rf{FW-11} and \rf{DL-M},
$$
\|f(x)-f(y)\|\le \|f(x)-\tf(x)\|+\|\tf(x)-\tf(y)\|+\|\tf(y)-f(y)\|
\le \delta/2+
\gone\lambda\,\rho(x,y)+\delta/2\le 2\gone\lambda\,\rho(x,y).
$$
Hence, $\|f\|_{\Lip(\Mc,\BS)}\le 2\gone\lambda$, and the proof of the lemma is complete.\bx
\medskip
\par The second auxiliary lemma provides additional properties of sets $\Gamma_\ell$ defined in Section 3.1 (see \rf{BC-3} and Definition \reff{GML-D}). We will need these properties in the proof of Lemma \reff{P-2} below.
\begin{lemma}\lbl{G-AP} Let $(\Mc,\rho)$ be a finite pseudometric space satisfying \rf{FIN-A}. Let $\ell\ge 0$ and let $F:\Mc\to\CNMY$. Suppose that for every subset $\Mc'\subset\Mc$ with $\#\Mc'\le k_{\ell+1}$ the restriction $F|_{\Mc'}$ of $F$ to $\Mc'$ has a Lipschitz selection $f_{\Mc'}:\Mc'\to Y$ with  $\|f_{\Mc'}\|_{\Lip(\Mc',Y)}\le\lambda$.
\par Let $x_0\in\Mc$, $\xi_0\in\GL(x_0)$, and let $1\le k\le \ell+1$. Let $S$ be a subset of $\Mc$ with $\#S=k$ containing $x_0$.
\par Then there exists a mapping $f^S:S\to Y$ such that
\smallskip
\par (a) $f^S(x_0)=\xi_0$.
\smallskip
\par (b) $f^S(y)\in \Gamma_{\ell+1-k}(y)$ for all $y\in S$.
\smallskip
\par (c) $\|f^S\|_{\Lip(S,Y)}\le 3^k\lambda$.
\end{lemma}
\par {\it Proof.} We recall that the sequence of positive integers $k_\ell$ is defined by the formula \rf{KEL}.
\par We proceed by induction on $k$. For $k=1$, we have $S=\{x_0\}$, and we can just set $f^S(x_0)=\xi_0$.
\par For the induction step, we fix $k\ge 2$ and suppose the lemma holds for $k-1$; we then prove it for $k$. Thus, let $\xi_0\in\GL(x_0)$, $x_0\in S$, $\# S=k\le \ell+1$.
\par Set $\hS=S\setminus\{x_0\}$. We pick $\hx_0\in\hS$ to minimize $\rho(\hx_0,x_0)$, and we pick $\hxi_0\in\GLO(\hx_0)$ such that $\|\hxi_0-\xi_0\|\le \lambda\,\rho(\hx_0,x_0)$. (See Lemma \reff{G-AB} (b).) For $y\in\hS$ we have $\rho(y,x_0)\ge \rho(\hx_0,x_0)$, hence
\bel{R-V}
\rho(y,\hx_0)+\rho(\hx_0,x_0)\le[\rho(y,x_0)
+\rho(x_0,\hx_0)]+\rho(\hx_0,x_0)\le 3\rho(y,x_0).
\ee
\par By the induction hypothesis, there exists
$\hf:\hS\to Y$ such that
\smallskip
\par $(\ha)$ $\hf(\hx_0)=\hxi_0$.
\smallskip
\par $(\hb)$ $\hf(y)\in \Gamma_{(\ell-1)+1-(k-1)}(y)=\Gamma_{\ell+1-k}(y)$ for all $y\in \hS$.
\smallskip
\par $(\hc)$ $\|\hf\|_{\Lip(\hS,Y)}\le 3^{k-1}\lambda$.
\smallskip
\par We now define $f:S\to Y$ by setting
$$
f(y)=\hf(y)~~~\text{for}~~~y\in\hS;~~~~f(x_0)=\xi_0.
$$
\par Then $f$ obviously satisfies $(a)$ and $(b)$. To see that $f$ satisfies (c), we first recall $(\hc)$; thus it is enough to check that
$$
\|f(y)-f(x_0)\|\le 3^k\lambda\,\rho(y,x_0)
$$
for $y\in\hS$, i.e.,
$$
\|\hf(y)-\xi_0\|\le 3^k\lambda\,\rho(y,x_0)~~~\text{for}~~~y\in\hS.
$$
\par However, for $y\in\hS$ we have
$$
\|\hf(y)-\xi_0\|\le \|\hf(y)-\hxi_0\|+\|\hxi_0-\xi_0\|=
\|\hf(y)-\hf(\hx_0)\|+\|\hxi_0-\xi_0\|\le
3^{k-1}\lambda\,\rho(y,\hx_0)+\lambda\,\rho(\hx_0,x_0),
$$
thanks to $(\hc)$ and the definition of $\hxi_0$.
\par Therefore,
$$
\|\hf(y)-\xi_0\|\le
3^{k-1}\lambda\,[\rho(y,\hx_0)+\rho(\hx_0,x_0)]\le 3^{k}\lambda\,\rho(y,x_0),
$$
by \rf{R-V}.
\par Thus, $f$ satisfies $(a)$, $(b)$, $(c)$, completing our induction.\bx
\medskip
\par We turn to the last auxiliary lemma. Let
\bel{KNN}
\tl=\ks~~~~\text{and let}~~~~\kn=k_{\tl+1}
\ee
where $k_\ell=(m+2)^\ell$, see \rf{KEL}.
\begin{lemma}\lbl{P-2} Let $(\Mc,\dm)$ be a finite pseudometric space satisfying \rf{FIN-A}, and let $x_0\in\Mc$ and $\lambda>0$.
\par Let $F:\MS\to\CNMY$ be a set-valued mapping such that for every subset $\Mc'\subset\MS$ consisting of at most $\kn$ points, the restriction $F|_{\Mc'}$ of $F$ to $\Mc'$ has a Lipschitz selection $f_{\Mc'}:\Mc'\to\BS$ with  $\|f_{\Mc'}\|_{\Lip(\Mc',\BS)}\le \lambda$.
\par Then there exists a point $\xi_0\in F(x_0)$ such that the following statement holds: For every subset $S\subset\Mc$ with $\#S\le \ks$, there exists a mapping $f_S:S\to Y$ with $\|f_S\|_{\Lip(S,Y)}\le C\lambda$ such that
\bel{KN-F}
\|f_S(x)-\xi_0\|\le C\lambda\,\rho(x,x_0) ~~~\text{for every}~~~x\in S,
\ee
and
\bel{FNB-4}
f_S(x)\in F(y)+\lambda\,\rho(x,y)\,B_Y~~~\text{for every}~~~x\in S, y\in\Mc\,.
\ee
Here $C$ is a constant depending only on $m$.
\end{lemma}
\par {\it Proof.} By the lemma's hypothesis, \rf{KNN} and
by Lemma \reff{G-AB} {\it(a)},
$$
\Gamma_{\tl}(x)\ne\emp ~~~\text{for every}~~~x\in\Mc\,.
$$
\par Let $\xi_0\in \Gamma_{\tl}(x_0)$. By \rf{GL-F2},
$$
\xi_0\in \Gamma_{\tl}(x_0)\subset F(x_0).
$$
\par Let $S\subset\Mc$, $\#S\le \ks$. Let $\tS=S\cup\{x_0\}$ and let $k=\#\tS=\#(S\cup\{x_0\})\,.$
Then
$$
1\le k\le\#S+1\le \ks+1=\tl+1.
$$
Therefore, by Lemma \reff{G-AP}, there exists a mapping $f^{\tS}:\tS\to Y$ with $\|f^{\tS}\|_{\Lip(\tS,Y)}\le 3^k\lambda$ such that $f^{\tS}(x_0)=\xi_0$ and
$$
f^{\tS}(x)\in \Gamma_{\tl+1-k}(x)~~~\text{for all}~~~x\in\tS.
$$
Recall that $k\le\tl+1=\ks+1$ so that
$$
\|f^{\tS}\|_{\Lip(\tS,Y)}\le C\lambda
$$
with $C=3^{\ks+1}$. Since $\ks$ depends only on $m$, the constant $C$ depends only on $m$ as well.
\par Hence, by \rf{BC-5},
\bel{FG-12}
f^{\tS}(x)\in \Gamma_{\tl+1-k}(x)\subset \Gamma_{0}(x)~~~\text{for every}~~~x\in\tS.
\ee
\par Let
$$
f_S=f^{\tS}|_S\,.
$$
\par Then $\|f_S\|_{\Lip(S,Y)}\le \|f^{\tS}\|_{\Lip(\tS,Y)}\le C\lambda$. Moreover, by \rf{FG-12},
\bel{G-27}
f_S(x)\in \Gamma_{0}(x)~~~\text{for all}~~~x\in S.
\ee
\par Since $\|f^{\tS}\|_{\Lip(\tS,Y)}\le C\lambda$ and $x_0\in\tS$,
$$
\|f_S(x)-\xi_0\|=\|f^{\tS}(x)-f^{\tS}(x_0)\|\le C\lambda\,\rho(x,x_0) ~~~\text{for every}~~~x\in S.
$$
Furthermore, by \rf{GM-ZR} and \rf{G-27}, for every
$x\in S$
$$
f_S(x)\in \Gamma_{0}(x)=\bigcap_{y\in\Mc}
\left(F(y)+\lambda\,\rho(x,y)\,B_Y\right)
$$
so that
$$
f_S(x)\in F(y)+\lambda\,\rho(x,y)\,B_Y~~~\text{for every}~~~x\in S, y\in\Mc\,.
$$
\par The proof of the lemma is complete.\bx
\begin{proposition}\lbl{FN-PM} Let $(\Mc,\rho)$ be a finite pseudometric space satisfying \rf{FIN-A}, and let $\lambda>0$.
\par Let $F:\MS\to\CNMY$ be a set-valued mapping such that for every subset $\Mc'\subset\MS$ with $\#\Mc'\le \kn$, the restriction $F|_{\Mc'}$ of $F$ to $\Mc'$ has a Lipschitz selection $f_{\Mc'}:\Mc'\to\BS$ with  $\|f_{\Mc'}\|_{\Lip(\Mc',\BS)}\le \lambda$.
\par Then $F$ has a Lipschitz selection $f:\Mc\to\BS$ with $\|f\|_{\Lip(\Mc,\BS)}\le \gamma_2\lambda$ where $\gamma_2$ is a constant depending only on $m$.
\end{proposition}
\par {\it Proof.} Let $x_0\in\Mc$. By Lemma \reff{P-2}, there exists a point $\xi_0\in F(x_0)$ such that for every
set $S\subset\Mc$ with $\#S\le\ks$ there exists a mapping $f_S:S\to Y$ with $\|f_S\|_{\Lip(S,Y)}\le C\lambda$ such that \rf{KN-F} and \rf{FNB-4} hold. Here $C$ is a constant depending only on $m$.
\par We introduce a new set-valued mapping $\tF:\Mc\to \CNMY$ by letting
\bel{TF-W}
\tF(x)=\left(\,\bigcap_{y\in\Mc}
\left[F(y)+\lambda\,\rho(x,y)\,B_Y\right]\right)
\bigcap\,B_Y(\xi_0,C\lambda\,\rho(x,x_0)),
~~~x\in \Mc\,.
\ee
\par By Lemma \reff{P-2} and definition \rf{TF-W}, for every set $S\subset\Mc$ consisting of at most $\ks$ points the restriction $\tF|_{S}$ of $\tF$ to $S$ has a Lipschitz selection $f_{S}:S\to\BS$ with  $\|f_{S}\|_{\Lip(S,Y)}\le C\lambda$. In particular, $\tF(x)\ne\emp$ for every $x\in
\Mc$.
\smallskip
\par Let us introduce a binary relation ``$\sim$'' on $\MS$ by letting
$$
x\sim y~~~\Longleftrightarrow~~~\rho(x,y)=0\,.
$$
Clearly, ``$\sim$'' satisfies the axioms of an equivalence relation, i.e., it is reflexive, symmetric and transitive. Given $x\in\MS$, by $[x]=\{y\in\MS:~y\sim x\}$ we denote the equivalence class of $x$. Let
$$
[\Mc]=\Mc\,/\sim\,\,=\,\{\,[x]: x\in\MS\,\}
$$
be the corresponding quotient set of $\MS$ by ``$\sim$'', i.e.,  the family of all equivalence classes of $\MS$ by ``$\sim$''. Finally, given an equivalence class $U\in[\Mc]$ let us choose a point $w_U\in U$ and put
$$
W=\{w_U:U\in[\Mc]\}.
$$
\par Clearly, $(W,\rho)$ is a {\it finite metric space}. Let
\bel{HF-7}
\hF=\tF|_W.
\ee
Then, by \rf{TF-W} and \rf{HF-7}, $\hF$ is a set-valued mapping defined on a finite metric space which takes values in the family of all nonempty convex {\it bounded} subsets of $Y$ of dimension at most $m$. Furthermore, this mapping satisfies the hypothesis of Lemma \reff{FM-B} with $C\lambda$ in place of $\lambda$.
\par Therefore, by this lemma, there exists a Lipschitz selection $\hf:W\to Y$ of $\hF$ on $W$ with
$$
\|\hf\|_{\Lip(W,Y)}\le 2\gone\,C\lambda =\gamma_2\lambda.
$$
Here $\gamma_2=2\gone C$ is a constant depending only on $m$ (because $\gone$ and $C$ depend on $m$ only).
\smallskip
\par We define a mapping $f:\Mc\to Y$ by letting
$$
f(x)=\hf(w_{[x]}), ~~~~~x\in\Mc.
$$
\par Then $f$ is a {\it selection of $F$} on $\Mc$. Indeed, let $x\in\Mc$. Since $\hf$ is a selection of $\hF=\tF|_W$, and $w_{[x]}\in W$,
$$
f(x)=\hf(w_{[x]})\in \tF(w_{[x]})
$$
so that, by \rf{TF-W},
$$
f(x)\in \tF(w_{[x]})\subset F(x)+\lambda\,\rho(w_{[x]},x)\,B_Y.
$$
But $w_{[x]}\sim x$ so that $\rho(w_{[x]},x)=0$, proving that $f(x)\in F(x)$.
\par Let us prove that $\|f\|_{\Lip(\Mc,Y)}\le \gamma_2\lambda$, i.e.,
\bel{LN-F2}
\|f(x)-f(y)\|\le \gamma_2\lambda\,\rho(x,y)~~~\text{for all}~~~x,y\in\Mc.
\ee
\par In fact, since $\|\hf\|_{\Lip(W,Y)}\le \gamma_2\lambda$,
\be
\|f(x)-f(y)\|&=&\|\hf(w_{[x]})-\hf(w_{[y]})\|\le \gamma_2\lambda\,\rho(w_{[x]},w_{[y]})\nn\\
&\le& \gamma_2\lambda\,(\rho(w_{[x]},x)+\rho(x,y)+
\rho(y,w_{[y]}))=
\gamma_2\lambda\,\rho(x,y),\nn
\ee
proving \rf{LN-F2}.
\par The proof of Proposition \reff{FN-PM} is complete.\bx
\medskip
\par {\it Proof of Theorem \reff{FINITE}}. We prove this theorem following the scheme of the proof of Theorem \reff{MAIN-FP}. In particular, to study a
pseudometric $\rho$ that takes only finite values, we use  Proposition \reff{FN-PM} and the constant $\kn$ rather than Proposition \reff{PM-SP} and $\ks$ respectively.
\par We note that \cite[Remark 1.3]{S02} implies a variant of Theorem \reff{KS-N} for the case of a finite pseudometric space $(\wtM,\trh)$ and a set-valued mapping $\tF$ with convex (not necessarily compact) images $\tF(x)$, $x\in\wtM$, of dimension at most $m$.
\par As in the proof of Theorem \reff{MAIN-FP}, the passage  from the case of finite pseudometrics $\rho:\Mc\times\Mc\to\R_+$ to the general case of an arbitrary pseudometric $\rho:\Mc\times\Mc\to\R_+\cup\{+\infty\}$ is an easy exercise.  \bx

\SECT{7. Further results and comments}{7}
\addtocontents{toc}{7. Further results and comments \hfill\thepage\par\VST}
\smallskip

\par \textbullet~ {\bf Generalization of the finiteness principle: set-valued mappings with closed images.}\smallskip
\par In Theorem \reff{MAIN-FP} we prove the finiteness principle for set-valued mappings $F$ whose values are convex {\it compact} sets with affine dimension bounded by $m$. The following result shows that this family of sets can be slightly extended.
\begin{theorem}\lbl{ST-KL} Theorem \reff{MAIN-FP} holds provided the requirement $F:\MS\to\KM$ in its formulation is replaced with the following one: for every $x\in\Mc$ the set $F(x)$ is a closed convex subset of $Y$ of dimension at most $m$, and there exists $x_0\in\Mc$ such that $F(x_0)$ is bounded.
\end{theorem}
\par For the proof of this statement we refer the reader to \cite[p. 74]{FS-2017}.
\smallskip
\par Theorem \reff{ST-KL} implies the following result.
\begin{theorem} Theorem \reff{MAIN-FP} holds provided the requirement $F:\MS\to\KM$ in its formulation is replaced with $F:\MS\to\KMY\cup\AM$.
\par Here $\AM$ denotes the family of all affine subspaces of $Y$ of dimension at most $m$.
\end{theorem}
\par {\it Proof.} The result follows from \cite{S04} whenever $F:\MS\to\AM$, and from Theorem \reff{ST-KL} whenever there exists $x_0\in\Mc$ such that $F(x_0)\in\KMY$.\bx
\bsk\medskip

\par \textbullet~ {\bf Steiner-type points as a special
case of the finiteness principle for Lipschitz selections.}\smallskip
\par Let $Y$ be a Banach space. Given $m\in\N$ let $\Mc=\Kc_m(Y)$ be the family of all nonempty convex compact subsets of $Y$ of affine dimension at most $m$ equipped with the Hausdorff distance $\rho=\dhf$.
\par Let $F:\Mc\to\KMY$ be the identity mapping on $\KMY$, i.e.,
$$
F(K)=K~~~~\text{for every}~~K\in\KMY.
$$
\par By Theorem \reff{ST-P}, this mapping has a selection $\SX:\Mc\to Y$ whose $\dhf$-Lipschitz seminorm is bounded by a constant $\gamma=\gamma(m)$ depending only on $m$.
\par The following claim asserts that the mapping $F$ satisfies the hypothesis of Theorem \reff{MAIN-FP}.
\begin{claim} \lbl{CL-DH} For every subset $\Mc'\subset\Mc$ with $\#\Mc'\le \FN$ the restriction $F|_{\Mc'}$ has a $\dhf$-Lipschitz selection $f_{\Mc'}:\Mc'\to Y$ with $\|f_{\Mc'}\|_{\Lip((\Mc',\dhf),Y)}\le \theta$ where $\theta=\theta(m)$ is a constant depending only on $m$.
\end{claim}
\par For a simple proof of this claim we refer the reader to \cite[p. 73]{FS-2017}.
\medskip
\par Claim \reff{CL-DH} shows that Theorem \reff{ST-P} can be considered as a particular case of our main result, Theorem \reff{MAIN-FP}, which is applied to the metric space $(\KMY,\dhf)$. (Note that the proof of Theorem \reff{MAIN-FP} uses Theorem \reff{ST-P}.) In general, this metric space has the same complexity as an $L_\infty$-space, and may have infinite Nagata dimension. For example, if $Y=\ell_\infty$ then $\KMY$ contains the set of one point subsets of $\ell_\infty$. As we have noted in the Introduction, $\ell_\infty$ has infinite Nagata dimension so that $\KMY$ has infinite Nagata dimension as well.
\par In this case we are unable to prove Theorem \reff{ST-P} using the ideas and methods developed in Sections 2-4.
\par Thus, analyzing the scheme of the proof of Theorem \reff{MAIN-FP}, we observe that this proof is actually based on solutions of the Lipschitz selection problem for two independent particular cases of this problem, namely, for  metric trees, see Corollary \reff{A-12}, and for the metric space $(\KMY,\dhf)$, see Theorem \reff{ST-P}. Theorem \reff{HDS-M} proven in Section 5 provides a certain ``bridge'' between these two independent results (i.e.,  Corollary \reff{A-12} and Theorem \reff{ST-P}). Combining all these results, we finally obtain a proof of Theorem \reff{MAIN-FP} in the general case.
\bsk\medskip
\par\centerline{\bf Acknowledgments}
\medskip
\par\noindent We are grateful to Alexander  Brudnyi, Arie Israel, Bo'az Klartag, Garving (Kevin) Luli and the participants of the 10th Whitney Problems Conference, Williamsburg, VA, for valuable conversations. We thank the referee for very careful reading and numerous suggestions, which led to improvements in our exposition.
\par We are grateful also to the College of William and Mary, Williamsburg, VA, USA, the American Institute of Mathematics, San Jose, CA, USA, the Fields Institute, Toronto, Canada, the University of Arkansas, AR, USA, the Banff International Research Station, Banff, Canada, the Centre International de Rencontres Math\'ematiques (CIRM), Luminy, Marseille, France, and the Technion, Haifa, Israel, for hosting and supporting workshops on the topic of this paper and closely related problems.
\par Finally, we thank the US-Israel Binational Science Foundation, the US National Science Foundation, the Office
of Naval Research and the Air Force Office of Scientific Research for generous support.
\bigskip

\par\noindent {\sc Charles Fefferman}, Department of Mathematics, Princeton University, Fine Hall Washington Road, Princeton, NJ 08544,  USA
{\hfill {\texttt cf@math.princeton.edu}}

\bigskip
\par\noindent {\sc Pavel Shvartsman}, Department of Mathematics, Technion - Israel Institute of Technology, 32000 Haifa, Israel
{\hfill {\texttt pshv@technion.ac.il}}

\end{document}